%% file: BVector.English.tex
\documentclass{amsbook}
\input{BVector}

\begin{document}
\title{Free Algebra with Countable Basis}

\begin{abstract}
In this book I treat the structure of
$D$\Hyph module which has countable basis.
If we do not care for topology of $D$\Hyph module,
then we consider Hamel basis.
If norm is defined in $D$\Hyph module,
then we consider Schauder basis.
In case of Schauder basis,
we consider vectors
whose expansion in the basis
converges normally.
\end{abstract}

\ShowEq{contents}
\end{document}

%% file: BVector.tex
\def\BookNumber{BVector}
\def\PrintBook{}
\input{Commands}
\xRefDef{0701.238}
\xRefDef{1011.3102}
\xRefDef{8433-5163}
\xRefDef{8443-0072}

\DefEq
{
\maketitle
\tableofcontents
\input{\FilePrefix Preface.\BookNumber.English}
\input{\FilePrefix Convention.English}

\input{\FilePrefix Hamel.Basis.English}
\input{\FilePrefix Schauder.Basis.English}

\input{\FilePrefix Biblio.English}
\input{\FilePrefix Index.English}
\input{\FilePrefix Symbol.English}
}
{contents}

%% file: Commands.tex
\def\Defined{}
\scrollmode
\ifx\FilePrefix\undefined
\newcommand{\FilePrefix}{}
\fi
\ifx\UseRussian\Defined
	\usepackage{cmap}
		\usepackage[T2A,T2B]{fontenc}
	\usepackage[cp1251]{inputenc}
	\usepackage[english,russian]{babel}
\fi
\ifx\GJSFRA\Defined
\usepackage[top=1.905cm,bottom=1.905cm,inner=1.65cm,outer=1.65cm]{geometry}
\fi
\ifx\Presentation\Defined
	\usepackage[margin=1cm]{geometry}
\fi
\ifx\CreateSpace\Defined
	\usepackage[top=1.905cm,bottom=1.905cm,inner=1.905cm,outer=1.27cm]{geometry}
\fi
\ifx\PublishBook\Defined
	\def\PrintPaper{}
	\usepackage{setspace}
	\ifx\UseRussian\undefined
		\usepackage{pslatex}
	\fi
	\usepackage[margin=2cm]{geometry}
\fi
\usepackage{footmisc}
\usepackage[all]{xy}
\usepackage{color}
\ifx\PrintPaper\undefined
	\definecolor{UrlColor}{rgb}{.9,0,.3}
	\definecolor{SymbColor}{rgb}{.4,0,.9}
	\definecolor{IndexColor}{rgb}{1,.3,.6}
	\newcommand\BlueText[1]{\textcolor{blue}{#1}}
	
\else
	\definecolor{UrlColor}{rgb}{.1,.1,.1}
	\definecolor{SymbColor}{rgb}{.1,.1,.1}
	\definecolor{IndexColor}{rgb}{.1,.1,.1}
	\newcommand\BlueText[1]{#1}
	
\fi

\usepackage{xr-hyper}
	\usepackage[unicode]{hyperref}
	\hypersetup{pdfdisplaydoctitle=true}
	\hypersetup{colorlinks}
	\hypersetup{citecolor=UrlColor}
	\hypersetup{urlcolor=UrlColor}
	\hypersetup{linkcolor=UrlColor}
	\hypersetup{pdffitwindow=true}
	\hypersetup{pdfnewwindow=true}
	\hypersetup{pdfstartview={FitH}}

\newcounter{Index}
\newcounter{Symbol}

\def\hyph{\penalty0\hskip0pt\relax-\penalty0\hskip0pt\relax}
\def\Hyph{-\penalty0\hskip0pt\relax}

\def\ValueOff{off}
\def\ValueOn{on}
\def\Items#1{\ItemList#1,LastItem,}%
\def\LastItem{LastItem}%
\def\ItemList#1,{\def\ViewBook{#1}%
	\ifx\ViewBook\LastItem%
	\else%
		\ifx\ViewBook\BookNumber%
			\def\Semafor{on}%
		\fi%
		\expandafter\ItemList%
	\fi%
}%

\newcommand{\ePrints}[1]
{%
	\def\Semafor{off}%
	\Items{#1}%
}%

\newcommand{\Basis}[1]{\overline{\overline{#1}}{}}
\newcommand{\Vector}[1]{\overline{#1}{}}
\ifx\PrintPaper\undefined
	\newcommand{\gi}[1]{\boldsymbol{\textcolor{IndexColor}{#1}}}
\else
	\newcommand{\gi}[1]{\boldsymbol{#1}}
\fi

\newcommand{\VX}[1]{\Vector{#1}_{[1]}}
\makeatletter
\newcommand{\NameDef}[1]{%
	\expandafter\gdef\csname #1\endcsname%
}%
\newcommand{\xNameDef}[1]{%
	\expandafter\xdef\csname #1\endcsname%
}%
\newcommand{\ShowSymbol}[1]{%
	\@nameuse{ViewSymbol#1}%
}%
\newcommand{\symb}[3]{%
	\@ifundefined{ViewSymbol#3}{%
		\NameDef{ViewSymbol#3}{\textcolor{SymbColor}{#1}}%
		\xNameDef{RefSymbol#3}{:}%
		\@namedef{LabeSymbol}{\label{symbol: #3::}}%
	}{%
		\addtocounter{Symbol}{1}%
		\xNameDef{RefSymbol#3}{\@nameuse{RefSymbol#3},\arabic{Symbol}}%
		\@namedef{LabeSymbol}{\label{symbol: #3:\arabic{Symbol}}}%
	}%
	\ifcase#2
	\or
		$\@nameuse{ViewSymbol#3}$%
	\or
		\[\@nameuse{ViewSymbol#3}\]%
	\else%
	\fi%
	\@nameuse{LabeSymbol}%
}%
\newcommand{\DefEq}[2]{%
	\@ifundefined{ViewEq#2}{%
		\NameDef{ViewEq#2}{#1}%
	}{%
	}%
}%
\newcommand{\DefEquation}[2]{%
	\DefEq
	{
	\begin{equation}
	#1
	\EqLabel{#2}
	\end{equation}
	}
	{#2}
}%
\newcommand\EqRef[1]{\eqref{eq: #1}}%
\newcommand\eqRef[2]{\EqRef{#1, #2}}%
\newcommand\EqLabel[1]{\label{eq: #1}}%
\newcommand\ShowEq[1]{%
	\@ifundefined{ViewEq#1}{%
		\message {error: missed ShowEq #1}
  }{%
	\@nameuse{ViewEq#1}%
	}%
}%
\newcommand\DrawEq[2]{%
	\@ifundefined{ViewEq#1}{%
		\message {error: missed ShowEq #1}%
  }{%
	\def\Temp{}%
	\def\Tempa{#2}%
	\ifx\Tempa\Temp%
		\[%
		\@nameuse{ViewEq#1}%
		\]%
	\else%
		\begin{equation}%
		\@nameuse{ViewEq#1}%
		\EqLabel{#1, #2}%
		\end{equation}%
	\fi%
	}%
}%
\makeatother

\newcommand{\subs}{${}_*$\Hyph}
\newcommand{\sups}{${}^*$\Hyph}

\newcommand{\CRstar}{{}^*{}_*}
\newcommand{\RCstar}{{}_*{}^*}
\newcommand{\CRcirc}{{}^{\circ}{}_{\circ}}
\newcommand{\RCcirc}{{}_{\circ}{}^{\circ}}

\newcommand{\RC}{$\RCstar$\Hyph}
\newcommand{\CR}{$\CRstar$\Hyph}
\newcommand{\drc}{$D\RCstar$\Hyph}
\newcommand{\Drc}{$\mathcal D\RCstar$\Hyph}
\newcommand{\dcr}{$D\CRstar$\hyph}
\newcommand{\rcd}{$\RCstar D$\Hyph}
\newcommand{\crd}{$\CRstar D$\Hyph}

\newcommand{\Acr}{$A\CRcirc$\Hyph}

\newcommand\sT{$\star T$\Hyph}%
\newcommand\Ts{$T\star$\Hyph}%
\newcommand\sD{$\star D$\Hyph}%
\newcommand\Ds{$D\star$\Hyph}%

\newcommand\pC[2]{{}_{#1\cdot #2}}%
\newcommand\DrcPartial[1]%
{%
	\def\tempa{}%
	\def\tempb{#1}%
	\ifx\tempa\tempc%
		(\partial\RCstar)%
	\else%
		({}_{#1}\partial\RCstar)%
	\fi%
}%
\newcommand\crDPartial[1]%
{%
	\def\tempa{}%
	\def\tempb{#1}%
	\ifx\tempa\tempc%
		(\CRstar\partial)%
	\else%
		(\CRstar\partial_{#1})%
	\fi%
}%
\newcommand\StandPartial[3]%
{%
	\frac{\partial^{\gi{#3}} #1}{\partial #2}%
}%

	\renewcommand{\uppercasenonmath}[1]{}

\makeatletter
\newcommand\@dotsep{4.5}
\def\@tocline#1#2#3#4#5#6#7
{\relax
		\par \addpenalty\@secpenalty\addvspace{#2}%
		\begingroup 
		\@ifempty{#4}{%
			\@tempdima\csname r@tocindent\number#1\endcsname\relax
		}{%
			\@tempdima#4\relax
		}%
		\parindent\z@ \leftskip#3\relax \advance\leftskip\@tempdima\relax
		\rightskip\@pnumwidth plus1em \parfillskip-\@pnumwidth
		#5\leavevmode\hskip-\@tempdima #6\relax
		\leaders\hbox{$\m@th
		\mkern \@dotsep mu\hbox{.}\mkern \@dotsep mu$}\hfill
		\hbox to\@pnumwidth{\@tocpagenum{#7}}\par
		\nobreak
		\endgroup
}
\makeatother 

\ifx\PrintBook\undefined

	\makeatletter
		\renewcommand{\@indextitlestyle}{%
			\twocolumn[\section{\indexname}]%
			\def\IndexSpace{off}%
		}
	\makeatother 
		\ifx\PrintPaper\undefined
			\thanks{\href{mailto:Aleks\_Kleyn@MailAPS.org}{Aleks\_Kleyn@MailAPS.org}}
			\ePrints{1102.1776,1201.4158}
			\ifx\Semafor\ValueOff
					\thanks{\ \ \ \url{http://sites.google.com/site/AleksKleyn/}}
					\thanks{\ \ \ \url{http://arxiv.org/a/kleyn\_a\_1}}
					\thanks{\ \ \ \url{http://AleksKleyn.blogspot.com/}}
			\fi
		\fi
\else

	\makeatletter
\def\@maketitle{%
  \cleardoublepage \thispagestyle{empty}%
  \begingroup \topskip\z@skip
  \null\vfil
  \begingroup
  \LARGE\bfseries \centering
  \openup\medskipamount
  \@title
  \par
  \ifx\subtitle\undefined
  \else
  \centerline{\emph\subtitle}
  \fi
  \par\vspace{24pt}%
  \def\and{\par\medskip}\centering
  \mdseries\authors\par\bigskip
  \endgroup
  \vfil
  \ifx\@empty\addresses \else \@setaddresses \fi
  \vfil
  \ifx\@empty\@dedicatory
  \else \begingroup
    \centering{\footnotesize\itshape\@dedicatory\@@par}%
    \endgroup
  \fi
  \vfill
  \newpage\thispagestyle{empty}
  \begin{center}
    \ifx\@empty\@subjclass\else\@setsubjclass\fi
    \ifx\@empty\@keywords\else\@setkeywords\fi
    \ifx\@empty\@translators\else\vfil\@settranslators\fi
    \ifx\@empty\thankses\else\vfil\@setthanks\fi
  \end{center}
  \vfil
  \@setabstract
\vfil
  \def\Temp{0000}
  \ifx\copyrightyear\Temp
  \else
  \begin{center}
\begin{tabular}{@{}c}
Copyright\ \copyright\ \copyrightyear\ \copyrightholder
\\
All rights reserved.
\end{tabular}
  \end{center}
  \fi
  \ifx\ISBN\undefined%
  \else%
  \begin{center}
\begin{tabular}{@{}r@{\ }l}
ISBN:&\ISBN
\\
ISBN-13:&\ISBNa
\end{tabular}
  \end{center}
  \fi%
  \endgroup}
	\renewcommand{\@indextitlestyle}{%
		\twocolumn[\chapter{\indexname}]%
		\def\IndexSpace{off}%
		\let\@secnumber\@empty
		\chaptermark{\indexname}%
	}
	\makeatother 
	\email{\href{mailto:Aleks\_Kleyn@MailAPS.org}{Aleks\_Kleyn@MailAPS.org}}
	\ePrints{1102.1776,1201.4158}
	\ifx\Semafor\ValueOff
		\urladdr{\url{http://sites.google.com/site/alekskleyn/}}
		\urladdr{\url{http://arxiv.org/a/kleyn\_a\_1}}
		\urladdr{\url{http://AleksKleyn.blogspot.com/}}
	\fi
\fi

\ifx\SelectlEnglish\undefined
	\ifx\UseRussian\undefined
		\def\SelectlEnglish{}
	\fi
\fi

\newcommand\arXivRef{http://arxiv.org/PS_cache/}
\newcommand\wRefDef[2]
	{
	\def\Tempa{#1}
	\def\Tempb{0405.027}
	\ifx\Tempa\Tempb
	\def\wRef{\arXivRef gr-qc/pdf/0405/0405027v3.pdf}
	\fi
	\def\Tempb{0405.028}
	\ifx\Tempa\Tempb
	\def\wRef{\arXivRef gr-qc/pdf/0405/0405028v5.pdf}
	\fi
	\def\Tempb{0412.391}
	\ifx\Tempa\Tempb
	\def\wRef{\arXivRef math/pdf/0412/0412391v4.pdf}
	\fi
	\def\Tempb{0612.111}
	\ifx\Tempa\Tempb
	\def\wRef{\arXivRef math/pdf/0612/0612111v2.pdf}
	\fi
	\def\Tempb{0701.238}
	\ifx\Tempa\Tempb
	\def\wRef{\arXivRef math/pdf/0701/0701238v4.pdf}
	\fi
	\def\Tempb{0702.561}
	\ifx\Tempa\Tempb
	\def\wRef{\arXivRef math/pdf/0702/0702561v3.pdf}
	\fi
	\def\Tempb{0707.2246}
	\ifx\Tempa\Tempb
	\def\wRef{\arXivRef arxiv/pdf/0707/0707.2246v2.pdf}
	\fi
	\def\Tempb{0803.3276}
	\ifx\Tempa\Tempb
	\def\wRef{\arXivRef arxiv/pdf/0803/0803.3276v3.pdf}
	\fi
	\def\Tempb{0812.4763}
	\ifx\Tempa\Tempb
	\def\wRef{\arXivRef arxiv/pdf/0812/0812.4763v6.pdf}
	\fi
 	\def\Tempb{0906.0135}
	\ifx\Tempa\Tempb
	\def\wRef{\arXivRef arxiv/pdf/0906/0906.0135v3.pdf}
	\fi
 	\def\Tempb{0909.0855}
	\ifx\Tempa\Tempb
	\def\wRef{\arXivRef arxiv/pdf/0909/0909.0855v5.pdf}
	\fi
 	\def\Tempb{0912.3315}
	\ifx\Tempa\Tempb
	\def\wRef{\arXivRef arxiv/pdf/0912/0912.3315v2.pdf}
	\fi
 	\def\Tempb{0912.4061}
	\ifx\Tempa\Tempb
	\def\wRef{\arXivRef arxiv/pdf/0912/0912.4061v2.pdf}
	\fi
 	\def\Tempb{1003.1544}
	\ifx\Tempa\Tempb
	\def\wRef{\arXivRef arxiv/pdf/1003/1003.1544v2.pdf}
	\fi
 	\def\Tempb{1104.5197}
	\ifx\Tempa\Tempb
	\def\wRef{\arXivRef arxiv/pdf/1104/1104.5197.pdf}
	\fi
 	\def\Tempb{1011.3102}
	\ifx\Tempa\Tempb
	\def\wRef{\arXivRef arxiv/pdf/1011/1011.3102.pdf}
	\fi
 	\def\Tempb{1105.4307}
	\ifx\Tempa\Tempb
	\def\wRef{\arXivRef arxiv/pdf/1105/1105.4307.pdf}
	\fi
 	\def\Tempb{1107.1139}
	\ifx\Tempa\Tempb
	\def\wRef{\arXivRef arxiv/pdf/1104/1107.1139.pdf}
	\fi
 	\def\Tempb{1107.5037}
	\ifx\Tempa\Tempb
	\def\wRef{\arXivRef arxiv/pdf/1107/1107.5037.pdf}
	\fi
 	\def\Tempb{1202.6021}
	\ifx\Tempa\Tempb
	\def\wRef{\arXivRef arxiv/pdf/1202/1202.6021.pdf}
	\fi
 	\def\Tempb{8433-5163}
	\ifx\Tempa\Tempb
	\def\wRef{http://www.amazon.com/}
	\fi
 	\def\Tempb{8443-0072}
	\ifx\Tempa\Tempb
	\def\wRef{http://www.amazon.com/}
	\fi
	\externaldocument[#1-#2-]{\FilePrefix #1.#2}[\wRef]
	}
\newcommand\LanguagePrefix{}%
\makeatletter
\newcommand\StartLabelItem[1]
{
\def\theenumi{\ref*{#1}.\@arabic\c@enumi}
\def\labelenumi{\theenumi:}
}
\newcommand\StopLabelItem
{
\def\theenumi{\@arabic\c@enumi}
\def\labelenumi{(\theenumi)}
}
\makeatother
\newcommand\RefItem[1]{\ref{item: #1}}
\ifx\SelectlEnglish\undefined
	\newcommand\CurrentLanguage{Russian.}%
	\author{Александр Клейн}
		\newtheorem{theorem}{Теорема}[section]
		\newtheorem{corollary}[theorem]{Следствие}
		\newtheorem{convention}[theorem]{Соглашение}
		\theoremstyle{definition}
		\newtheorem{definition}[theorem]{Определение}
		\newtheorem{example}[theorem]{Пример}
		\newtheorem{xca}[theorem]{Exercise}
		\theoremstyle{remark}
		\newtheorem{remark}[theorem]{Замечание}
		\newtheorem{lemma}[theorem]{Лемма}
		
	\newcommand\xRefDef[1]
		{
		\wRefDef{#1}{Russian}
		\NameDef{xRefDef#1}{}%
		}
	\makeatletter
	\newcommand\xRef[2]%
	{%
		\@ifundefined{xRefDef#1}{%
		\ref{#2}%
		}{%
		\citeBib{#1}-\ref{#1-Russian-#2}%
		}%
	}%
	\newcommand\xEqRef[2]%
	{%
		\@ifundefined{xRefDef#1}{%
		\eqref{eq: #2}%
		}{%
		\citeBib{#1}-\eqref{#1-Russian-eq: #2}%
		}%
	}%
	\makeatother
	\ifx\PrintBook\undefined
		\newcommand{\BibTitle}{%
			\section{Список литературы}%
		}
	\else
		\newcommand{\BibTitle}{%
			\chapter{Список литературы}%
		}
	\fi
\else
	\newcommand\CurrentLanguage{English.}%
	\author{Aleks Kleyn}
		\newtheorem{theorem}{Theorem}[section]
		\newtheorem{corollary}[theorem]{Corollary}
		\newtheorem{convention}[theorem]{Convention}
		\theoremstyle{definition}
		\newtheorem{definition}[theorem]{Definition}

		\theoremstyle{remark}
		\newtheorem{remark}[theorem]{Remark}
		\newtheorem{lemma}[theorem]{Lemma}
	\newcommand\xRefDef[1]
		{
		\wRefDef{#1}{English}
		\NameDef{xRefDef#1}{}%
		}
	\makeatletter
	\newcommand\xRef[2]%
	{%
		\@ifundefined{xRefDef#1}{%
		\ref{#2}%
		}{%
		\citeBib{#1}-\ref{#1-English-#2}%
		}%
	}%
	\newcommand\xEqRef[2]%
	{%
		\@ifundefined{xRefDef#1}{%
		\eqref{eq: #2}%
		}{%
		\citeBib{#1}-\eqref{#1-English-eq: #2}%
		}%
	}%
	\makeatother
	\ifx\PrintBook\undefined
		\newcommand{\BibTitle}{%
			\section{References}%
		}
	\else
		\newcommand{\BibTitle}{%
			\chapter{References}%
		}
	\fi
\fi

\ifx\PrintBook\undefined
		\numberwithin{Hfootnote}{section}
\else
	\numberwithin{section}{chapter}
	\numberwithin{footnote}{chapter}
	\numberwithin{Hfootnote}{chapter}
\fi

\ifx\Presentation\undefined
	\numberwithin{equation}{section}
	\numberwithin{figure}{section}
	\numberwithin{table}{section}
	\numberwithin{Item}{section}
\fi

\makeatletter
\newcommand\org@maketitle{}
\let\org@maketitle\maketitle
\def\maketitle{%
	\hypersetup{pdftitle={\@title}}%
	\hypersetup{pdfauthor={\authors}}%
	\hypersetup{pdfsubject=\@keywords}%
	\ifx\UseRussian\Defined
		\pdfbookmark[1]{\@title}{TitleRussian}
	\else
		\pdfbookmark[1]{\@title}{TitleEnglish}
	\fi
	\org@maketitle
}
\def\make@stripped@name#1{%
	\begingroup
		\escapechar\m@ne
		\global\let\newname\@empty
		\protected@edef\Hy@tempa{\CurrentLanguage #1}%
		\edef\@tempb{%
			\noexpand\@tfor\noexpand\Hy@tempa:=%
			\expandafter\strip@prefix\meaning\Hy@tempa
		}%
		\@tempb\do{%
			\if\Hy@tempa\else
				\if\Hy@tempa\else
					\xdef\newname{\newname\Hy@tempa}%
				\fi
			\fi
		}%
	\endgroup
}%
\newenvironment{enumBib}{%
	\BibTitle
	\advance\@enumdepth \@ne
	\edef\@enumctr{enum\romannumeral\the\@enumdepth}\list
	{\csname biblabel\@enumctr\endcsname}{\usecounter
	{\@enumctr}\def\makelabel##1{\hss\llap{\upshape##1}}}
}{%
	\endlist
}

\makeatletter

\newcommand{\BiblioItem}[2]
{
	\def\Semafor{off}
	\@ifundefined{\LanguagePrefix ViewCite#1}{}{%
		\def\Semafor{on}%
	}%
	\ifx\Semafor\ValueOff
		\@ifundefined{xRefDef#1}{}{%
		\def\Semafor{on}%
		}%
	\fi
	\ifx\Semafor\ValueOn
		\ifx\IndexState\ValueOff
			\begin{enumBib}
			\def\IndexState{on}
		\fi
		\item \label{\LanguagePrefix bibitem: #1}#2%
	\fi
}
\makeatother
\newcommand{\OpenBiblio}
{
	\def\IndexState{off}
}
\newcommand{\CloseBiblio}
{
	\ifx\IndexState\ValueOn
		\end{enumBib}
		\def\IndexState{off}
	\fi
}

\makeatletter
\def\StartCite{[}%
\def\citeBib#1{\protect\showCiteBib#1,endCite,}%
\def\endCite{endCite}%
\def\showCiteBib#1,{\def\temp{#1}%
\ifx\temp\endCite
]%
\def\StartCite{[}%
\else
	\StartCite\LanguagePrefix \ref{\LanguagePrefix bibitem: #1}%
	\@ifundefined{\LanguagePrefix ViewCite#1}{%
		\NameDef{\LanguagePrefix ViewCite#1}{}%
	}{%
	}%
	\def\StartCite{, }%
\expandafter\showCiteBib%
\fi}%
\makeatother

\newcommand{\arp}{\ar @{-->}}
\newcommand{\ars}{\ar @{.>}}
\newcommand\Bundle[1]{{\mathbb #1}}
\newcommand{\bundle}[4]%
{%
	\def\tempa{}%
	\def\tempb{#3}%
	\def\tempc{#1}%
	\ifx\tempa\tempb%
		\ifx\tempa\tempc%
			#2%
		\else%
			\xymatrix{#2:#1\arp[r]&#4}%
		\fi%
	\else%
		\ifx\tempa\tempc%
			#2[#3]%
		\else%
			\xymatrix{#2[#3]:#1\arp[r]&#4}%
		\fi%
	\fi%
}%
\makeatletter
\newcommand{\AddIndex}[2]%
{%
	\@ifundefined{RefIndex#2}{%
		\xNameDef{RefIndex#2}{:}%
		\@namedef{LabelIndex}{\label{index: #2::}}%
	}{%
		\addtocounter{Index}{1}%
		\xNameDef{RefIndex#2}{\@nameuse{RefIndex#2},\arabic{Index}}%
		\@namedef{LabelIndex}{\label{index: #2:\arabic{Index}}}%
	}%
	\@nameuse{LabelIndex}%
	{\bf #1}%
}%
\newcommand{\Index}[2]%
{%
	\def\Semafor{off}%
	\@ifundefined{RefIndex#2}{%
	}{%
		\def\Semafor{on}
	}%
	\ifx\Semafor\ValueOn%
		\def\tempa{}%
		\def\tempb{#2}%
		\ifx\IndexState\ValueOff%
			\begin{theindex}%
			\def\IndexState{on}%
		\fi%
		\ifx\IndexSpace\ValueOn%
			\indexspace%
			\def\IndexSpace{off}%
		\fi%
		\item #1%
		\ifx\tempa\tempb%
		\else%
			\edef\PageRefs{\@nameuse{RefIndex#2}}
			\def\Sep{\ }%
			\@for\PageRef:=\PageRefs\do{%
			\Sep
			\pageref{index: #2:\PageRef}%
			\def\Sep{,\ }%
			}%
		\fi%
	\fi%
}%

\newcommand{\Symb}[2]
{
	\def\Semafor{off}
	\@ifundefined{ViewSymbol#2}{%
	}{%
		\def\Semafor{on}
	}%
	\ifx\Semafor\ValueOn
		\ifx\IndexState\ValueOff
			\begin{theindex}
			\def\IndexState{on}
		\fi
		\ifx\IndexSpace\ValueOn
			\indexspace
			\def\IndexSpace{off}
		\fi
		\item $\displaystyle\@nameuse{ViewSymbol#2}$\ \ #1
		\edef\PageRefs{\@nameuse{RefSymbol#2}}
		\def\Sep{}%
		\@for\PageRef:=\PageRefs\do{%
		\Sep
		\pageref{symbol: #2:\PageRef}%
		\def\Sep{,\ }%
		}%
	\fi
}

\makeatother

\newcommand{\SetIndexSpace}%
{%
	\def\IndexSpace{on}%
}%

\newcommand{\OpenIndex}
{
	\def\IndexState{off}
}
\newcommand{\CloseIndex}
{
	\ifx\IndexState\ValueOn
		\end{theindex}
		\def\IndexState{off}
	\fi
}

\def\LastMemo{LastMemo}%
\def\MemoList#1//{\def\temp{#1}%
	\ifx\temp\LastMemo
	\else%
	\begin{flushright}
\begin{minipage}{200pt}
\setlength{\parindent}{5mm}
		\BlueText{#1}%
\end{minipage}
	\end{flushright}
		\expandafter\MemoList%
	\fi%
}%

%% file: Convention.English.tex
\input{\FilePrefix Convention.Eq}

\section{Conventions}

\ePrints{0812.4763,0906.0135,0908.3307,0909.0855,0912.3315,1003.1544}
\Items{1006.2597}
\ifx\Semafor\ValueOn
\begin{convention}
Function and mapping are synonyms. However according to
tradition, correspondence between either rings or vector
spaces is called mapping and a mapping of
either real field or quaternion algebra is called function.
\qed
\end{convention}
\fi

\ePrints{0701.238,0812.4763,0908.3307,0912.4061,1001.4852}
\Items{1003.1544,4776-3181}
\ifx\Semafor\ValueOn
\begin{convention}
In any expression where we use index I assume
that this index may have internal structure.
For instance, considering the algebra $A$ we enumerate coordinates of
$a\in A$ relative to basis $\Basis e$ by an index $i$.
This means that $a$ is a vector. However, if $a$
is matrix, then we need two indexes, one enumerates
rows, another enumerates columns. In the case, when index has
structure, we begin the index from symbol $\cdot$ in
the corresponding position. 
For instance, if I consider the matrix $a^i_j$ as an element of a vector
space, then I can write the element of matrix as $a^{\cdot}{}^i_j$.
\qed
\end{convention}
\fi

\ePrints{0701.238,0812.4763,0908.3307,0912.4061,1006.2597,1011.3102}
\Items{Calculus.Paper,BVector}
\ifx\Semafor\ValueOn
\begin{convention}
I assume sum over index $s$ in expression like
\ShowEq{Sum over repeated index}
\qed
\end{convention}
\fi

\ePrints{0701.238,0812.4763,0906.0135,0908.3307,0909.0855}
\ifx\Semafor\ValueOn
\begin{convention}
We can consider division ring $D$ as $D$\Hyph vector space
of dimension $1$. According to this statement, we can explore not only
homomorphisms of division ring $D_1$ into division ring $D_2$,
but also linear maps of division rings.
This means that map is multiplicative over
maximum possible field. In particular, linear map
of division ring $D$ is multiplicative over center $Z(D)$. This statement
does not contradict with
definition of linear map of field because for field $F$ is true
$Z(F)=F$.
When field $F$ is different from
maximum possible, I explicit tell about this in text.
\qed
\end{convention}
\fi

\ePrints{0912.4061}
\ifx\Semafor\ValueOn
\begin{convention}
For given field $F$, unless otherwise stated,
we consider finite dimensional $F$\Hyph algebra.
\qed
\end{convention}
\fi

\ePrints{0701.238,0812.4763,0906.0135,0908.3307,4776-3181}
\ifx\Semafor\ValueOn
\begin{convention}
In spite of noncommutativity of product a lot of statements
remain to be true if we substitute, for instance, right representation by
left representation or right vector space by left
vector space.
To keep this symmetry in statements of theorems
I use symmetric notation.
For instance, I consider \Ds vector space
and \sD vector space.
We can read notation \Ds vector space
as either D\Hyph star\Hyph vector space or
left vector space.
We can read notation \Ds linear dependent vectors
as either D\Hyph star\Hyph linear dependent vectors or
vectors that are linearly dependent from left.
\qed
\end{convention}
\fi

\ePrints{0701.238,0812.4763,0906.0135,0908.3307,0909.0855,0912.4061}
\Items{1001.4852,1003.1544,1006.2597,1104.5197,1105.4307,1107.1139}
\Items{1202.6021,MQuater,BVector,CACAA.07.195}
\ifx\Semafor\ValueOn
\begin{convention}
Let $A$ be free finite
dimensional algebra.
Considering expansion of element of algebra $A$ relative basis $\Basis e$
we use the same root letter to denote this element and its coordinates.
However we do not use vector notation in algebra.
In expression $a^2$, it is not clear whether this is component
of expansion of element
$a$ relative basis, or this is operation $a^2=aa$.
To make text clearer we use separate color for index of element
of algebra. For instance,
\ShowEq{Expansion relative basis in algebra}
\qed
\end{convention}
\fi

\ePrints{0701.238,0812.4763,0906.0135,0908.3307,0909.0855,0912.4061}
\Items{1001.4852,1003.1544,1006.2597,1104.5197,1105.4307,1107.1139}
\Items{1202.6021,MQuater,CACAA.07.195}
\ifx\Semafor\ValueOn
\begin{convention}
If free finite dimensional algebra has unit, then we identify
the vector of basis $\Vector e_{\gi 0}$ with unit of algebra.
\qed
\end{convention}
\fi

\ePrints{1104.5197,1105.4307}
\ifx\Semafor\ValueOn
\begin{convention}
Although the algebra is a free module over some
ring, we do not use the vector notation
to write elements of algebra. In the case when I consider the
matrix of coordinates of element of algebra, I will use vector
notation to write corresponding element.
In order to avoid ambiguity when I use conjugation,
I denote $a^*$ element conjugated to element $a$.
\qed
\end{convention}
\fi

\ePrints{0906.0135,0912.3315,8443-0072,1111.6035,1102.5168}
\ifx\Semafor\ValueOn
\begin{convention}
In \citeBib{Cohn: Universal Algebra},
an arbitrary operation of algebra is denoted by letter $\omega$,
and $\Omega$ is the set of operations of some universal algebra.
Correspondingly, the universal algebra with the set of operations
$\Omega$ is denoted as $\Omega$\Hyph algebra.
Similar notations we see in
\citeBib{Burris Sankappanavar} with small difference
that an operation in the algebra is denoted by letter $f$
and $\mathcal F$ is the set of operations.
I preferred first case of notations because in this case it is
easier to see where I use operation.
\qed
\end{convention}
\fi

\ePrints{0906.0135,0912.3315,8443-0072}
\ifx\Semafor\ValueOn
\begin{convention}
Since the number of universal algebras
in the tower of representations is varying,
then we use vector notation for a tower of
representations. We denote the set
$(A_1,...,A_n)$ of $\Omega_i$\Hyph algebras $A_i$, $i=1$, ..., $n$
as $\Vector A$. We denote the set of representations
$(f_{1,2},...,f_{n-1,n})$ of these algebras as $\Vector f$.
Since different algebras have different type, we also
talk about the set of $\Vector{\Omega}$\Hyph algebras.
\ePrints{8443-0072}
\ifx\Semafor\ValueOn
We
\else
In relation to the set $\Vector A$,
we also use matrix notations 
that we discussed
in section \xRef{0701.238}{section: Concept of Generalized Index}.
For instance, we
\fi
use the symbol $\Vector A_{[1]}$ to denote the
set of $\Vector{\Omega}$\Hyph algebras $(A_2,...,A_n)$.
In the corresponding notation $(\VX A,\Vector f)$ of tower
of representation, we assume that $\Vector f=(f_{2,3},...,f_{n-1,n})$.
\qed
\end{convention}

\begin{convention}
Since we use vector notation for elements of the
tower of representations, we need convention about notation of operation.
We assume that we get result of operation componentwise. For instance,
\ShowEq{vector notation in tower of representations}
\qed
\end{convention}
\fi

\ePrints{8443-0072,1111.6035,0906.0135,NewAffine,1102.5168}
\ifx\Semafor\ValueOn
\begin{convention}
Let $A$ be $\Omega_1$\Hyph algebra.
Let $B$ be $\Omega_2$\Hyph algebra.
Notation
\ShowEq{A->*B}
means that there is representation of $\Omega_1$\Hyph algebra $A$
in $\Omega_2$\Hyph algebra $B$.
\qed
\end{convention}
\fi

\ePrints{0702.561,0707.2246,0803.2620}
\ifx\Semafor\ValueOn
\begin{convention}
I use arrow $\xymatrix{\arp[r]&}$ to represent
projection of bundle on diagram.
I use arrow $\xymatrix{\ars[r]&}$ to represent
section of bundle on diagram.
\qed
\end{convention}
\fi

\ePrints{0912.3315}
\ifx\Semafor\ValueOn
\begin{remark}
I believe that diagrams of maps are an important tool.
However, sometimes I want
to see the diagram as three dimensional figure
and I expect that this would increase its expressive
power. Who knows what surprises the future holds.
In 1992, at a conference in Kazan, I have described to my colleagues
what advantages the computer preparation of papers has.
8 years later I learned from the letter from Kazan that now we can
prepare paper using LaTeX.
\qed
\end{remark}
\fi

\ePrints{1001.4852,1003.1544,1006.2597,1011.3102}
\Items{Calculus.Paper}
\ifx\Semafor\ValueOn
\begin{convention}
If, in a certain expression, we use several operations
which include the operation $\circ$, then
it is assumed that the operation $\circ$ is executed first.
Below is an example of equivalent expressions.
\ShowEq{list circ expressions}
\qed
\end{convention}
\fi


\ePrints{1107.1139}
\ifx\Semafor\ValueOn
\begin{convention}
For given $D$\Hyph algebra $A$
we define left shift
\ShowEq{left shift, D algebra}
by the equation
\ShowEq{left shift 1, D algebra}
and right shift
\ShowEq{right shift, D algebra}
by the equation
\ShowEq{right shift 1, D algebra}
\qed
\end{convention}
\fi

\ifx\PrintPaper\undefined
Without a doubt, the reader may have questions,
comments, objections. I will appreciate any response.
\fi

%% file: Convention.Eq.tex

\DefEq
{
\[
\Vector r(\Vector a)=(r_1(a_1),...,r_n(a_n))
\]
}
{vector notation in tower of representations}

\DefEq
{
\[
\xymatrix
{
A\ar[r]|{*}&B
}
\]
}
{A->*B}

\DefEq
{
\[
\begin{array}{r@{\ }lr@{\ }l}
f\circ xy&\equiv f(x)y
&
f\circ(xy)&\equiv f(xy)
\\
f\circ x+y&\equiv f(x)+y
&
f\circ (x+y)&\equiv f(x+y)
\end{array}
\]
}
{list circ expressions}

\DefEq
{
\[
a\pC s0xa\pC s1
\]
}
{Sum over repeated index}

\DefEq
{
\[
a=a^{\gi i}\Vector e_{\gi i}
\]
}
{Expansion relative basis in algebra}

\DefEq
{
\symb{a\circ}1{left shift, D algebra}
}
{left shift, D algebra}

\DefEq
{
\symb{a\star}1{right shift, D algebra}
}
{right shift, D algebra}

\DefEq
{
\[
\ShowSymbol{left shift, D algebra}x=ax
\]
}
{left shift 1, D algebra}

\DefEq
{
\[
\ShowSymbol{right shift, D algebra}x=xa
\]
}
{right shift 1, D algebra}

%% file: Hamel.Basis.English.tex
\input{Hamel.Basis.Eq}
\ifx\PrintBook\Defined
\chapter{Hamel Basis}
\fi

\section{Module}
\label{Section: Module}

\begin{theorem}
\label{theorem: effective representation of the ring}
Let ring $D$ has unit $e$.
Representation\footnote{This subsection
is written on the base of the section
\xRef{8433-5163}{Section: Module}.}
\ShowEq{representation of the ring}
of the ring $D$
in an Abelian group $A$ is
\AddIndex{effective}{effective representation of ring}
iff $a=0$
follows from equation $f(a)=0$.
\end{theorem}
\begin{proof}
We define the sum of transformations $f$ and $g$ of an Abelian group
according to rule
\ShowEq{sum of transformations of Abelian group}
Therefore, considering the representation of the ring $D$ in
the Abelian group $A$, we assume
\ShowEq{sum of transformations of Abelian group, 1}
We define the product of transformation of representation
according to rule
\ShowEq{product of transformations of representation}

Suppose $a$, $b\in R$
cause the same transformation. Then
\ShowEq{representation of ring, 1}
for any $m\in A$.
From the equation
\EqRef{representation of ring, 1}
it follows that $a-b$ generates zero transformation
\ShowEq{representation of ring, 2}
Element $e+a-b$ generates an identity transformation.
Therefore, the representation $f$ is effective
iff $a=b$.
\end{proof}

\begin{definition}
\label{definition: module over ring}
Let $D$ be commutative ring.
Effective representation of ring $D$
in an Abelian group $A$ is called
Abelian group $A$ is called either 
\AddIndex{module over ring}{module over ring} $D$
or
\AddIndex{$D$\Hyph module}{D module}.
\qed
\end{definition}

\begin{theorem}
\label{theorem: definition of module}
Following conditions hold for $D$\Hyph module $A$:
\begin{itemize}
\item 
\AddIndex{associative law}{associative law, D module}
\ShowEq{associative law, D module}
\item 
\AddIndex{distributive law}{distributive law, D module}
\ShowEq{distributive law, D module}
\item 
\AddIndex{unitarity law}{unitarity law, D module}
\ShowEq{unitarity law, D module}
\end{itemize}
for any
\ShowEq{definition of module}
\end{theorem}
\begin{proof}
Since transformation $a$ is endomorphism of the Abelian group,
we obtain the equation \EqRef{distributive law, D module, 1}.
Since representation is homomorphism of the aditive group of the ring $D$,
we obtain the equation \EqRef{distributive law, D module, 2}.
Since the representation of the ring $D$ is representation
of the multiplicative group of
the ring $D$,
we obtain the equations \EqRef{associative law, D module} and
\EqRef{unitarity law, D module}.
\end{proof}

Vectors
\ShowEq{Vector A subs row, 1}
of $D$\Hyph module $A$ are \AddIndex{$D$\Hyph linearly independent}
{D linearly independent, module}\footnote{I follow to the
definition in \citeBib{Serge Lang}, p. 130.}
if $c=0$ follows from the equation
\ShowEq{rcd linearly independent}
Otherwise vectors \ShowEq{Vector A subs row, 1}
are \AddIndex{$D$\Hyph linearly dependent}{D linearly dependent, module}.

\begin{definition}
\label{definition: free module over ring}
We call set of vectors
\ShowEq{basis, module}
a \AddIndex{$D$\Hyph basis for module}{D basis, module}
if vectors $\Veb$ are
$D$\Hyph linearly independent and adding to this system any other vector
we get a new system which is $D$\Hyph linearly dependent.
$A$ is \AddIndex{free module over ring}{free module over ring} $D$,
if $A$ has basis
over ring $D$.\footnote{I follow to the
definition in \citeBib{Serge Lang}, p. 135.}
\qed
\end{definition}

Following definition is consequence of definitions
\ref{definition: module over ring}
and \xRef{8443-0072}{definition: morphism of representations of F algebra}.

\begin{definition}
\label{definition: linear map from A1 to A2, module}
Let $A_1$ and
$A_2$ be modules over the ring $R$.
Morphism
\ShowEq{linear map from A1 to A2}
of representation of the ring $D$ in the Abelian group $A_1$
into representation of the ring $D$ in the Abelian group $A_2$
is called
\AddIndex{linear map}{linear map}
of $D$\Hyph module $A_1$ into $D$\Hyph module $A_2$.
\qed
\end{definition}

\begin{theorem}
\label{theorem: linear map from A1 to A2}
Linear map
\ShowEq{linear map from A1 to A2}
of $D$\Hyph module $A_1$
into $D$\Hyph module $A_2$
satisfies to equations\footnote{In some books
(for instance, \citeBib{Serge Lang}, p. 119)
the theorem \ref{theorem: linear map from A1 to A2}
is considered as a definition.}
\ShowEq{linear map from A1 to A2, 1}
\end{theorem}
\begin{proof}
From definition
\ref{definition: linear map from A1 to A2, module}
and theorem
\xRef{8443-0072}{theorem: morphism of representations of algebra, reduce}
it follows that
the map $g$ is a homomorphism of the Abelian group $A_1$
into the Abelian group $A_2$ (the equation
\EqRef{linear map from A1 to A2, 1 1}).
The equation
\EqRef{linear map from A1 to A2, 1 2}
follows from the equation
\xEqRef{8443-0072}{morphism of representations of F algebra}.
\end{proof}

\begin{theorem}
Let $A_1$ and
$A_2$ be modules over ring $D$.
The set $\LDA$
is an Abelian group
relative composition law
\ShowEq{sum of linear maps}
\ShowEq{sum of linear maps, definition}
which is called
\AddIndex{sum of linear maps}{sum of linear maps}.
\end{theorem}
\begin{proof}
According to the definition
\ref{definition: linear map from A1 to A2, module}
\ShowEq{sum of linear maps 1, module}
Therefore, the map
defined by equation \EqRef{sum of linear maps, definition}
is linear map of $D$\Hyph module $A_1$
into $D$\Hyph module $A_2$.
Commutativity and associativity of sum
follow from equations
\ShowEq{sum of linear maps 2, module}

Let us define map
\ShowEq{0 circ x}
It is evident that
\ShowEq{0 in L}
From the equation
\ShowEq{0+f circ}
it follows that
\ShowEq{0+f}

Let us define map
\ShowEq{-f circ x}
It is evident that
\ShowEq{-f in L}
From the equation
\ShowEq{-f+f circ}
it follows that
\ShowEq{-f+f}

Therefore, the set $\LDA$
is an Abelian group.
\end{proof}

\begin{theorem}
Let $A_1$ and
$A_2$ be modules over ring $D$.
The representation of the ring $D$ in the Abelian group $\LDA$
which is defined by the equation
\ShowEq{product of map over scalar}
\ShowEq{product of map over scalar, definition}
is called
\AddIndex{product of map over scalar}
{product of map over scalar}.
This representation generates the structure of $D$\Hyph module
in the Abelian group $\LDA$.
\end{theorem}
\begin{proof}
From equations
\ShowEq{f->df 1}
it follows that map
\ShowEq{f->df}
is a transformation of the set $\LDA$.
From equation
\ShowEq{f->df 2}
it follows that map
\EqRef{f->df}
is a homomorphism of the Abelian group $\LDA$.
According to the definition
\xRef{8433-5163}{definition: module over ring},
the Abelian group $\LDA$ is $D$\Hyph modules.
\end{proof}

\begin{definition}
\label{definition: conjugated D module}
Let $A$ be $D$\Hyph  module.
$D$\Hyph  module
\ShowEq{conjugated D module}
is called
\AddIndex{conjugated $D$\Hyph module}{conjugated D module}.
\qed
\end{definition}

According to the definition
\ref{definition: conjugated D module},
elements of conjugated $D$\Hyph module
are $D$\Hyph linear maps
\ShowEq{conjugated D module 1}
$D$\Hyph linear map
\EqRef{conjugated D module 1}
is called
\AddIndex{linear functional on $D$\Hyph module}
{linear functional, D module}
$A$ or just
\AddIndex{$D$\Hyph linear functional}{D linear functional}.

\begin{theorem}
Let
\ShowEq{basis of module A}
be basis of $D$\Hyph module $A$.
Let $A'$ be $D$\Hyph module, conjugated to $D$\Hyph module $A$.
The set of vectors
\ShowEq{basis of module A'}
such that
\ShowEq{basis of module A A'}
is basis of $D$\Hyph module $A'$.
\end{theorem}
\begin{proof}
From the equation
\EqRef{basis of module A A'}
it follows that
\ShowEq{basis of module A A' 1}
Let
\ShowEq{f:A->D}
be linear map. Then
\ShowEq{basis of module A A' 2}
where
\ShowEq{fi}
From equations
\EqRef{basis of module A A' 1},
\EqRef{basis of module A A' 2}
it follows that
\ShowEq{basis of module A A' 3}
According to definitions
\EqRef{sum of linear maps, definition},
\EqRef{product of map over scalar, definition}
\ShowEq{basis of module A' 1}
Therefore, the set
\ShowEq{basis of module A'}
is basis of $D$\Hyph module $A'$.
\end{proof}

\begin{corollary}
\label{corollary: basis of module A'}
Let
\ShowEq{basis of module A}
be basis of $D$\Hyph module $A$.
Then
\ShowEq{ai=ei a}
\qed
\end{corollary}

Basis
\ShowEq{basis of module A'}
is called
\AddIndex{basis dual to basis}{dual basis}
\ShowEq{basis of module A}.

\begin{theorem}
\label{theorem: basis of module A1 A2}
Let
\ShowEq{basis of module A1}
be basis of $D$\Hyph module $A_1$.
Let
\ShowEq{basis of module A2}
be basis of $D$\Hyph module $A_2$.
The set of vectors
\ShowEq{basis of module A1 A2}
defined by equation
\ShowEq{basis L(A1,A2)}
is basis of $D$\Hyph module $\LDA$.
\end{theorem}
\begin{proof}
Let
\ShowEq{f:A1->A2}
be map of $D$\Hyph module $A_1$ with basis $\Basis e_1$
into $D$\Hyph module $A_2$ with basis $\Basis e_2$.
Let
\ShowEq{a in A1}
According to the corollary
\ref{corollary: basis of module A'}
\ShowEq{fa=,1}
Since
\ShowEq{fa=,2}
then from the equation
\EqRef{fa=,1}
it follows that
\ShowEq{fa=,3}
From equations
\EqRef{basis L(A1,A2)},
\EqRef{fa=,3}
it follows that
\ShowEq{fa=,4}
Therefore
\ShowEq{fa=,5}
Since maps
\EqRef{basis L(A1,A2)}
are linearly independent,
then the set of these maps is basis.
\end{proof}

\begin{definition}
\label{definition: polylinear map of modules}
Let $D$ be the commutative ring.
Let $A_1$, ..., $A_n$, $S$ be $D$\Hyph modules.
We call map
\ShowEq{polylinear map of algebras}
\AddIndex{polylinear map}{polylinear map} of modules
$A_1$, ..., $A_n$
into module
$S$,
if
\ShowEq{polylinear map of algebras, 1}
\qed
\end{definition}

\section{Algebra over Ring}
\label{section: Algebra over Ring}

\begin{definition}
\label{definition: algebra over ring}
Let $D$ be commutative ring.\footnote{This section
is written on the base of the section
\xRef{8433-5163}{section: Algebra over Ring}.}
$A$ is an \AddIndex{algebra over ring}{algebra over ring} $D$
or
\AddIndex{$D$\Hyph algebra}{D algebra},
if $A$ is $D$\Hyph module and
we defined product\footnote{I follow the definition given in
\citeBib{Richard D. Schafer}, p. 1,
\citeBib{0105.155}, p. 4.  The statement which
is true for any $D$\Hyph module,
is true also for $D$\Hyph algebra.} in $A$
\ShowEq{product in algebra, definition 2}
where $f$ is bilinear map
\ShowEq{product in algebra, definition 1}
If $A$ is free
$D$\Hyph module, then $A$ is called
\AddIndex{free algebra over ring}{free algebra over ring} $D$.
\qed
\end{definition}

According to construction that was done in subsections
\xRef{8443-0072}{subsection: Free Algebra over Ring},
\xRef{8443-0072}{subsection: n algebra},
a diagram of representations of $D$\Hyph algebra has form
\ShowEq{diagram of representations of D algebra}
On the diagram of representations
\EqRef{diagram of representations of D algebra},
$D$ is ring, $A$ is Abelian group.
We initially consider the vertical representation,
and then we consider the horizontal representation.

\begin{theorem}
The multiplication in the algebra $A$ is distributive over addition.
\end{theorem}
\begin{proof}
The statement of the theorem follows from the chain of equations
\ShowEq{product distributive in algebra}
\end{proof}

\begin{definition}
\label{definition: linear map from A1 to A2, algebra}
Let $A_1$ and
$A_2$ be algebras over ring $D$.
The linear map
\ShowEq{linear map from A1 to A2}
of the $D$\hyph module $A_1$
into the $D$\hyph module $A_2$
is called
\AddIndex{linear map}{linear map}
of $D$\Hyph algebra $A_1$ into $D$\Hyph algebra $A_2$.
Let us denote
\ShowEq{set linear maps, algebra}
set of linear maps
of algebra
$A_1$
into algebra
$A_2$.
\qed
\end{definition}

\begin{theorem}
If we define product
\ShowEq{(fg)a=f(ga)}
on $D$\Hyph module $\LDa$,
then $\LDa$ is $D$\Hyph algebra.
\end{theorem}
\begin{proof}
From equations
\ShowEq{fg 1}
it follows that map
\ShowEq{fg}
is bilinear map.
According to the definition
\ref{definition: algebra over ring},
$D$\Hyph module $\LDa$ is $D$\Hyph algebra.
\end{proof}

The multiplication in algebra can be neither commutative
nor associative. Following definitions are based
on definitions given in \citeBib{Richard D. Schafer}, p. 13.

\begin{definition}
\label{definition: commutator of algebra}
The \AddIndex{commutator}{commutator of algebra}
\ShowEq{commutator of algebra}
measures commutativity in $D$\Hyph algebra $A$.
$D$\Hyph algebra $A$ is called
\AddIndex{commutative}{commutative D algebra},
if
\ShowEq{commutative D algebra}
\qed
\end{definition}

\begin{definition}
\label{definition: associator of algebra}
The \AddIndex{associator}{associator of algebra}
\ShowEq{associator of algebra}
measures associativity in $D$\Hyph algebra $A$.
$D$\Hyph algebra $A$ is called
\AddIndex{associative}{associative D algebra},
if
\ShowEq{associative D algebra}
\qed
\end{definition}

\begin{definition}
\label{definition: nucleus of algebra}
The set\footnote{The definition is based on the similar definition in
\citeBib{Richard D. Schafer}, p. 13}
\ShowEq{nucleus of algebra}
is called the
\AddIndex{nucleus of $D$\Hyph algebra $A$}{nucleus of algebra}.
\qed
\end{definition}

\begin{definition}
\label{definition: center of algebra}
The set\footnote{The definition is based on
the similar definition in
\citeBib{Richard D. Schafer}, p. 14}
\ShowEq{center of algebra}
is called the
\AddIndex{center of $D$\Hyph algebra $A$}{center of algebra}.
\qed
\end{definition}

\begin{theorem}
Let $\Basis e$ be the basis of free finite dimensional algebra $A$ over ring $D$.
Let
\ShowEq{a b in basis of algebra}
We can get the product of $a$, $b$ according to rule
\DrawEq{product ab in algebra}{1}
where
\ShowEq{structural constants of algebra}
are \AddIndex{structural constants}{structural constants}
of algebra $A$ over ring $D$.
The product of basis vectors in the algebra $A$ is defined according to rule
\DrawEq{product of basis vectors}{algebra}
\end{theorem}
\begin{proof}
The equation
\eqRef{product of basis vectors}{algebra}
is corollary of the statement that $\Basis e$
is the basis of the algebra $A$.
Since the product in the algebra is a bilinear map,
then we can write the product of $a$ and $b$ as
\DrawEq{product ab in algebra, 1}{1}
From equations
\eqRef{product of basis vectors}{algebra},
\eqRef{product ab in algebra, 1}{1},
it follows that
\DrawEq{product ab in algebra, 2}{1}
Since $\Basis e$ is a basis of the algebra $A$, then the equation
\eqRef{product ab in algebra}{1}
follows from the equation
\eqRef{product ab in algebra, 2}{1}.
\end{proof}

\section{\texorpdfstring{$D$}{D}\hyph algebra with Hamel Basis}

If $D$\Hyph module $A$ has countable basis $\Basis e$,
then, in general,
infinite sum in $D$\Hyph module $A$ is not defined.
If continuity is not defined
in $D$\Hyph module $A$, then we use next definition
(\citeBib{0 521 59180 5}, p. 223).

\begin{definition}
\label{definition: module, Hamel basis}
Let $D$\Hyph module $A$ have countable basis
\ShowEq{Hamel basis}
If any element
of $D$\Hyph module $A$ has finite expansion relative to basis $\Basis e$,
namely, in the equation
\ShowEq{a=ai ei}
the set of values $a^{\gi i}\in D$, which are different from $0$, is finite,
then basis $\Basis e$ is called
\AddIndex{Hamel basis}{Hamel basis}.
The sequence of scalars 
\ShowEq{Hamel basis, coordinates}
is called
\AddIndex{coordinates of vector
\ShowEq{a=ai ei}
relative to Hamel basis}
{coordinates of vector, Hamel basis}
$\Basis e$.
\qed
\end{definition}

\begin{theorem}
\label{theorem: linear map, module, Hamel basis}
Let
\ShowEq{f:A1->A2}
be map of $D$\Hyph module $A_1$ with Hamel basis $\Basis e_1$
into $D$\Hyph module $A_2$ with Hamel basis $\Basis e_2$.
Let $f^{\gi i}_{\gi j}$ be coordinates of the map $f$
relative to bases $\Basis e_1$ and $\Basis e_2$.
Then for any $\gi j$,
the set of values $f^{\gi i}_{\gi j}$, which are different from $0$, is finite.
\end{theorem}
\begin{proof}
The theorem follows from the equation
\ShowEq{f:A1->A2, 1}
\end{proof}

\begin{theorem}
Let
\ShowEq{f:A1->A2}
be linear map of $D$\Hyph module $A_1$ with basis $\Basis e_1$ 
into $D$\Hyph module $A_2$ with Hamel basis $\Basis e_2$.
Then for any $a_1\in A_1$, the image
\DrawEq{f:A1->A2, 2}{Hamel}
is defined properly.
\end{theorem}
\begin{proof}
Let
\ShowEq{a1 in A1}
According to the definition
\ref{definition: module, Hamel basis},
the set of values $a^{\gi j}$, which are different from $0$, is finite.
Let $f^{\gi i}_{\gi j}$ be coordinates of the map $f$
relative to bases $\Basis e_1$ and $\Basis e_2$.
According to the theorem
\ref{theorem: linear map, module, Hamel basis},
for any $\gi j$,
the set of values $f^{\gi i}_{\gi j}$, which are different from $0$, is finite.
The union of finite set of finite sets is
finite set.
Therefore,
the set of values $a_1^{\gi j}f^{\gi i}_{\gi j}$, which are different from $0$,
is finite.
According to the definition
\ref{definition: module, Hamel basis},
expression
\eqRef{f:A1->A2, 2}{Hamel}
is expansion of the element $a_2$ relative to Hamel basis $\Basis e_2$.
\end{proof}

\begin{convention}
Let $\Basis e$ be Hamel basis of free $D$\Hyph algebra $A$.
The product of basis vectors in $D$\Hyph algebra $A$
is defined according to rule
\DrawEq{product of basis vectors}{Hamel basis}
where
\ShowEq{structural constants of algebra}
are \AddIndex{structural constants}{structural constants}
of $D$\Hyph algebra $A$.
Since the product of vectors of the basis $\Basis e$ of
$D$\Hyph algebra $A$ is a vector of $D$\Hyph algebra $A$, then
we require that
for any $\gi i$, $\gi j$,
the set of values $C^{\gi k}_{\gi{ij}}$, which are different from $0$,
is finite.
\qed
\end{convention}

\begin{theorem}
Let $\Basis e$ be Hamel basis of free $D$\Hyph algebra $A$.
Then for any
\ShowEq{a b in basis of algebra}
product defined according to rule
\DrawEq{product ab in algebra}{Hamel basis}
is defined properly.
\end{theorem}
\begin{proof}
Since the product in the algebra is a bilinear map,
then we can write the product of $a$ and $b$ as
\DrawEq{product ab in algebra, 1}{Hamel basis}
From equations
\eqRef{product of basis vectors}{Hamel basis},
\eqRef{product ab in algebra, 1}{Hamel basis},
it follows that
\DrawEq{product ab in algebra, 2}{Hamel basis}
Since $\Basis e$ is a basis of the algebra $A$, then the equation
\eqRef{product ab in algebra}{Hamel basis}
follows from the equation
\eqRef{product ab in algebra, 2}{Hamel basis}.

Since the basis $\Basis e$ is Hamel basis, then
\begin{itemize}
\item the set of values $a^{\gi i}$, which are different from $0$, is finite;
\item the set of values $b^{\gi j}$, which are different from $0$, is finite.
\end{itemize}
Therefore, the set of products $a^{\gi i}b^{\gi j}$,
which are different from $0$, is finite.
For any $\gi i$, $\gi j$,
the set of values $C^{\gi k}_{\gi{ij}}$, which are different from $0$, is finite.
Therefore, the product is properly defined by the equation
\eqRef{product ab in algebra}{Hamel basis}.
\end{proof}

\begin{theorem}
\label{theorem: standard component of tensor, algebra, Hamel basis}
Let $A_1$, ..., $A_n$ be
free algebras over commutative ring $D$.
Let $\Basis e_i$ be Hamel basis of $D$\Hyph algebra $A_i$.
Then the set of vectors
\ShowEq{tensor, algebra, Hamel basis}
is Hamel basis of tensor product
\ShowEq{tensor, algebra, Hamel basis, 1}
\end{theorem}
\begin{proof}
To prove the theorem, we need to consider the diagram
\xEqRef{8433-5163}{diagram top, algebra, tensor product}
which we used to prove the theorem
\xRef{8433-5163}{theorem: tensor product of algebras is module}.
\ShowEq{standard component of tensor, algebra, Hamel basis, diagram}
Let $M_1$ be module over ring $D$ generated
by product $\Times$
of $D$\Hyph algebras $A_1$, ..., $A_n$.
\begin{itemize}
\item
Let vector $b\in M_1$ have finite expansion
relative to the basis $\Times$
\ShowEq{b in M_1, Hamel basis}
where $I_1$ is finite set.
Let vector $c\in M_1$ have finite expansion
relative to the basis $\Times$
\ShowEq{c in M_1, Hamel basis}
where $I_2$ is finite set.
The set
\ShowEq{I=I1+I2}
is finite set.
Let
\ShowEq{b I2, c I1}
Then
\ShowEq{b+c in M_1, Hamel basis}
where $I$ is finite set.
Similarly, for $d\in D$
\ShowEq{db in M_1, Hamel basis}
where $I_1$ is finite set.
Therefore, we proved the following statement.\footnote{The set
$\Times$ cannot be Hamel basis because
this set is not countable.}
\begin{lemma}
\label{lemma: M is submodule}
The set $M$ of vectors of module $M_1$, which have finite expansion
relative to the basis $\Times$, is submodule of module $M_1$.
\end{lemma}
\end{itemize}
Injection
\ShowEq{map i, algebra, tensor product, 1}
is defined according to rule
\ShowEq{map i, algebra, Hamel basis, tensor product}
Let $N\subset M$ be submodule generated by elements of the following type
\ShowEq{kernel, algebra, Hamel basis, tensor product}
where $d_i\in A_i$, $c_i\in A_i$, $a\in D$.
Let
\[
j:M\rightarrow M/N
\]
be canonical map on factor module.
Since elements \EqRef{kernel 1, algebra, Hamel basis, tensor product}
and \EqRef{kernel 2, algebra, Hamel basis, tensor product}
belong to kernel of linear map $j$,
then, from equation
\EqRef{map i, algebra, Hamel basis, tensor product},
it follows
\ShowEq{f, algebra, Hamel basis, tensor product}
From equations \EqRef{f 1, algebra, Hamel basis, tensor product}
and \EqRef{f 2, algebra, Hamel basis, tensor product}
it follows that map $f$ is polylinear over ring $D$.

The module $M/N$ is tensor product $A_1\otimes...\otimes A_n$;
the map $j$ has form
\ShowEq{map j, Hamel basis}
and the set of
tensors like
\ShowEq{tensor, algebra, Hamel basis}
is countable basis of the module $M/N$.
According to the lemma
\ref{lemma: M is submodule},
arbitrary vector
\ShowEq{b in M, Hamel basis, 0}
has representation
\ShowEq{b in M, Hamel basis}
where $I$ is finite set.
According to the definition
\EqRef{map j, Hamel basis}
of the map $j$
\ShowEq{j(b), Hamel basis}
where $I$ is finite set.
Since $\Basis e_k$ is Hamel basis of $D$\Hyph algebra $A_k$,
then for any set of indexes $k\cdot i$, in equation
\ShowEq{aki, Hamel basis}
the set of values
\ShowEq{aki 1, Hamel basis}
which are different from $0$, is finite.
Therefore, the equation
\EqRef{j(b), Hamel basis}
has form
\ShowEq{j(b), 1, Hamel basis}
where the set of values
\ShowEq{j(b), 2, Hamel basis}
which are different from $0$, is finite.
\end{proof}

\begin{corollary}
\label{corollary: tensor product, Hamel basis}
Let $A_1$, ..., $A_n$ be
free algebras over commutative ring $D$.
Let $\Basis e_i$ be Hamel basis of $D$\Hyph algebra $A_i$.
Then any tensor $a\in A_1\otimes...\otimes A_n$
has finite set of standard components
different from $0$.
\qed
\end{corollary}

\begin{theorem}
Let $A_1$ be algebra
over the ring $D$.
Let $A_2$ be free associative algebra
over the ring $D$ with Hamel basis $\Basis e$.
The map
\ShowEq{standard representation of map A1 A2, 1, Hamel basis}
generated by the map $f\in\LDA$
through the tensor $a\in\ATwo$, has the standard representation
\ShowEq{standard representation of map A1 A2, 2, Hamel basis}
\end{theorem}
\begin{proof}
According to theorem
\ref{theorem: standard component of tensor, algebra, Hamel basis},
the standard representation of the tensor $a$ has form
\ShowEq{standard representation of map A1 A2, 3, Hamel basis}
The equation
\EqRef{standard representation of map A1 A2, 2, Hamel basis}
follows from equations
\EqRef{standard representation of map A1 A2, 1, Hamel basis},
\EqRef{standard representation of map A1 A2, 3, Hamel basis}.
\end{proof}

%% file: Hamel.Basis.Eq.tex

\def\Eij{e_{\gi i}e_{\gi j}}
\def\ATwo{A_2\otimes A_2}
\def\Times{A_1\times...\times A_n}
\def\LDA{\mathcal L(D;A_1;A_2)}
\def\LDa{\mathcal L(D;A;A)}

\def\AcF{a\circ f}
\def\GAcF{g=\AcF}

\def\Veb{e_i}
\def\Ii{i\in I}

\DefEq
{
\begin{align*}
((f_1+f_2)\circ g)\circ a&=(f_1+f_2)\circ (g\circ a)
=f_1\circ (g\circ a)+f_2\circ (g\circ a)
\\
&=(f_1\circ g)\circ a+(f_2\circ g)\circ a
=(f_1\circ g+f_2\circ g)\circ a
\\
((df)\circ g)\circ a&=(df)\circ (g\circ a)=d(f\circ (g\circ a))
\\
&=d((f\circ g)\circ a)=(d(f\circ g))\circ a
\\
(f\circ (g_1+g_2))\circ a&=f\circ ((g_1+g_2)\circ a)
=f\circ (g_1\circ a+g_2\circ a)
\\
&=f\circ (g_1\circ a)+f\circ (g_2\circ a)
\\
&=(f\circ g_1)\circ a+(f\circ g_2)\circ a
=(f\circ g_1+f\circ g_2)\circ a
\\
(f\circ (dg))\circ a&=f\circ ((dg)\circ a)
=f\circ (d(g\circ a))=d(f\circ (g\circ a))
\\
&=d((f\circ g)\circ a)=(d(f\circ g))\circ a
\end{align*}
}
{fg 1}

\DefEq
{
$\Basis e=(e_{\gi i}\in A,\gi i\in\gi I)$
}
{basis of module A}

\DefEq
{
$\Basis e_1=(e_{1\cdot\gi i}\in A_1,\gi i\in\gi I)$
}
{basis of module A1}

\DefEq
{
$\Basis e_2=(e_{2\cdot\gi j}\in A_2,\gi j\in\gi J)$
}
{basis of module A2}

\DefEq
{
\symb{(e^{\gi i}_1,e_{2\cdot\gi j})}1{basis L(A1,A2)},
$\gi i\in\gi I$, $\gi j\in\gi J$,
}
{basis of module A1 A2}

\DefEq
{
$a\in A_1$, $a=a^{\gi i}e_{1\cdot\gi i}$.
}
{a in A1}

\DefEquation
{
f\circ a=f\circ(a^{\gi i}e_{1\cdot\gi i})
=a^{\gi i}(f\circ e_{1\cdot\gi i})
=(e^{\gi i}_1\circ a)(f\circ e_{1\cdot\gi i})
}
{fa=,1}

\DefEquation
{
f\circ a=(e^{\gi i}_1\circ a)f^{\gi j}_{\gi i} e_{2\cdot\gi j}
}
{fa=,3}

\DefEq
{
\[f\circ a=f^{\gi j}_{\gi i}(e^{\gi i}_1,e_{2\cdot\gi j})\circ a\]
}
{fa=,4}

\DefEquation
{
f:A\rightarrow D
}
{conjugated D module 1}

\DefEq
{
\[f=f^{\gi j}_{\gi i}(e^{\gi i}_1,e_{2\cdot\gi j})\]
}
{fa=,5}

\DefEquation
{
\ShowSymbol{basis L(A1,A2)}\circ a=(e^{\gi i}_1\circ a) e_{2\cdot\gi j}
}
{basis L(A1,A2)}

\DefEq
{
$f\circ e_{1\cdot\gi i}\in A_2$,
}
{fa=,2}

\DefEq
{
\[a^{\gi i}=e^{\gi i}\circ a\]
}
{ai=ei a}

\DefEq
{
$\Basis e=(e^{\gi i}\in A',\gi i\in\gi I)$
}
{basis of module A'}

\DefEquation
{
e^{\gi i}\circ e_{\gi j}=\delta^{\gi i}_{\gi j}
}
{basis of module A A'}

\DefEq
{
\[f:A\rightarrow D\]
}
{f:A->D}

\DefEquation
{
f\circ a=f\circ(a^{\gi i}e_{\gi i})=a^{\gi i}(f\circ e_{\gi i})
=a^{\gi i}f_{\gi i}
}
{basis of module A A' 2}

\DefEquation
{
f\circ a=f_{\gi i}(e^{\gi i}\circ a)
}
{basis of module A A' 3}

\DefEq
{
\[f=f_{\gi i}e^{\gi i}\]
}
{basis of module A' 1}

\DefEq
{
$f_{\gi i}=f\circ e_{\gi i}$
}
{fi}

\DefEquation
{
e^{\gi i}\circ a=e^{\gi i}\circ (a^{\gi j}e_{\gi j})
=a^{\gi j}(e^{\gi i}\circ e_{\gi j})
=a^{\gi j}\delta^{\gi i}_{\gi j}
=a^{\gi i}
}
{basis of module A A' 1}

\DefEq
{
\symb{A'}0{conjugated D module}
$\ShowSymbol{conjugated D module}=\mathcal L(D;A;D)$
}
{conjugated D module}

\DefEq
{
\begin{align}
f\circ(a+b)&=f\circ a+f\circ b
\EqLabel{linear map from A1 to A2, 1 1}
\\
f\circ(pa)&=p(f\circ a)
\EqLabel{linear map from A1 to A2, 1 2}
\end{align}
\[
\begin{matrix}
a,b\in A_1
&
p\in D
\end{matrix}
\]
}
{linear map from A1 to A2, 1}

\DefEq
{
\begin{align*}
f\circ(
a_1, ...,
a_i+ b_i, ...,
a_n)
&=
f\circ(
a_1, ...,
a_i, ...,
a_n)
+
f\circ(
a_1, ...,
b_i, ...,
a_n)
\\
f\circ(
a_1, ...,
pa_i, ...,
a_n)
&=
pf\circ(
a_1, ...,
a_i, ...,
a_n)
\end{align*}
\[
\begin{matrix}
1\le i\le n
&
a_i, b_i \in A_i
&
p\in D
\end{matrix}
\]
}
{polylinear map of algebras, 1}

\DefEq
{
$a$, $b \in D$, $m$, $n \in A$.
}
{definition of module}

\DefEq
{
\[
f:A_1\times...\times A_n\rightarrow
S
\]
}
{polylinear map of algebras}

\DefEq
{
\[
f:A_1\rightarrow A_2
\]
}
{linear map from A1 to A2}

\DefEq
{
$f\circ g$
}
{fg}

\DefEquation
{
(f\circ g)\circ a=f\circ (g\circ a)
}
{(fg)a=f(ga)}

\DefEq
{
\symb{\LDA}1{set linear maps, scalar}
}
{set linear map, algebra}

\DefEq
{
\symb{f+g}0{sum of linear maps}
}
{sum of linear maps}

\DefEquation
{
(\ShowSymbol{sum of linear maps})\circ x=f\circ x+g\circ x
}
{sum of linear maps, definition}

\DefEq
{
\symb{df}0{product of map over scalar}
}
{product of map over scalar}

\DefEq
{
\begin{align*}
(df)\circ(d_1a)&=d(f\circ(d_1a))=d_1(d(f\circ a))=d_1((df)\circ a)
\\
(df)\circ(a_1+a_2)&=d(f\circ(a_1+a_2))=d(f\circ a_1+f\circ a_2)
\\
&=d(f\circ a_1)+d(f\circ a_2)=(df)\circ a_1+(df)\circ a_2
\end{align*}
}
{f->df 1}

\DefEq
{
\begin{align*}
(d(f+g))\circ a&=d((f+g)\circ a)=d(f\circ a+g\circ a)
\\
&=d(f\circ a)+d(g\circ a)=(df)\circ a+(dg)\circ a
\end{align*}
}
{f->df 2}

\DefEquation
{
f\rightarrow df
}
{f->df}

\DefEquation
{
(\ShowSymbol{product of map over scalar})\circ x=d(f\circ x)
}
{product of map over scalar, definition}

\DefEq
{
\begin{align*}
(f+g)\circ(a+b)&=f\circ(a+b)+g\circ(a+b)
\\
&=f\circ a+f\circ b+g\circ a+g\circ b
\\
&=(f+g)\circ a+(f+g)\circ b
\\
(f+g)\circ(da)&=f\circ(da)+g\circ(da)
\\
&=df\circ a+dg\circ a
\\
&=d(f+g)\circ a
\end{align*}
}
{sum of linear maps 1, module}

\DefEq
{
\begin{align*}
(f+g)\circ a&=f\circ a+g\circ a
=g\circ a+f\circ a
\\
&=(g+f)\circ a
\\
((f+g)+h)\circ a&=(f+g)\circ a+h\circ a
=(f\circ a+g\circ a)+h\circ a
\\
&=f\circ a+(g\circ a+h\circ a)
=f\circ a+(g+h)\circ a
\\
&=(f+(g+h))\circ a
\end{align*}
}
{sum of linear maps 2, module}

\DefEq
{
\[(0+f)\circ x=0\circ x+f\circ x=0+f\circ x=f\circ x\]
}
{0+f circ}

\DefEq
{
\[0+f=f\]
}
{0+f}

\DefEq
{
$0\in\LDA$.
}
{0 in L}

\DefEq
{
$0\circ x=0$.
}
{0 circ x}

\DefEq
{
\[((-f)+f)\circ x=(-f)\circ x+f\circ x=(-(f\circ x))+f\circ x=0=0\circ x\]
}
{-f+f circ}

\DefEq
{
\[(-f)+f=0\]
}
{-f+f}

\DefEq
{
$-f\in\LDA$.
}
{-f in L}

\DefEq
{
\[(-f)\circ x=-(f\circ x)\]
}
{-f circ x}

\DefEq
{
\[
f:A_1\rightarrow A_2
\]
}
{linear map from A1 to A2}

\newcommand\Mapj
{
j(a_1,...,a_n)=a_1\otimes...\otimes a_n
}
\newcommand\mapj[1]
{
\begin{align}
&(d_1,...,d_i+c_i,...,d_n)-
(d_1,..., d_i,..., d_n)-
(d_1,..., c_i,..., d_n)
\EqLabel{kernel 1, algebra, #1, tensor product}
\\
&(d_1,...,ad_i,..., d_n)-
a(d_1,..., d_i,..., d_n)
\EqLabel{kernel 2, algebra, #1, tensor product}
\end{align}
}

\newcommand\mapf[1]
{
\begin{align}
f\circ(d_1,...,d_i+c_i,..., d_n)
=&f\circ(d_1,..., d_i,..., d_n)+
f\circ(d_1,..., c_i,..., d_n)
\EqLabel{f 1, algebra, #1, tensor product}
\\
f\circ(d_1,...,ad_i,..., d_n)=&
a\ f\circ(d_1,..., d_i,..., d_n)
\EqLabel{f 2, algebra, #1, tensor product}
\end{align}
}

\newcommand\Dpr[3]
{
\begin{matrix}
\vcenter
{
\xymatrix
{
D_{#3}\ar[r]|(.4){*}^{f_{#2 1,2}}&
A_{#1}\ar[r]|{*}^{f_{#2 2,3}}&
A_{#1}
\\
&&D_{#3}\ar[u]|{*}_{f_{#2 1,2}}
}
}
&
\begin{array}{r@{\,}l}
f_{#2 1,2}(d):v&\rightarrow d\,v
\\
f_{#2 2,3}(v):w&\rightarrow
C_{#1}(v, w)
\\
C_{#1}&\in\mathcal L(A_{#1}^2;A_{#1})
\end{array}
\end{matrix}
}

\DefEq
{
\[
\begin{matrix}
b=b^i(a_{1\cdot i},...,a_{n\cdot i})
&i\in I_1
\end{matrix}
\]
}
{b in M_1, Hamel basis}

\DefEq
{
\[
\begin{matrix}
b=b^i(a_{1\cdot i},...,a_{n\cdot i})
&i\in I
\end{matrix}
\]
}
{b in M, Hamel basis}

\DefEq
{
\[
a_{k\cdot i}=a^{\gi{p_k}}_{k\cdot i}e_{k\cdot\gi{p_k}}
\]
}
{aki, Hamel basis}

\DefEq
{
$a^{\gi{p_k}}_{k\cdot i}$,
}
{aki 1, Hamel basis}

\DefEquation
{
\begin{matrix}
j\circ b=b^i(a_{1\cdot i}\otimes...\otimes a_{n\cdot i})
&i\in I
\end{matrix}
}
{j(b), Hamel basis}

\DefEquation
{
\begin{matrix}
j\circ b=b^ia^{\gi{p_1}}_{1\cdot i}...a^{\gi{p_n}}_{n\cdot i}
(e_{1\cdot\gi{p_1}}\otimes...\otimes e_{n\cdot\gi{p_n}})
&i\in I
\end{matrix}
}
{j(b), 1, Hamel basis}

\DefEq
{
\[
b^ia^{\gi{p_1}}_{1\cdot i}...a^{\gi{p_n}}_{n\cdot i}
\]
}
{j(b), 2, Hamel basis}

\DefEq
{
\symb{\{e_{\gi i}\}_{\gi i=\gi 1}^{\gi\infty}}0{Hamel basis}
$\Basis e=\ShowSymbol{Hamel basis}$.
}
{Hamel basis}

\DefEq
{
\symb{\{a^{\gi i}\}_{\gi i=\gi 1}^{\gi\infty}}1{Hamel basis, coordinates}
}
{Hamel basis, coordinates}

\DefEq
{
\[
\sum_{\gi i=\gi 1}^{\gi\infty}|f^{\gi i}_{\gi j}|<\infty
\]
}
{fij in D}

\DefEq
{
$e_{1\cdot\gi i}\in\Basis e_1$
}
{e1 in}

\DefEq
{
$f\circ e_{1\cdot\gi i}$
}
{f e1 in}

\DefEquation
{
\begin{matrix}
a_{1pq}=\sum_{i=p}^qa_1^{\gi i}\,e_{\gi i}
&
|a_{1pq}|<\epsilon
\end{matrix}
}
{a1 in A1, 1}

\DefEq
{
$b\in M$
}
{b in M, Hamel basis, 0}

\DefEq
{
\mapj{Hamel basis}
}
{kernel, algebra, Hamel basis, tensor product}

\DefEq
{
\mapf{Hamel basis}
}
{f, algebra, Hamel basis, tensor product}

\DefEq
{
\[
\xymatrix{
i:\Times
\ar[r]&M
}
\]
}
{map i, algebra, tensor product, 1}

\DefEq
{
i\circ(d_1,...,d_n)=(d_1,...,d_n)
}
{map i, algebra, tensor product}

\DefEquation
{
\ShowEq{map i, algebra, tensor product}
}
{map i, algebra, Hamel basis, tensor product}

\DefEq
{
\[
\begin{matrix}
db=db^i(a_{1\cdot i},...,a_{n\cdot i})
&i\in I_1
\end{matrix}
\]
}
{db in M_1, Hamel basis}

\DefEq
{
\[
\begin{matrix}
b+c=(b^i+c^i)(a_{1\cdot i},...,a_{n\cdot i})
&i\in I
\end{matrix}
\]
}
{b+c in M_1, Hamel basis}

\DefEq
{
$I=I_1\cup I_2$
}
{I=I1+I2}

\DefEq
{
\[
\begin{matrix}
b_i=0&i\in I\setminus I_1
\\
c_i=0&i\in I\setminus I_2
\end{matrix}
\]
}
{b I2, c I1}

\DefEq
{
\[
\begin{matrix}
c=c^i(a_{1\cdot i},...,a_{n\cdot i})
&i\in I_2
\end{matrix}
\]
}
{c in M_1, Hamel basis}

\DefEquation
{
\Mapj
}
{map j, Hamel basis}

\DefEquation
{
g=
a^{\gi{ij}}(e_{\gi i}\otimes e_{\gi j})\circ f
=
a^{\gi{ij}}e_{\gi i}fe_{\gi j}
}
{standard representation of map A1 A2, 2, Hamel basis}

\DefEquation
{
a=a^{\gi{ij}}e_{\gi i}\otimes e_{\gi j}
}
{standard representation of map A1 A2, 3, Hamel basis}

\DefEquation
{
\GAcF
}
{standard representation of map A1 A2, 1, Hamel basis}

\DefEquation
{
\begin{matrix}
a_1\in A_1&
a_1=a_1^{\gi i}\,e_{1\cdot\gi i}
\end{matrix}
}
{a1 in A1}

\DefEq
{
\[
a=a^{\gi i}\,e_{\gi i}
\]
}
{a=ai ei}

\DefEq
{
\[
f:A_1\rightarrow A_2
\]
}
{f:A1->A2}

\DefEq
{
\[
f\circ e_{1\cdot\gi j}=
f^{\gi i}_{\gi j}\,e_{2\cdot\gi i}
\]
}
{f:A1->A2, 1}

\DefEq
{
\begin{matrix}
a_2=f\circ a_1&
a_2^{\gi i}=a_1^{\gi j}f^{\gi i}_{\gi j}&
a_2=a_2^{\gi i}\,e_{2\cdot\gi i}
\end{matrix}
}
{f:A1->A2, 2}

\DefEquation
{
\xymatrix
{
f:D\ar[r]|{*}&A
}
}
{representation of the ring}

\DefEq
{
\[
(f+g)\circ a=f\circ a+g\circ a
\]
}
{sum of transformations of Abelian group}

\DefEq
{
\[f(ab)=f(a)\circ f(b)\]
}
{product of transformations of representation}

\DefEq
{
\[
f(a-b)\circ m=0
\]
}
{representation of ring, 2}

\DefEquation
{
(ab)\circ m=a\circ(b\circ m)
}
{associative law, D module}

\DefEq
{
\begin{align}
\EqLabel{distributive law, D module, 1}
a\circ(m+n)&=a\circ m+a\circ n\\
\EqLabel{distributive law, D module, 2}
(a+b)\circ m&=a m+b m
\end{align}
}
{distributive law, D module}

\DefEq
{
$a_i$, $\Ii$,
}
{Vector A subs row, 1}

\DefEq
{
\symb{\Basis{e}=(\Veb,\Ii)}1
{basis, module}
}
{basis, module}

\DefEq
{
\[c^i a_i =0\]
}
{rcd linearly independent}

\DefEquation
{
1 m=m
}
{unitarity law, D module}

\DefEquation
{
f(a)\circ m=f(b)\circ m
}
{representation of ring, 1}

\DefEq
{
\[
f(a+b)\circ x=f(a)\circ x+f(b)\circ x
\]
}
{sum of transformations of Abelian group, 1}

\DefEq
{
\begin{align*}
(a+b)c&=f\circ(a+b,c)=f\circ(a,c)+f\circ(b,c)=ac+bc
\\
a(b+c)&=f\circ(a,b+c)=f\circ(a,b)+f\circ(a,c)=ab+ac
\end{align*}
}
{product distributive in algebra}

\DefEq
{
\symb{[a,b]}
0{commutator of algebra}
\[
\ShowSymbol{commutator of algebra}=ab-ba
\]
}
{commutator of algebra}

\DefEq
{
\[
[a,b]=0
\]
}
{commutative D algebra}

\DefEq
{
\symb{(a,b,c)}
0{associator of algebra}
\begin{equation}
\ShowSymbol{associator of algebra}=(ab)c-a(bc)
\EqLabel{associator of algebra}
\end{equation}
}
{associator of algebra}

\DefEq
{
\[
(a,b,c)=0
\]
}
{associative D algebra}

\DefEq
{
\symb{N(A)}0
{nucleus of algebra}
\[
\ShowSymbol{nucleus of algebra}=
\{
a\in A:
\forall b, c\in A,
(a,b,c)=(b,a,c)=(b,c,a)=0
\}
\]
}
{nucleus of algebra}

\DefEq
{
\symb{Z(A)}0
{center of algebra}
\[
\ShowSymbol{center of algebra}=
\{
a\in A:
a\in N(A),
\forall b\in A,
ab=ba
\}
\]
}
{center of algebra}

\DefEq
{
\[
\begin{matrix}
a=a^{\gi i}\,e_{\gi i}&b=b^{\gi i}\,e_{\gi i}&a, b\in A
\end{matrix}
\]
}
{a b in basis of algebra}

\DefEq
{
(ab)^{\gi k}=C^{\gi k}_{\gi{ij}}a^{\gi i}b^{\gi j}
}
{product ab in algebra}

\DefEq
{
\symb{C^{\gi k}_{\gi{ij}}}1{structural constants}
}
{structural constants of algebra}

\DefEq
{
e_{\gi i}e_{\gi j}=
C^{\gi k}_{\gi{ij}}e_{\gi k}
}
{product of basis vectors}

\DefEq
{
ab=a^{\gi i}b^{\gi j}\Eij
}
{product ab in algebra, 1}

\DefEq
{
ab=a^{\gi i}b^{\gi j}
C^{\gi k}_{\gi{ij}}e_{\gi k}
}
{product ab in algebra, 2}

\DefEq
{
\[
f:A\times A\rightarrow A
\]
}
{product in algebra, definition 1}

\DefEquation
{
\xymatrix{
&&&M/N
\\
\Times\ar[rrru]^f\ar[rr]_(.7)i
&&M\ar[ur]_j
}
}
{standard component of tensor, algebra, Hamel basis, diagram}

\DefEq
{
$A_1\otimes...\otimes A_n$.
}
{tensor, algebra, Hamel basis, 1}

\DefEquation
{
ab=f\circ(a,b)
}
{product in algebra, definition 2}

\DefEq
{
$e_{1\cdot i_1}\otimes...\otimes e_{n\cdot i_n}$
}
{tensor, algebra, Hamel basis}

\DefEquation
{
\Dpr {}{}{}
}
{diagram of representations of D algebra}

%% file: Schauder.Basis.English.tex
\input{Schauder.Basis.Eq}
\ifx\PrintBook\Defined
\chapter{Schauder Basis}
\fi

\section{Topological Ring}

\input{\FilePrefix Center.Ring.English}

\begin{definition}
Ring $D$ is called 
\AddIndex{topological ring}{topological ring}\footnote{I
made definition according to definition
from \citeBib{Pontryagin: Topological Group},
chapter 4.}
if $D$ is topological space and the algebraic operations
defined in $D$ are continuous in the topological space $D$.
\qed
\end{definition}

According to definition, for arbitrary elements $a, b\in D$
and for arbitrary neighborhoods $W_{a-b}$ of the element $a-b$,
$W_{ab}$ of the element $ab$ there exists neighborhoods $W_a$
of the element $a$ and $W_b$ of the element $b$ such
that $W_a-W_b\subset W_{a-b}$,
$W_aW_b\subset W_{ab}$.

\begin{definition}
\label{definition: norm on ring}
\AddIndex{Norm on ring}
{norm on ring} $D$\footnote{I
made definition according to definition
from \citeBib{Bourbaki: General Topology: Chapter 5 - 10},
IX, \S 3.2 and definition \citeBib{Arnautov Glavatsky Mikhalev}-1.1.12,
p. 23.} is a map
\[d\in D\rightarrow |d|\in R\]
which satisfies the following axioms
\begin{itemize}
\item $|a|\ge 0$
\item $|a|=0$ if, and only if, $a=0$
\item $|ab|=|a|\ |b|$
\item $|a+b|\le|a|+|b|$
\end{itemize}

The ring $D$, endowed with the structure defined by a given norm on
$D$, is called \AddIndex{normed ring}{normed ring}.
\qed
\end{definition}

Invariant distance on additive group of ring $D$
\[d(a,b)=|a-b|\]
defines topology of metric space,
compatible with ring structure of $D$.

\begin{definition}
\label{definition: limit of sequence, normed ring}
Let $D$ be normed ring.
Element $a\in D$ is called 
\AddIndex{limit of a sequence}{limit of sequence}
$\{a_n\}$
\ShowEq{limit of sequence, normed ring}
if for every $\epsilon\in R$, $\epsilon>0$
there exists positive integer $n_0$ depending on $\epsilon$ and such,
that $|a_n-a|<\epsilon$ for every $n>n_0$.
\qed
\end{definition}

\begin{definition}
Let $D$ be normed ring.
The sequence $\{a_n\}$, $a_n\in D$ is called 
\AddIndex{fundamental}{fundamental sequence}
or \AddIndex{Cauchy sequence}{Cauchy sequence},
if for every $\epsilon\in R$, $\epsilon>0$,
there exists positive integer $n_0$ depending on $\epsilon$ and such,
that $|a_p-a_q|<\epsilon$ for every $p$, $q>n_0$.
\qed
\end{definition}

\begin{definition}
Normed ring $D$ is called
\AddIndex{complete}{complete ring}
if any fundamental sequence of elements
of ring $D$ converges, i.e.
has limit in ring $D$.
\qed
\end{definition}

Let $D$ be complete ring of characteristic $0$.
Since division in the ring, in general, is not defined,
we cannot state that the ring $D$ contains
rational field.
We will assume that considered ring $D$ contains
rational field.
Under this assumption, it is evident that the ring has characteristic $0$.

\begin{theorem}
Let $D$ be ring containing
rational field and let $d\in D$.
Then for any integer $n\in Z$
\ShowEq{inverse integer in ring}
\end{theorem}
\begin{proof}
According to theorem \ref{theorem: integer subring of ring}
following chain of equation is true
\ShowEq{inverse integer in ring, 1}
Let us multiply right and left sides of equation
\EqRef{inverse integer in ring, 1} by $n^{-1}$.
We get
\ShowEq{inverse integer in ring, 2}
\EqRef{inverse integer in ring} follows
from \EqRef{inverse integer in ring, 2}.
\end{proof}

\begin{theorem}
\label{theorem: ring, rational number}
Let $D$ be ring containing
rational field and let $d\in D$.
Then every rational number $p\in Q$ commutes with $d$.
\end{theorem}
\begin{proof}
Let us represent rational number $p\in Q$ as
$p=mn^{-1}$, $m$, $n\in Z$.
Statement of theorem follows from chain of equations
\[
pd=mn^{-1}d=n^{-1}dm=dmn^{-1}=dp
\]
based on the statement of theorem \ref{theorem: integer subring of ring}
and equation \EqRef{inverse integer in ring}.
\end{proof}

\begin{theorem}
\label{theorem: rational subring of center}
Let $D$ be ring containing
rational field.
Then field of rational numbers $Q$ is subfield of center $Z(D)$ of ring $D$.
\end{theorem}
\begin{proof}
Corollary of theorem \ref{theorem: ring, rational number}.
\end{proof}

Later on, speaking about normed ring of characteristic $0$,
we will assume that homeomorphism of field of rational numbers $Q$
into ring $D$ is defined.

\begin{theorem}
\label{theorem: limit of sequence, normed ring, product on scalar}
Let $D$ be normed ring of characteristic $0$ and let $d\in D$.
Let $a\in D$ be limit of a sequence $\{a_n\}$.
Then
\[
\lim_{n\rightarrow\infty}(a_nd)=ad
\]
\[
\lim_{n\rightarrow\infty}(da_n)=da
\]
\end{theorem}
\begin{proof}
Statement of the theorem is trivial, however I give this proof
for completeness sake. 
Since $a\in D$ is limit of the sequence $\{a_n\}$,
then according to definition
\ref{definition: limit of sequence, normed ring}
for given $\epsilon\in R$, $\epsilon>0$,
there exists positive integer $n_0$ such, that
\ShowEq{an-a d}
for every $n>n_0$.
According to definition \ref{definition: norm on ring}
the statement of theorem follows from inequalities
\ShowEq{an-a d1}
for any $n>n_0$.
\end{proof}

\begin{theorem}
\label{theorem: complete ring contains real number}
Complete ring $D$ of characteristic $0$
contains as subfield an isomorphic image of the field $R$ of
real numbers.
It is customary to identify it with $R$.
\end{theorem}
\begin{proof}
Consider fundamental sequence of rational numbers $\{p_n\}$.
Let $p'$ be limit of this sequence in ring $D$.
Let $p$ be limit of this sequence in field $R$.
Since immersion of field $Q$ into division ring $D$ is homeomorphism,
then we may identify $p'\in D$ and $p\in R$.
\end{proof}

\begin{theorem}
\label{theorem: complete ring and real number}
Let $D$ be complete ring of characteristic $0$ and let $d\in D$.
Then any real number $p\in R$ commute with $d$.
\end{theorem}
\begin{proof}
Let us represent real number $p\in R$ as
fundamental sequence of rational numbers $\{p_n\}$.
Statement of theorem follows from chain of equations
\ShowEq{complete ring and real number}
based on statement of theorem
\ref{theorem: limit of sequence, normed ring, product on scalar}.
\end{proof}

\section{Normed \texorpdfstring{$D$}{D}-Algebra}

\ifx\texFuture\Defined
\begin{definition}
Given a topological commutative ring $D$ and
$D$\Hyph algebra $A$ such that $A$ has
a topology compatible with the structure of the additive
group of $A$ and maps
\ShowEq{topological D algebra}
are continuous,
then $V$ is called a
\AddIndex{topological $D$\Hyph algebra}
{topological D algebra}\footnote{I
made definition according to definition
from \citeBib{Bourbaki: Topological Vector Space},
p. TVS I.1}.
\qed
\end{definition}
\fi

\begin{definition}
\label{definition: norm on D module}
Let $D$ be valued commutative ring.\footnote{I
made definition according to definition
from \citeBib{Bourbaki: General Topology: Chapter 5 - 10},
IX, \S 3.3.
We use notation either $|a|$ or $\|a\|$ for norm.
}
\AddIndex{Norm on $D$\Hyph module}
{norm on D module} $A$
is a map
\ShowEq{norm on D module}
which satisfies the following axioms
\StartLabelItem{definition: norm on D module}
\begin{enumerate}
\ShowEq{norm on D module 1}
\ShowEq{norm on D module 2, 1}
if, and only if,
\ShowEq{norm on D module 2, 2}
\ShowEq{norm on D module 3}
\end{enumerate}
$D$\Hyph module $A$,
endowed with the structure defined by a given norm on
$A$, is called
\AddIndex{normed $D$\Hyph module}{normed D module}.
\qed
\end{definition}

\begin{theorem}
\label{theorem: normed module, norm of difference}
Norm in $D$\Hyph module $A$ satisfies to equation
\ShowEq{normed module, norm of difference}
\end{theorem}
\begin{proof}
From the equation
\ShowEq{normed module, norm of difference, 1}
and statement
\ref{item: norm on D module 3, 1},
it follows that
\ShowEq{normed module, norm of difference, 2}
The equation
\EqRef{normed module, norm of difference}
follows from the equation
\EqRef{normed module, norm of difference, 2}.
\end{proof}

\begin{definition}
\label{definition: normal basis}
The basis $\Basis e$ is called
\AddIndex{normal basis}{normal basis},
if
\ShowEq{Banach algebra, normal basis}
for any vector $e_{\gi i}$  of the basis $\Basis e$.
\qed
\end{definition}

\begin{definition}
\label{definition: limit of sequence, normed module}
Let $A$ be normed $D$\Hyph module.
Element $a\in A$ is called 
\AddIndex{limit of a sequence}{limit of sequence}
\ShowEq{limit of sequence, normed module}
if for every $\epsilon\in R$, $\epsilon>0$
there exists positive integer $n_0$ depending on $\epsilon$ and such, that
\ShowEq{an-a}
for every $n>n_0$.
\qed
\end{definition}

\begin{definition}
\label{definition: fundamental sequence, normed module}
Let $A$ be normed $D$\Hyph module.
The sequence $\{a_n\}$, $a_n\in A$, is called 
\AddIndex{fundamental}{fundamental sequence}
or \AddIndex{Cauchy sequence}{Cauchy sequence},
if for every $\epsilon\in R$, $\epsilon>0$,
there exists positive integer $n_0$ depending on $\epsilon$ and such, that
\ShowEq{ap-aq}
for every $p$, $q>n_0$.
\qed
\end{definition}

\begin{definition}
Normed $D$\Hyph module $A$ is called
\AddIndex{Banach $D$\Hyph module}{Banach module}
if any fundamental sequence of elements
of module $A$ converges, i.e.
has limit in module $A$.
\qed
\end{definition}

\begin{definition}
\label{definition: norm on d algebra}
Let $D$ be valued commutative ring.
Let $A$ be $D$\Hyph algebra.
The norm
\ShowEq{norm on d algebra}
on $D$\Hyph module $A$ such that
\ShowEq{norm on d algebra 3}
is called
\AddIndex{norm on $D$\Hyph algebra}
{norm on D algebra} $A$.
$D$\Hyph algebra $A$,
endowed with the structure defined by a given norm on
$A$, is called
\AddIndex{normed $D$\Hyph algebra}
{normed D algebra}.
\qed
\end{definition}

\begin{definition}
Normed $D$\Hyph algebra $A$ is called
\AddIndex{Banach $D$\Hyph algebra}{Banach algebra}
if any fundamental sequence of elements
of algebra $A$ converges, i.e.
has limit in algebra $A$.
\qed
\end{definition}

\begin{definition}
Let $A$ be Banach $D$\Hyph algebra.
Set of elements
\ShowEq{unit sphere in algebra}
is called
\AddIndex{unit sphere in algebra}{unit sphere in algebra} $A$.
\qed
\end{definition}

\section{Normed \texorpdfstring{$D$\Hyph}{D-}Module
\texorpdfstring{$\mathcal L(D;A_1;A_2)$}{L(D;A1;A2)}}
\label{section: Module L(D;A1;A2)}

\begin{definition}
\label{definition: continuous map, module}
Map
\ShowEq{f A1 A2}
of $D$\Hyph module $A_1$ with norm $\|x\|_1$
into $D$\Hyph module $A_2$ with norm $\|y\|_2$
is called \AddIndex{continuous}{continuous map}, if
for every as small as we please $\epsilon>0$
there exist such $\delta>0$, that
\DrawEq{continuous map, x}{}
implies
\DrawEq{continuous map, y}{}
\qed
\end{definition}

\begin{theorem}
Let\footnote{This theorem is based
on the theorem
\citeBib{Kolmogorov Fomin}\Hyph 1, page 77.}
\ShowEq{f A1 A2}
linear map
of $D$\Hyph module $A_1$ with norm $\|x\|_1$
into $D$\Hyph module $A_2$ with norm $\|y\|_2$.
Since the linear map $f$ is continuous at $x\in A_1$,
then the linear map $f$ is continuous everywhere in $D$\Hyph module $A_1$.
\end{theorem}
\begin{proof}
Let $\epsilon>0$.
According to the definition
\ref{definition: continuous map, module},
there exist such $\delta>0$, that
\DrawEq{continuous map, x}{at x}
implies
\ShowEq{continuous map, y at x}
Let $b\in A_1$.
From the equation
\eqRef{continuous map, x}{at x},
it follows that
\ShowEq{continuous map, x 1}
From the equation
\EqRef{continuous map, y at x},
it follows that
\ShowEq{continuous map, y 1}
Therefore, the linear map $f$ is continuous
at point $x+b$.
\end{proof}

\begin{corollary}
\label{corollary: linear map is continuous at 0}
Linear map
\ShowEq{f A1 A2}
of normed $D$\Hyph module $A_1$
into normed $D$\Hyph module $A_2$
is continuous iff it is continuous at point $0\in A_1$.
\qed
\end{corollary}

\begin{theorem}
\label{theorem: sum of continuous linear maps}
The sum of continuous linear maps
of $D$\Hyph module $A_1$ with norm $\|x\|_1$
into $D$\Hyph module $A_2$ with norm $\|y\|_2$
is continuous linear map.
\end{theorem}
\begin{proof}
Let
\ShowEq{f A1 A2}
be continuous linear map.
According to the corollary
\ref{corollary: linear map is continuous at 0}
and the definition
\ref{definition: continuous map, module},
for given $\epsilon>0$
there exist such $\delta_1>0$, that
\ShowEq{continuous map, x, 1}
implies
\ShowEq{continuous map, y, 1}
Let
\ShowEq{g A1 A2}
be continuous linear map.
According to the corollary
\ref{corollary: linear map is continuous at 0}
and the definition
\ref{definition: continuous map, module},
for given $\epsilon>0$
there exist such $\delta_2>0$, that
\ShowEq{continuous map, x, 2}
implies
\ShowEq{continuous map, y, 2}
Let
\ShowEq{delta 1 2}
From inequalities
\EqRef{continuous map, y, 1},
\EqRef{continuous map, y, 2}
and the statement
\ref{item: norm on D module 3, 1},
it follows that
\ShowEq{continuous map, x, 3}
implies
\ShowEq{continuous map, y, 3}
Therefore,
according to the corollary
\ref{corollary: linear map is continuous at 0}
and the definition
\ref{definition: continuous map, module},
linear map $f+g$ is continuous.
\end{proof}

\begin{theorem}
\label{theorem: Product of continuous map over scalar}
Let
\ShowEq{f A1 A2}
be continuous linear map
of $D$\Hyph module $A_1$ with norm $\|x\|_1$
into $D$\Hyph module $A_2$ with norm $\|y\|_2$.
Product of the map $f$
over scalar $d\in D$
is continuous linear map.
\end{theorem}
\begin{proof}
According to the corollary
\ref{corollary: linear map is continuous at 0}
and the definition
\ref{definition: continuous map, module},
for given $\epsilon>0$
there exist such $\delta>0$, that
\ShowEq{continuous map, x, 3}
implies
\ShowEq{continuous map, y, d}
From inequality
\EqRef{continuous map, y, d}
and the statement
\ref{item: norm on D module 3, 2},
it follows that
\ShowEq{continuous map, x, 3}
implies
\ShowEq{continuous map, y, 4}
Therefore,
according to the corollary
\ref{corollary: linear map is continuous at 0}
and the definition
\ref{definition: continuous map, module},
linear map $d\,f$ is continuous.
\end{proof}

\begin{theorem}
The set
\ShowEq{set continuous linear maps, module}
of continuous linear maps
of normed $D$\Hyph module $A_1$
into normed $D$\Hyph module $A_2$
is $D$\Hyph module.
\end{theorem}
\begin{proof}
The theorem follows from theorems
\ref{theorem: sum of continuous linear maps},
\ref{theorem: Product of continuous map over scalar}.
\end{proof}

\begin{theorem}
\label{theorem: norm of map, module}
Let $A_1$ be $D$\Hyph module with norm $\|x\|_1$.
Let $A_2$ be $D$\Hyph module with norm $\|y\|_2$.
The map
\ShowEq{LDA->R}
determined by the equation
\ShowEq{norm of map}
\ShowEq{norm of map, module}
is norm of $D$\Hyph module $\LDA$
and is called
\AddIndex{norm of map $f$}{norm of map}.
\end{theorem}
\begin{proof}
The statement
\ref{item: norm on D module 1}
is evident.

Let $\|f\|=0$.
According to the definition
\EqRef{norm of map, module}
\ShowEq{|f|=0 1}
for any $x\in A_1$.
According to the statement
\ref{item: norm on D module 2},
\ShowEq{|f|=0 2}
for any $x\in A_1$.
Therefore, the statement
\ref{item: norm on D module 2}
is true for $\|f\|$.

According to the definition
\EqRef{sum of linear maps, definition}
and the statement
\ref{item: norm on D module 3, 1},
\ShowEq{|f1+f2| 1}
From the inequality
\EqRef{|f1+f2| 1}
and the definition
\EqRef{norm of map, module},
it follows that
\ShowEq{|f1+f2| 2}
Therefore, the statement
\ref{item: norm on D module 3, 1}
is true for $\|f\|$.

According to the definition
\EqRef{product of map over scalar, definition}
and the statement
\ref{item: norm on D module 3, 2},
\ShowEq{|df| 1}
From the inequality
\EqRef{|df| 1}
and the definition
\EqRef{norm of map, module},
it follows that
\ShowEq{|df| 2}
Therefore, the statement
\ref{item: norm on D module 3, 2}
is true for $\|f\|$.
\end{proof}

\begin{theorem}
Let $D$ be ring with norm $|d|$.
Let $A$ be $D$\Hyph module with norm $\|x\|_0$.
The map
\ShowEq{A'->R}
determined by the equation
\ShowEq{norm of functional}
is norm of $D$\Hyph module $A'$
and is called
\AddIndex{norm of functional}{norm of functional} $f$.
\end{theorem}
\begin{proof}
The theorem follows from the theorem
\ref{theorem: norm of map, module}.
\end{proof}

\begin{theorem}
Let $D$ be ring with norm $|d|$.
Let $\Basis e$ be basis of $D$\Hyph module $A$ with norm $\|x\|_1$.
Let $A'$ be conjugated $D$\Hyph module with norm $\|x\|_2$.
Then
\ShowEq{norm of conjugate basis}
\end{theorem}
\begin{proof}
Let index $\gi i$ have given value.
Let $a\in A$.
Since
\ShowEq{a=a-+ai}
then according to the statement
\ref{item: norm on D module 3, 1}
\ShowEq{a=a-+ai 1}
If
\ShowEq{a=a-+ai 2}
then according to the statements
\ref{item: norm on D module 2},
\ref{item: norm on D module 3, 2}
\ShowEq{a=a-+ai 3}
Therefore,
\ShowEq{a=a-+ai 4}
The equation
\EqRef{norm of conjugate basis}
follows from the equation
\EqRef{a=a-+ai 4}.
\end{proof}

\begin{corollary}
\label{corollary: basis dual to normal basis is normal}
Let $\Basis e$ be normal basis of normed $D$\Hyph module $A$.
The basis dual to basis $\Basis e$
also is normal basis of normed $D$\Hyph module $A'$.
\qed
\end{corollary}

\begin{theorem}
Let $\Basis e_1$ be basis of $D$\Hyph module $A_1$ with norm $\|x\|_1$.
Let $\Basis e_2$ be basis of $D$\Hyph module $A_2$ with norm $\|x\|_2$.
Then
\ShowEq{norm of basis L(A1,A2)}
\end{theorem}
\begin{proof}
Let indices $\gi i$, $\gi j$ have given values.
Let $a\in A_1$.
Since
\ShowEq{a=a-+ai L(A1,A2)}
then according to the statement
\ref{item: norm on D module 3, 1}
\ShowEq{a=a-+ai 1 L(A1,A2)}
If
\ShowEq{a=a-+ai 2 L(A1,A2)}
then according to the statements
\ref{item: norm on D module 2},
\ref{item: norm on D module 3, 2}
\ShowEq{a=a-+ai 3 L(A1,A2)}
Since
\ShowEq{a=a-+ai 5}
then according to the statement
\ref{item: norm on D module 3, 2}
\ShowEq{a=a-+ai 6}
From equations
\EqRef{a=a-+ai 3 L(A1,A2)},
\EqRef{a=a-+ai 6}
it follows that
\ShowEq{a=a-+ai 7}
The equation
\EqRef{norm of basis L(A1,A2)}
follows from the equation
\EqRef{a=a-+ai 7}.
\end{proof}

\begin{corollary}
\label{corollary: norm of basis L(A1,A2)}
Let $\Basis e_1$
be normal basis of $D$\Hyph module $A_1$ with norm $\|x\|_1$.
Let $\Basis e_2$
be normal basis of $D$\Hyph module $A_2$ with norm $\|x\|_2$.
Then
\ShowEq{normal basis L(A1,A2)}
\qed
\end{corollary}

\begin{theorem}
\label{theorem: norm of linear map, module}
Let
\ShowEq{f A1 A2}
be linear map of $D$\Hyph module $A_1$ with norm $\|x\|_1$
into $D$\Hyph module $A_2$ with norm $\|y\|_2$.
Then
\ShowEq{norm of linear map, module}
\end{theorem}
\begin{proof}
From the definition
\ref{definition: linear map from A1 to A2, module}
and the theorems
\ref{theorem: complete ring contains real number},
\ref{theorem: complete ring and real number},
it follows that
\ShowEq{linear map and real number}
From the equation \EqRef{linear map and real number}
and the statement \ref{item: norm on D module 3, 2},
it follows that
\ShowEq{linear map and real number, 1}
Assuming $\displaystyle r=\frac 1{\|x\|_1}$, we get
\ShowEq{norm of linear map, module, 1}
Equation \EqRef{norm of linear map, module}
follows from equations \EqRef{norm of linear map, module, 1}
and \EqRef{norm of map, module}.
\end{proof}

\begin{theorem}
\label{theorem: |fx|<|f||x|}
Let
\ShowEq{f A1 A2}
be linear map of $D$\Hyph module $A_1$ with norm $\|x\|_1$
into $D$\Hyph module $A_2$ with norm $\|y\|_2$.
Then
\ShowEq{|fx|<|f||x|}
\end{theorem}
\begin{proof}
According to the statement
\ref{item: norm on D module 3, 2}
\ShowEq{x/|x|}
From the theorem
\ref{theorem: norm of linear map, module}
and the equation
\EqRef{x/|x|},
it follows that
\ShowEq{f x/|x|}
From the statement
\ref{item: norm on D module 3, 2}
and the equation
\EqRef{f x/|x|},
it follows that
\ShowEq{f x/|x| 1}
The inequality
\EqRef{|fx|<|f||x|}
follows from the inequality
\EqRef{f x/|x| 1}.
\end{proof}

\begin{theorem}
\label{theorem: map is continuous iff}
Let\footnote{This theorem
is based on the theorem
\citeBib{Kolmogorov Fomin}-2, pages 77 - 78.}
\ShowEq{f A1 A2}
linear map of $D$\Hyph module $A_1$ with norm $\|x\|_1$
into $D$\Hyph module $A_2$ with norm $\|y\|_2$.
The map $f$ is continuous iff
\ShowEq{|f|<infty}
\end{theorem}
\begin{proof}
Let
\ShowEq{|f|<infty}
Since map $f$ is linear, then
according to the theorem \ref{theorem: norm of map, module}
\ShowEq{|f|<infty 1}
Let us assume arbitrary $\epsilon>0$. Assume $\displaystyle\delta=\frac\epsilon{\|f\|}$.
Then
\ShowEq{|f|<infty 3}
follows from inequality
\ShowEq{|f|<infty 2}
According to definition \ref{definition: continuous map, module},
the map $f$ is continuous.

Let
\ShowEq{|f|=infty}
According to the theorem
\ref{theorem: norm of map, module},
for any $n$, there exists $x_n$ such that
\ShowEq{|f(x)|>n}
Let
\ShowEq{xn=}
According to the definition
\ref{definition: linear map from A1 to A2, algebra},
the statement \ref{item: norm on D module 3, 2},
equation
\EqRef{xn=},
inequality
\EqRef{|f(x)|>n}
\ShowEq{|f(y)|=1}
Therefore, the map $f$ is not continuous at the point $0\in A_1$.
\end{proof}

$D$\Hyph module $\LcDA$ is submodule of
$D$\Hyph module $\LDA$.
According to the theorem
\ref{theorem: map is continuous iff},
since
\ShowEq{LcDA-LDA}
then
\ShowEq{|f|=infty}

\ifx\texFuture\Defined
\begin{theorem}
Let $A_1$ be $D$\Hyph module with norm $\|x\|_1$.
Let $A_2$ be Banach $D$\Hyph module with norm $\|y\|_2$.
$D$\Hyph module $\LcDA$ is Banach $D$\Hyph module.
\end{theorem}
\begin{proof}
\end{proof}
\fi

\begin{theorem}
\label{theorem: norm of product of maps}
Let $A_1$ be $D$\Hyph module with norm $\|x\|_1$.
Let $A_2$ be $D$\Hyph module with norm $\|x\|_2$.
Let $A_3$ be $D$\Hyph module with norm $\|x\|_3$.
Let
\ShowEq{g A1 A2}
\ShowEq{f A2 A3}
be continuous linear maps.
The map
\ShowEq{fg}
is continuous linear map
\ShowEq{|fg|}
\end{theorem}
\begin{proof}
According to the definitions
\EqRef{(fg)a=f(ga)},
\EqRef{norm of map, module}\footnote{Let
\ShowEq{fi A1 R}
In general, maps $f_1$, $f_2$ have maximum
in different points of the set $A_1$.
Therefore,
\ShowEq{fi A1 R 1}
}
\ShowEq{|fg| 1}
Since, in general,
\ShowEq{g A1}
then
\ShowEq{|fg| 0}
From the inequalities
\EqRef{|fg| 1},
\EqRef{|fg| 0},
it follows that
\ShowEq{|fg| 2}
The inequality
\EqRef{|fg|}
follows from the inequality
\EqRef{|fg| 2}
and the definition
\EqRef{norm of map, module}.
\end{proof}

\begin{theorem}
Normed $D$\Hyph module $\mathcal{LC}(D;A;A)$
is normed $D$\Hyph algebra,
where product is defined according to rule
\ShowEq{f g ->}
\end{theorem}
\begin{proof}
The proof of statement that $D$\Hyph module $\mathcal{LC}(D;A;A)$
is $D$\Hyph algebra is similar to the proof of the theorem
\xRef{8433-5163}{theorem: module L(A;A) is algebra}.
According to the definition
\ref{definition: norm on d algebra}
and the theorem
\ref{theorem: norm of product of maps},
the norm
\EqRef{norm of map, module}
is norm on $D$\Hyph algebra $\mathcal{LC}(D;A;A)$.
\end{proof}

\section{Normed \texorpdfstring{$D$\Hyph}{D-}Module
\texorpdfstring{$\mathcal L(D;A_1,...,A_n;A)$}{L(D;A1,...,An;A)}}
\label{section: Module L(D;A1,...,An;A)}

\begin{definition}
\label{definition: continuous multivariable map, module}
Let
\ShowEq{Ai 1n}
be Banach $D$\Hyph algebra with norm $\|x\|_i$.
Let $A$
be Banach $D$\Hyph algebra with norm $\|x\|$.
Multivariable map
\ShowEq{A1...An->A}
is called \AddIndex{continuous}{continuous multivariable map}, if
for every as small as we please $\epsilon>0$
there exist such $\delta>0$, that
\DrawEq{|xi|<e, 1-n}{}
implies
\DrawEq{continuous multivariable map, y}{}
\qed
\end{definition}

\begin{theorem}
\label{theorem: sum of continuous multivariable maps}
The sum of continuous multivariable maps
is continuous multivariable map.
\end{theorem}
\begin{proof}
Let
\ShowEq{A1...An->A}
be continuous multivariable map.
According to the definition
\ref{definition: continuous multivariable map, module},
for given $\epsilon>0$
there exist such $\delta_1>0$, that
\ShowEq{continuous multivariable map, x, 1}
implies
\ShowEq{continuous multivariable map, y, 1}
Let
\ShowEq{g A1...An->A}
be continuous multivariable map.
According to the definition
\ref{definition: continuous multivariable map, module},
for given $\epsilon>0$
there exist such $\delta_2>0$, that
\ShowEq{continuous multivariable map, x, 2}
implies
\ShowEq{continuous multivariable map, y, 2}
Let
\ShowEq{delta 1 2}
From inequalities
\EqRef{continuous multivariable map, y, 1},
\EqRef{continuous multivariable map, y, 2}
and the statement
\ref{item: norm on D module 3, 1},
it follows that
\ShowEq{continuous multivariable map, x, 3}
implies
\ShowEq{continuous multivariable map, y, 3}
Therefore,
according to the definition
\ref{definition: continuous multivariable map, module},
multivariable map $f+g$ is continuous.
\end{proof}

\begin{theorem}
\label{theorem: sum of continuous polylinear maps}
The sum of continuous polylinear maps
is continuous polylinear map.
\end{theorem}
\begin{proof}
The theorem follows from theorems
\xRef{1011.3102}{theorem: sum of polylinear maps, algebra},
\ref{theorem: sum of continuous multivariable maps}.
\end{proof}

\begin{theorem}
\label{theorem: Product of continuous multivariable map over scalar}
Product of the continuous multivariable map $f$
over scalar $d\in D$
is continuous multivariable map.
\end{theorem}
\begin{proof}
Let
\ShowEq{A1...An->A}
be continuous multivariable map.
According to the definition
\ref{definition: continuous multivariable map, module},
for given $\epsilon>0$
there exist such $\delta>0$, that
\ShowEq{continuous multivariable map, x, 3}
implies
\ShowEq{continuous multivariable map, y, d}
From inequality
\EqRef{continuous multivariable map, y, d}
and the statement
\ref{item: norm on D module 3, 2},
it follows that
\ShowEq{continuous multivariable map, x, 3}
implies
\ShowEq{continuous multivariable map, y, 4}
Therefore,
according to the definition
\ref{definition: continuous multivariable map, module},
multivariable map $d\,f$ is continuous.
\end{proof}

\begin{theorem}
\label{theorem: Product of continuous polylinear map over scalar}
Product of the continuous polylinear map $f$
over scalar $d\in D$
is continuous polylinear map.
\end{theorem}
\begin{proof}
The theorem follows from theorems
\xRef{1011.3102}{theorem: polylinear map times scalar, algebra},
\ref{theorem: Product of continuous multivariable map over scalar}.
\end{proof}

\begin{theorem}
The set
\ShowEq{set continuous multivariable maps, module}
of continuous multivariable maps
is $D$\Hyph module.
\end{theorem}
\begin{proof}
The theorem follows from theorems
\ref{theorem: sum of continuous multivariable maps},
\ref{theorem: Product of continuous multivariable map over scalar}.
\end{proof}

\begin{theorem}
The set
\ShowEq{set continuous polylinear maps, module}
of continuous polylinear maps
is $D$\Hyph module.
\end{theorem}
\begin{proof}
The theorem follows from theorems
\ref{theorem: sum of continuous polylinear maps},
\ref{theorem: Product of continuous polylinear map over scalar}.
\end{proof}

Let $A_1$ be $D$\Hyph module with norm $\|x\|_1$.
Let $A_2$ be $D$\Hyph module with norm $\|x\|_2$.
Let $A_3$ be $D$\Hyph module with norm $\|x\|_3$.
Since $\mathcal L(D;A_2;A_3)$ is $D$\Hyph module with norm
\ShowEq{|f|23}
then we can consider continuous map
\ShowEq{h A1 A2->A3}
Since $a_1\in A_1$, then
\ShowEq{h a_1}
is continuous map.
According to the theorem
\ref{theorem: |fx|<|f||x|}
\ShowEq{|a3|=|ha1||a2|}
Since
\ShowEq{L(A1,L(A2,A3))}
is normed $D$\Hyph module,
then according to the theorem
\ref{theorem: |fx|<|f||x|}
\ShowEq{|ha1|}
From inequalities
\EqRef{|a3|=|ha1||a2|},
\EqRef{|ha1|},
it follows that
\ShowEq{a3 h a1 a2}
We can consider the map
\EqRef{h A1 A2->A3}
as bilinear map
\ShowEq{f A1,A2->A3}
defined by rule
\ShowEq{f a1,a2}
Based on theorems
\ref{theorem: norm of map, module},
\ref{theorem: |fx|<|f||x|}
and inequality
\EqRef{a3 h a1 a2},
we define the norm of the bilinear map $f$ by the equation
\ShowEq{|f| a1a2}

If we proceed by induction over number of variables,
we can generalize the definition
of norm of bilinear map.

\begin{definition}
\label{definition: norm of polylinear map}
Let
\ShowEq{Ai 1n}
be Banach $D$\Hyph algebra with norm $\|x\|_i$.
Let $A$
be Banach $D$\Hyph algebra with norm $\|x\|$.
Let
\ShowEq{f A1...An->A}
be polylinear map.
Value
\ShowEq{norm of polymap, module}
is called
\AddIndex{norm of polylinear map $f$}
{norm of polymap}.
\qed
\end{definition}

\begin{theorem}
\label{theorem: norm of polylinear map, module}
Let
\ShowEq{Ai 1n}
be Banach $D$\Hyph module with norm $\|x\|_i$.
Let $A$
be Banach $D$\Hyph module with norm $\|x\|$.
Let
\ShowEq{f A1...An->A}
be polylinear map.
Then
\ShowEq{norm of polylinear map, module}
\end{theorem}
\begin{proof}
From the definition
\ref{definition: polylinear map of modules}
and the theorems
\ref{theorem: complete ring contains real number},
\ref{theorem: complete ring and real number},
it follows that
\ShowEq{polylinear map and real number}
From the equation \EqRef{polylinear map and real number}
and the statement \ref{item: norm on D module 3, 2},
it follows that
\ShowEq{polylinear map and real number, 1}
Assuming $\displaystyle r=\frac 1{\|x\|_1}$, we get
\ShowEq{norm of polylinear map, module, 1}
Equation \EqRef{norm of polylinear map, module}
follows from equations \EqRef{norm of polylinear map, module, 1}
and \EqRef{norm of polymap, module}.
\end{proof}

\begin{theorem}
\label{theorem: |fx|<|f||x|1n}
Let
\ShowEq{Ai 1n}
be Banach $D$\Hyph module with norm $\|x\|_i$.
Let $A$
be Banach $D$\Hyph module with norm $\|x\|$.
Let
\ShowEq{f A1...An->A}
be polylinear map.
Then
\DrawEq{|fx|<|f||x|1n}{theorem}
\end{theorem}
\begin{proof}
According to the statement
\ref{item: norm on D module 3, 2}
\ShowEq{x/|x|1n}
From the theorem
\ref{theorem: norm of polylinear map, module}
and the equation
\EqRef{x/|x|1n},
it follows that
\ShowEq{f x/|x|1n}
From the statement
\ref{item: norm on D module 3, 2}
and the equation
\EqRef{f x/|x|1n},
it follows that
\ShowEq{f x/|x|1n 1}
The inequality
\eqRef{|fx|<|f||x|1n}{theorem}
follows from the inequality
\EqRef{f x/|x|1n 1}.
\end{proof}

Let
\ShowEq{Ai 1n}
be Banach $D$\Hyph module with norm $\|x\|_i$.
Let $A$
be Banach $D$\Hyph module with norm $\|x\|$.
We can represent polylinear map
\ShowEq{f A1...An->A}
in the following form
\DrawEq{f1n=hx}{convention}
where
\ShowEq{h 1n-1->A}
is polylinear map.

\begin{theorem}
\label{theorem: polylinear map f h()}
Since the map $f$ is continues, then the map
\ShowEq{polymap 3, algebra}
is also continues.
\end{theorem}
\begin{proof}
According to the definition
\ref{definition: continuous multivariable map, module},
for every as small as we please $\epsilon>0$
there exist such $\delta>0$, that
\DrawEq{|xi|<e, 1-n}{}
implies
\ShowEq{polymap 30, y}
Therefore,
for every as small as we please $\epsilon>0$
there exist such $\delta>0$, that
\ShowEq{polymap 3, x}
implies
\ShowEq{polymap 3, y}
From the equation
\eqRef{f1n=hx}{convention}
and inequality
\EqRef{polymap 3, y}
it follows that
for every as small as we please $\epsilon>0$
there exist such $\delta>0$, that
\ShowEq{polymap 3, x}
implies
\ShowEq{polymap 31, y}
According to the definition
\ref{definition: continuous map, module},
the map
\ShowEq{polymap 32, y}
is continuous.
\end{proof}

\begin{theorem}
\label{theorem: polylinear map f h}
Since the map $f$ is continues, then the map
$h$ is also continues.
\end{theorem}
\begin{proof}
According to the definition
\ref{definition: continuous multivariable map, module},
for every as small as we please $\epsilon>0$
there exist such $\delta>0$, that
\DrawEq{|xi|<e, 1-n}{}
implies
\ShowEq{polymap 30, y}
Therefore,
for every as small as we please $\epsilon>0$
there exist such $\delta>0$, that
\ShowEq{polymap 4, x}
implies
\ShowEq{polymap 4, y}
From the equation
\eqRef{f1n=hx}{convention}
and inequality
\EqRef{polymap 4, y}
it follows that
for every as small as we please $\epsilon>0$
there exist such $\delta>0$, that
\ShowEq{polymap 4, x}
implies
\ShowEq{polymap 41, y}
From the inequality
\EqRef{polymap 41, y}
it follows that
for every as small as we please $\epsilon>0$
there exist such $\delta>0$, that
\ShowEq{polymap 4, x}
implies
\ShowEq{polymap 42, y}
According to the definition
\EqRef{norm of map, module}
\ShowEq{polymap 43, y}
From inequalities
\EqRef{polymap 42, y},
\EqRef{polymap 43, y},
it follows that
for every as small as we please $\epsilon>0$
there exist such $\delta>0$, that
\ShowEq{polymap 4, x}
implies
\ShowEq{polymap 44, y}
According to the definition
\ref{definition: continuous map, module},
the map
$h$
is continuous.
\end{proof}

\begin{theorem}
\label{theorem: polylinear map is continuous, submaps}
\ShowEq{polymap 5, f}
iff
\ShowEq{polymap 5, h}
\end{theorem}
\begin{remark}
In other words,
the polylinear map $f$ is continuous iff
the map $h$ is continuous and for any
\ShowEq{polymap 5, a}
the map
\ShowEq{polymap 3, algebra}
is continuous.
\end{remark}
\begin{proof}
From theorems
\ref{theorem: polylinear map f h()},
\ref{theorem: polylinear map f h},
it follows that continuity of maps $h$ and
\ShowEq{polymap 3, algebra}
follows from continuity of the map $f$.

Let maps $h$ and
\ShowEq{polymap 3, algebra}
be continuous.
According to the definition
\ref{definition: continuous multivariable map, module},
to prove continuity of the map $f$,
we need to estimate the difference
\ShowEq{polymap 5, y}
provided that
\DrawEq{|xi|<e, 1-n}{4}
According to the equation
\eqRef{f1n=hx}{convention},
\ShowEq{polymap 51, y}
According to the equation
\EqRef{polymap 51, y}
and statement
\ref{item: norm on D module 3, 2},
\ShowEq{polymap 52, y}
According to the definition
\ref{definition: continuous map, module},
for every as small as we please $\epsilon>0$
there exist such $\delta_1>0$, that
\ShowEq{|xn|<d1}
implies
\ShowEq{h(x)dx<e/2}
Consider second term in right part of the inequality
\EqRef{polymap 52, y}.
\StartLabelItem{theorem: polylinear map is continuous, submaps}
\begin{enumerate}
\item Since $x_n=0$, then
\ShowEq{dh(x)x<e/2 1}
\label{item: dh(x)x<e/2 1}
\item So we assume $x_n\ne 0$.
According to the definition
\ref{definition: continuous multivariable map, module},
for every as small as we please $\epsilon>0$
there exist such $\delta_2>0$, that
\ShowEq{|xi|<d2}
implies
\ShowEq{dh(x)<e/2}
From inequalities
\EqRef{dh(x)<e/2},
\EqRef{|fx|<|f||x|},
it follows that
\ShowEq{dh(x)x<e/2 2}
\label{item: dh(x)x<e/2 2}
\end{enumerate}
Therefore, in both cases
\RefItem{dh(x)x<e/2 1},
\RefItem{dh(x)x<e/2 2},
for every as small as we please $\epsilon>0$
there exist such $\delta_2>0$, that
\ShowEq{|xi|<d2}
implies
\ShowEq{dh(x)x<e/2 3}
Let
\ShowEq{delta 1 2}
From inequalities
\EqRef{polymap 52, y},
\EqRef{h(x)dx<e/2},
\EqRef{dh(x)x<e/2 3},
it follows that
for every as small as we please $\epsilon>0$
there exist such $\delta_2>0$, that
\DrawEq{|xi|<e, 1-n}{}
implies
\ShowEq{polymap 30, y}
According to the definition
\ref{definition: continuous multivariable map, module},
the map $f$ is continues.
\end{proof}

\begin{theorem}
\label{theorem: |f|=|h|}
Let
\ShowEq{Ai 1n}
be Banach $D$\Hyph module with norm $\|x\|_i$.
Let $A$ be Banach $D$\Hyph module with norm $\|x\|$.
Let
\ShowEq{f A1...An->A}
be polylinear map.
Let
\ShowEq{h 1n-1->A}
be polylinear map such that
\DrawEq{f1n=hx}{|f|=|h|}
Then
\ShowEq{|f|=|h|}
\end{theorem}
\begin{proof}
According to the definition
\ref{definition: norm of polylinear map}
\ShowEq{|hx1n-1|}
According to the theorem
\ref{theorem: |fx|<|f||x|}
\ShowEq{|x|=|hx1n-1||xn|}
From inequalities
\EqRef{|hx1n-1|},
\EqRef{|x|=|hx1n-1||xn|},
it follows that
\ShowEq{x h x1 xn}
According to the theorem
\ref{theorem: |fx|<|f||x|1n}
\DrawEq{|fx|<|f||x|1n}{|f|=|h|}
Based on theorem
\ref{theorem: |fx|<|f||x|1n}
and inequalities
\EqRef{x h x1 xn},
\eqRef{|fx|<|f||x|1n}{|f|=|h|},
we get inequality
\ShowEq{|f|<|h|}
From the equation
\eqRef{f1n=hx}{|f|=|h|}
and inequality
\eqRef{|fx|<|f||x|1n}{|f|=|h|},
it follows that
\ShowEq{|hx|/|xn|}
From the definition
\ref{definition: norm of polylinear map}
and inequality
\EqRef{|hx|/|xn|},
it follows that
\ShowEq{|hx|<|f||x|}
Based on theorem
\ref{theorem: |fx|<|f||x|1n}
and inequalities
\EqRef{|hx1n-1|},
\EqRef{|hx|<|f||x|},
we get inequality
\ShowEq{|h|<|f|}
The equation
\EqRef{|f|=|h|}
follows from inequalities
\EqRef{|f|<|h|},
\EqRef{|h|<|f|}.
\end{proof}

\begin{theorem}
The polylinear map $f$ is continues iff
$\|f\|<\infty$.
\end{theorem}
\begin{proof}
We prove the theorem by induction over number $n$ of arguments of the map $f$.
For $n=1$, the theorem follows from the theorem
\ref{theorem: map is continuous iff}.

Let the theorem be true for $n=k-1$.
Let
\ShowEq{A1k}
be Banach $D$\Hyph module with norm $\|x\|_i$.
Let $A$ be Banach $D$\Hyph module with norm $\|x\|$.
We can represent polylinear map
\ShowEq{f A1...Ak->A}
in the following form
\DrawEq{f1k=ha}{}
where
\ShowEq{h 1k-1->A}
is polylinear map.
According to the theorem
\ref{theorem: polylinear map f h},
the map $h$ is continuous polylinear map of $k-1$ variables.
According to the induction assumption,
\ShowEq{|h|<infty}
According to the theorem
\ref{theorem: |f|=|h|},
\ShowEq{|f|=|h|<infty}
\end{proof}

\section{\texorpdfstring{$D$}{D}\hyph algebra with Schauder basis}

\begin{definition}
\label{definition: Schauder basis}
Let $A$ be Banach $D$\Hyph module.\footnote{The definition
\ref{definition: Schauder basis}
is based on the definition
\citeBib{4419-7514}-4.6, p. 182,
and lemma
\citeBib{4419-7514}-4.7, p. 183.}
A sequence of vectors
\ShowEq{Schauder basis}
is called
\AddIndex{Schauder basis}{Schauder basis}
if
\begin{itemize}
\item The set of vectors $e_{\gi i}$ is linear independent.
\item For each vector $a\in A$, there exists
unique sequence
\ShowEq{Schauder basis, coordinates}
such that
\DrawEq{Schauder basis, coordinates 0}{}
\end{itemize}
The sequence 
\ShowEq{ai in D}
is called
\AddIndex{coordinates of vector
\ShowEq{a=ai ei}
relative to Schauder basis}
{coordinates of vector, Schauder basis}
$\Basis e$.
\qed
\end{definition}

Let $\Basis e$ be Schauder basis of Banach $D$\Hyph module $A$.
We say that the {\bf expansion}
\ShowEq{a=ae}
\AddIndex{of vector $a\in A$ relative to the basis $\Basis e$ converges}
{expansion converges}.

\begin{theorem}
Let $\Basis e$ be Schauder basis of Banach $D$\Hyph module $A$.
Then
\ShowEq{Schauder basis, norm is finite, 1}
for any vector $a\in A$.
\end{theorem}
\begin{proof}
From the theorem
\ref{theorem: normed module, norm of difference},
it follows that we cannot define
\ShowEq{Schauder basis, norm is finite, 3}
if
\ShowEq{Schauder basis, norm is finite, 2}.
Therefore, we cannot expand $a$ relative to Schauder basis.
\end{proof}

\begin{theorem}
\label{theorem: Schauder basis, limit}
Let $\Basis e$ be Schauder basis of Banach $D$\Hyph module $A$.
Let $a_i$ be coordinates of vector $a$ relative to the basis $\Basis e$.
Then for every
\ShowEq{epsilon in R}
there exists positive integer $n_0$ depending on $\epsilon$ and such,
that
\ShowEq{Schauder basis, limit}
for every $p$, $q>n_0$.
\end{theorem}
\begin{proof}
The inequality
\EqRef{Schauder basis, limit, p}
follows from definitions
\ref{definition: limit of sequence, normed module},
\ref{definition: Schauder basis}.
The inequality
\EqRef{Schauder basis, limit, p, q}
follows from definitions
\ref{definition: fundamental sequence, normed module},
\ref{definition: Schauder basis}.
\end{proof}

\begin{theorem}
Let
\ShowEq{ai in I}
finite set of vectors of Banach $D$\Hyph module $A$
with Schauder basis $\Basis e$.
Then\footnote{See the definition
\xRef{0701.238}{definition: linear span, vector space}
of linear span in vector space.}
\ShowEq{span in A}
\end{theorem}
\begin{proof}
To prove the statement
\EqRef{span in A},
it suffices to prove the following statements.
\begin{itemize}
\item
If
\ShowEq{a1 a2 A}
then $a_1+a_2\in A$.

Since 
\ShowEq{a1 a2 A}
then, according to the theorem
\ref{theorem: Schauder basis, limit},
there exists positive integer $n_0$ depending on $\epsilon$ and such,
that
\ShowEq{Schauder basis, limit p12}
for every $p$, $q>n_0$.
From inequalities
\EqRef{Schauder basis, limit p12}
it follows that
\ShowEq{Schauder basis, limit p12 1}
From inequality
\EqRef{Schauder basis, limit p12 1}
it follows that sequence
\ShowEq{Schauder basis, limit p12 2}
is fundamental sequence.
\item
If
\ShowEq{aA dD}
then $d\,a\in A$.

Since 
\ShowEq{aA dD}
then, according to the theorem
\ref{theorem: Schauder basis, limit},
there exists positive integer $n_0$ depending on $\epsilon$ and such,
that
\ShowEq{Schauder basis, limit p22}
for every $p$, $q>n_0$.
From inequality
\EqRef{Schauder basis, limit p22}
it follows that
\ShowEq{Schauder basis, limit p22 1}
From inequality
\EqRef{Schauder basis, limit p22 1}
it follows that sequence
\ShowEq{Schauder basis, limit p22 2}
is fundamental sequence.
\end{itemize}
\end{proof}

\begin{theorem}
\label{theorem: Schauder basis, norm is finite}
Let $\Basis e$ be Schauder basis of Banach $D$\Hyph module $A$.
Then
\ShowEq{Schauder basis, norm is finite}
for any vector $e_{\gi i}$.
\end{theorem}
\begin{proof}
Let for
\ShowEq{Schauder basis, norm is finite, i=j}
For any
\ShowEq{Schauder basis, norm is finite, n>j}
if for sequence
\ShowEq{ai in D}
it is true that
\ShowEq{Schauder basis, norm is finite, j}
Therefore, we cannot say whether the vector
\ShowEq{a=ai ei} is defined.
\end{proof}

Without loss of generality, we can assume that the basis $\Basis e$
is normal basis. If we assume that the norm of the vector
$e_{\gi i}$ is different from $1$, then we can replace this vector
by the vector
\ShowEq{Banach algebra, normal basis, 1}
According to the corollary
\ref{corollary: basis dual to normal basis is normal},
dual basis also is normal basis.

\begin{theorem}
\label{theorem: coordinates converges absolutely}
Let $\Basis e$ be
normal Schauder basis of Banach $D$\Hyph module $A$.
Let
\ShowEq{ai in D}
be such sequence that
\ShowEq{ai in D 1}
Then there exists limit\footnote{See
similar theorems
\citeBib{Shilov single 3}, page 60,
\citeBib{Fihtengolts: Calculus volume 2}, pages 264, 295,
\citeBib{Smirnov vol 1}, pages 319, 329.}
\DrawEq{Schauder basis, coordinates 0}{1}
\end{theorem}
\begin{proof}
Existence of the limit
\eqRef{Schauder basis, coordinates 0}{1}
follows from the inequality
\ShowEq{Schauder basis, coordinates 2}
since the inequality
\EqRef{Schauder basis, coordinates 2}
is preserved in passing to the limit
\ShowEq{Schauder basis, coordinates 3}
\end{proof}

Let $\Basis e$ be normal Schauder basis of Banach $D$\Hyph module $A$.
If
\ShowEq{ai in D 1}
then we say that the {\bf expansion}
\ShowEq{a=ae}
\AddIndex{of vector $a\in A$ relative to the basis $\Basis e$
converges normally}
{expansion converges normally}.\footnote{The definition of
normal convergence of the expansion of vector relative to basis
is similar to the definition of normal convergence of the series.
See \citeBib{Cartan differential form}, page 12.}
We denote
\ShowEq{A plus Schauder}
the set of vectors
whose expansion relative to the basis $\Basis e$
converges normally.

\begin{theorem}
\label{theorem: |a|<|ai|}
Let $\Basis e$ be normal Schauder basis of Banach $D$\Hyph module $A$.
If expansion
of vector $a\in A$ relative to the basis $\Basis e$
converges normally, then
\ShowEq{|a|<|ai|}
\end{theorem}
\begin{proof}
From the statement
\RefItem{norm on D module 3, 1},
it follows that
\ShowEq{|a|<|ai| 1}
The inequation
\EqRef{|a|<|ai|}
follows from the inequation
\EqRef{|a|<|ai| 1}
and the statement that Schauder basis $\Basis e$ is a normal basis.
\end{proof}


\begin{theorem}
\label{theorem: linear map, module, Schauder basis}
Let
\ShowEq{f:A1->A2}
be map of $D$\Hyph module $A_1$ with basis $\Basis e_1$ 
into $D$\Hyph module $A_2$ with Schauder basis $\Basis e_2$.
Let $f^{\gi i}_{\gi j}$ be coordinates of the map $f$
relative to bases $\Basis e_1$ and $\Basis e_2$.
Then sequence
\ShowEq{linear map, algebra, Schauder basis}
has limit for any $\gi j$.
\end{theorem}
\begin{proof}
The theorem follows from the equation
\ShowEq{f:A1->A2, 1}
\end{proof}

\begin{theorem}
\label{theorem: |fe|<|f|}
Let
\ShowEq{f:A1->A2}
be map of $D$\Hyph module $A_1$ with norm $\|x\|_1$
and normal basis $\Basis e_1$ 
into $D$\Hyph module $A_2$ with norm $\|y\|_2$
and Schauder basis $\Basis e_2$.
Then
\ShowEq{|fe|<|f|}
for any $\gi i$.
\end{theorem}
\begin{proof}
According to the theorem
\ref{theorem: norm of map, module}
and the definition
\ref{definition: normal basis},
the inequality
\EqRef{|fe|<|f|}
follows from the inequality
\ShowEq{|fe|<|f|, 1}
\end{proof}

\begin{remark}
\label{remark: module L(A1,A2) has Schauder basis}
The theorem
\ref{theorem: linear map, module, Schauder basis}
determines constraint on coordinates of map of
$D$\Hyph module with Schauder basis.
However, we can make this constraint more stronger.
Let $A_1$ be $D$\Hyph module with normal Schauder basis $\Basis e_1$.
Let $A_2$ be $D$\Hyph module with normal Schauder basis $\Basis e_2$.
According to the theorem
\ref{theorem: basis of module A1 A2},
$D$\Hyph module $\LDA$ has basis
\ShowEq{L(A1,A2) Schauder basis}
Since the basis of $D$\Hyph module $\LDA$
is countable basis and $D$\Hyph module $\LDA$
has norm, then we require considered basis
to be Schauder basis.
According to the definition
\ref{definition: Schauder basis},
there exists limit
\ShowEq{map in Schauder basis}
The existence of the limit
\EqRef{map in Schauder basis}
implies existence of limit of sequence
\EqRef{linear map, algebra, Schauder basis}.
However, the existence of the limit
\EqRef{map in Schauder basis}
also implies that
\ShowEq{map in Schauder basis 1}
\qed
\end{remark}

\begin{theorem}
\label{theorem: |f|<epsilon}
Let
\ShowEq{f:A1->A2}
be linear map of $D$\Hyph module $A_1$ with norm $\|x\|_1$
and normal Schauder basis $\Basis e_1$
into $D$\Hyph module $A_2$ with norm $\|y\|_2$
and normal Schauder basis $\Basis e_2$.
For any
\ShowEq{epsilon in R}
there exist
$\gi N$, $\gi M$
such that
\DrawEq{|f|<epsilon}{theorem}
\end{theorem}
\begin{proof}
According to the remark
\ref{remark: module L(A1,A2) has Schauder basis},
the set of maps
\ShowEq{basis of L(A1,A2)}
is Schauder basis of $D$\Hyph module $\LDA$.
Therefore, expansion
\ShowEq{fa=,5}
of map $f$ convergies.
According to the theorem
\ref{theorem: Schauder basis, limit},
for any
\ShowEq{epsilon in R}
there exist
$\gi N$, $\gi M$
such that
\ShowEq{|f|<epsilon 1}
According to the corollary
\ref{corollary: norm of basis L(A1,A2)},
inequation
\eqRef{|f|<epsilon}{theorem}
follows from the inequation
\EqRef{|f|<epsilon 1}.
\end{proof}

\begin{theorem}
\label{theorem: |f|<F}
Let
\ShowEq{f:A1->A2}
be linear map of $D$\Hyph module $A_1$ with norm $\|x\|_1$
and normal Schauder basis $\Basis e_1$
into $D$\Hyph module $A_2$ with norm $\|y\|_2$
and normal Schauder basis $\Basis e_2$.
There exists
\ShowEq{F<infty}
such that
\ShowEq{|f|<<F}
\end{theorem}
\begin{proof}
According to the theorem
\ref{theorem: |f|<epsilon},
for given
\ShowEq{epsilon in R}
there exist
$\gi N$, $\gi M$
such that
\DrawEq{|f|<epsilon}{}
Since $\gi N$, $\gi M$ are finite,
then there exists
\ShowEq{F1=max f}
We get the inequation
\EqRef{|f|<<F},
if assume
\ShowEq{F F1 epsilon}
\end{proof}

\begin{theorem}
\label{theorem: Schauder basis, image of map}
Let
\ShowEq{f:A1->A2}
be linear map of $D$\Hyph module $A_1$ with norm $\|x\|_1$
and normal Schauder basis $\Basis e_1$
into $D$\Hyph module $A_2$ with norm $\|y\|_2$
and Schauder basis $\Basis e_2$.
Let
\ShowEq{|f|<infty}
Then for any
\ShowEq{a1 in A1 plus}
the image
\DrawEq{f:A1->A2, 2}{Schauder}
is defined properly,
\ShowEq{a2 in A2 plus}
\end{theorem}
\begin{proof}
\ifx\texFuture\Defined
Let
\ShowEq{a1 in A1}
According to the theorem
\ref{theorem: Schauder basis, limit},
for given
\ShowEq{epsilon in R}
there exists positive integer $n_0$ depending on $\epsilon$ and such,
that
\ShowEq{a1 in A1, 1}
for every $p$, $q>n_0$.
\fi
From the equation
\EqRef{a1 in A1 plus}
and the theorem
\ref{theorem: |fe|<|f|},
it follows that
\ShowEq{|a2|<infty}
From inequality
\EqRef{|a2|<infty}
it follows that
\ShowEq{|a2|<infty 1}
According to the theorem
\ref{theorem: coordinates converges absolutely},
image of $a_1\in A_1$ under the map $f$ is defined properly.
\end{proof}

\begin{remark}
From the proof of the theorem
\ref{theorem: Schauder basis, image of map},
we see that the requirement of normal convergence
of expansion of vector relative to normal basis
is essential.
According to the remark
\ref{remark: module L(A1,A2) has Schauder basis},
if $A_i$, $i=1$, $2$,
is $D$\Hyph module with normal Schauder basis $\Basis e_i$,
then the set $\mathcal L(D;A_1;A_2)$ is
$D$\Hyph module with normal Schauder basis
\ShowEq{Basis e1 e2}.
We denote
\ShowEq{L plus A1 A2}
the set of linear maps whose expansion relative to the basis
\ShowEq{Basis e1 e2}
converges normally.
\qed
\end{remark}

\begin{theorem}
Let $A_1$
be $D$\Hyph module with norm $\|x\|_1$
and normal Schauder basis $\Basis e_1$.
Let $A_2$
be $D$\Hyph module with norm $\|x\|_2$
and normal Schauder basis $\Basis e_2$.
Let a map
\ShowEq{f+A1->A2}
Then
\ShowEq{|f|<ij|f|}
\end{theorem}
\begin{proof}
According to the corollary
\ref{corollary: norm of basis L(A1,A2)},
the basis
\ShowEq{Basis e1 e2}
is normal Schauder basis.
The theorem follows from the theorem
\ref{theorem: |a|<|ai|}.
\end{proof}

\begin{corollary}
\label{corollary: |f|<infty}
Let $A_1$
be $D$\Hyph module with norm $\|x\|_1$
and normal Schauder basis $\Basis e_1$.
Let $A_2$
be $D$\Hyph module with norm $\|x\|_2$
and normal Schauder basis $\Basis e_2$.
Let a map
\ShowEq{f+A1->A2}
Then
\ShowEq{|f|<infty}
\qed
\end{corollary}

\begin{theorem}
\label{theorem: s|ab|<s|a|s|b|}
Let
\ShowEq{a,b>0}
Then
\ShowEq{s|ab|<s|a|s|b|}
\end{theorem}
\begin{proof}
We prove the theorem by induction over $n$.

The inequation
\EqRef{s|ab|<s|a|s|b|}
for $n=2$ follows from the inequation
\ShowEq{s|ab|<s|a|s|b| n=2}

Let the inequation
\EqRef{s|ab|<s|a|s|b|}
is true for $n=k-1$.
The inequation
\ShowEq{s|ab|<s|a|s|b| n=k}
follows from the inequation
\EqRef{s|ab|<s|a|s|b| n=2}.
From the inequation
\EqRef{s|ab|<s|a|s|b| n=k}
it follows that the inequation
\EqRef{s|ab|<s|a|s|b|}
is true for $n=k$.
\end{proof}

\begin{theorem}
\label{theorem: s|ab|<s|a|s|b| infty}
Let
\ShowEq{a,b>0 infty}
Then
\ShowEq{s|ab|<s|a|s|b| infty}
\end{theorem}
\begin{proof}
The theorem follows from the theorem
ref{theorem: s|ab|<s|a|s|b|}
when $n\rightarrow\infty$.
\end{proof}

\begin{theorem}
Let $A_i$, $i=1$, $2$, $3$,
be $D$\Hyph module with norm $\|x\|_i$
and normal Schauder basis $\Basis e_i$.
Let a map
\ShowEq{f+A1->A2}
Let a map
\ShowEq{g+A2->A3}
Then a map
\ShowEq{gf+A1->A3}
\end{theorem}
\begin{proof}
According to the statement
\RefItem{norm on D module 3, 1},
\ShowEq{gf+A1->A3 1}
From the theorem
\ref{theorem: s|ab|<s|a|s|b| infty}
and the inequation
\EqRef{gf+A1->A3 1},
it follows that
\ShowEq{gf+A1->A3 2}
From the inequation
\EqRef{gf+A1->A3 2},
it follows that
\ShowEq{gf+A1->A3 3}
The theorem follows from the inequation
\EqRef{gf+A1->A3 3}.
\end{proof}

\begin{theorem}
\label{theorem: polylinear map, Schauder basis}
Let
\ShowEq{Ai 1n}
be $D$\Hyph module with norm $\|x\|_i$
and normal Schauder basis $\Basis e_i$.
Let
$A$
be $D$\Hyph module with norm $\|x\|$
and normal Schauder basis $\Basis e$.
Let
\ShowEq{A1...An->A}
be polylinear map,
\ShowEq{|f|<infty}
Let
\ShowEq{ai in Ai+}
Then
\ShowEq{a=f(a1...an) A+}
\end{theorem}
\begin{proof}
We will prove the theorem
by induction on $n$.

For $n=1$, the statement of the theorem follows from the theorem
\ref{theorem: Schauder basis, image of map}.

Let the statement of the theorem is true for $n=k-1$.
Let
\ShowEq{f A1...Ak-1 ->A}
be polylinear map,
\ShowEq{|f|<infty}
We can represent the map $f$ as composition of maps
\DrawEq{f1k=ha}{}
According to the theorem
\ref{theorem: |f|=|h|},
\ShowEq{|h|<infty}
According to the induction assumption
\ShowEq{h() in L+}
According to the theorem
\ref{theorem: |fx|<|f||x|1n},
\ShowEq{|h()|<infty}
According to the theorem
\ref{theorem: Schauder basis, image of map}
\ShowEq{h()a in A+}
\end{proof}

According to the definition
\xRef{8433-5163}{definition: algebra over ring},
algebra is a module in which the product is defined
as bilinear map
\ShowEq{xy=Cxy}
We require
\ShowEq{|C|<infty}

\begin{convention}
Let $\Basis e$ be Schauder basis of free $D$\Hyph algebra $A$.
The product of basis vectors in $D$\Hyph algebra $A$
is defined according to rule
\DrawEq{product of basis vectors}{Schauder basis}
where
\ShowEq{structural constants of algebra}
are \AddIndex{structural constants}{structural constants}
of $D$\Hyph algebra $A$.
Since the product of vectors of the basis $\Basis e$ of
$D$\Hyph algebra $A$ is a vector of $D$\Hyph algebra $A$, then
we require that
the sequence
\ShowEq{structural constants of algebra in Schauder basis}
has limit for any $\gi i$, $\gi j$.
\qed
\end{convention}

\begin{theorem}
Let $\Basis e$ be Schauder basis of free $D$\Hyph algebra.
Then for any
\ShowEq{a b in basis of algebra plus}
product defined according to rule
\DrawEq{product ab in algebra}{Schauder basis}
is defined properly.
Under this condition
\ShowEq{ab in A+}
\end{theorem}
\begin{proof}
Since the product in the algebra is a bilinear map,
then we can write the product of $a$ and $b$ as
\DrawEq{product ab in algebra, 1}{Schauder basis}
From equations
\eqRef{product of basis vectors}{Schauder basis},
\eqRef{product ab in algebra, 1}{Schauder basis},
it follows that
\DrawEq{product ab in algebra, 2}{Schauder basis}
Since $\Basis e$ is a basis of the algebra $A$, then the equation
\eqRef{product ab in algebra}{Schauder basis}
follows from the equation
\eqRef{product ab in algebra, 2}{Schauder basis}.

From the theorem
\ref{theorem: polylinear map, Schauder basis},
it follows that
\ShowEq{ab in A+}
\end{proof}

\ifx\texFuture\Defined
\begin{theorem}
\label{theorem: standard component of tensor, algebra, Hamel basis}
Let $A_1$, ..., $A_n$ be
free algebras over commutative ring $D$.
Let $\Basis e_i$ be Hamel basis of $D$\Hyph algebra $A_i$.
Then the set of vectors
\ShowEq{tensor, algebra, Hamel basis}
is Hamel basis of tensor product
\ShowEq{tensor, algebra, Hamel basis, 1}
\end{theorem}
\begin{proof}
To prove the theorem, we need to consider the diagram
\xEqRef{8433-5163}{diagram top, algebra, tensor product}
which we used to prove the theorem
\xRef{8433-5163}{theorem: tensor product of algebras is module}.
\ShowEq{standard component of tensor, algebra, Hamel basis, diagram}
Let $M_1$ be module over ring $D$ generated
by product $\Times$
of $D$\Hyph algebras $A_1$, ..., $A_n$.
\begin{itemize}
\item
Let vector $b\in M_1$ have finite expansion
relative to the basis $\Times$
\ShowEq{b in M_1, Hamel basis}
where $I_1$ is finite set.
Let vector $c\in M_1$ have finite expansion
relative to the basis $\Times$
\ShowEq{c in M_1, Hamel basis}
where $I_2$ is finite set.
The set
\ShowEq{I=I1+I2}
is finite set.
Let
\ShowEq{b I2, c I1}
Then
\ShowEq{b+c in M_1, Hamel basis}
where $I$ is finite set.
Similarly, for $d\in D$
\ShowEq{db in M_1, Hamel basis}
where $I_1$ is finite set.
Therefore, we proved the following statement.\footnote{The set
$\Times$ cannot be Hamel basis because
this set is not countable.}
\begin{lemma}
\label{lemma: M is submodule}
The set $M$ of vectors of module $M_1$, which have finite expansion
relative to the basis $\Times$, is submodule of module $M_1$.
\end{lemma}
\end{itemize}
Injection
\ShowEq{map i, algebra, tensor product, 1}
is defined according to rule
\ShowEq{map i, algebra, Hamel basis, tensor product}
Let $N\subset M$ be submodule generated by elements of the following type
\ShowEq{kernel, algebra, Hamel basis, tensor product}
where $d_i\in A_i$, $c_i\in A_i$, $a\in D$.
Let
\[
j:M\rightarrow M/N
\]
be canonical map on factor module.
Since elements \EqRef{kernel 1, algebra, Hamel basis, tensor product}
and \EqRef{kernel 2, algebra, Hamel basis, tensor product}
belong to kernel of linear map $j$,
then, from equation
\EqRef{map i, algebra, Hamel basis, tensor product},
it follows
\ShowEq{f, algebra, Hamel basis, tensor product}
From equations \EqRef{f 1, algebra, Hamel basis, tensor product}
and \EqRef{f 2, algebra, Hamel basis, tensor product}
it follows that map $f$ is polylinear over ring $D$.

The module $M/N$ is tensor product $A_1\otimes...\otimes A_n$;
the map $j$ has form
\ShowEq{map j, Hamel basis}
and the set of
tensors like
\ShowEq{tensor, algebra, Hamel basis}
is countable basis of the module $M/N$.
According to the lemma
\ref{lemma: M is submodule},
arbitrary vector
\ShowEq{b in M, Hamel basis, 0}
has representation
\ShowEq{b in M, Hamel basis}
where $I$ is finite set.
According to the definition
\EqRef{map j, Hamel basis}
of the map $j$
\ShowEq{j(b), Hamel basis}
where $I$ is finite set.
Since $\Basis e_k$ is Hamel basis of $D$\Hyph algebra $A_k$,
then for any set of indexes $k\cdot i$, in equation
\ShowEq{aki, Hamel basis}
the set of values
\ShowEq{aki 1, Hamel basis}
which are different from $0$, is finite.
Therefore, the equation
\EqRef{j(b), Hamel basis}
has form
\ShowEq{j(b), 1, Hamel basis}
where the set of values
\ShowEq{j(b), 2, Hamel basis}
which are different from $0$, is finite.
\end{proof}

\begin{corollary}
\label{corollary: tensor product, Hamel basis}
Let $A_1$, ..., $A_n$ be
free algebras over commutative ring $D$.
Let $\Basis e_i$ be Hamel basis of $D$\Hyph algebra $A_i$.
Then any tensor $a\in A_1\otimes...\otimes A_n$
has finite set of standard components
different from $0$.
\qed
\end{corollary}

\begin{theorem}
Let $A_1$ be algebra
over the ring $D$.
Let $A_2$ be free associative algebra
over the ring $D$ with Hamel basis $\Basis e$.
The map
\ShowEq{standard representation of map A1 A2, 1, Hamel basis}
generated by the map $f\in\LDA$
through the tensor $a\in\ATwo$, has the standard representation
\ShowEq{standard representation of map A1 A2, 2, Hamel basis}
\end{theorem}
\begin{proof}
According to theorem
\ref{theorem: standard component of tensor, algebra, Hamel basis},
the standard representation of the tensor $a$ has form
\ShowEq{standard representation of map A1 A2, 3, Hamel basis}
The equation
\EqRef{standard representation of map A1 A2, 2, Hamel basis}
follows from equations
\EqRef{standard representation of map A1 A2, 1, Hamel basis},
\EqRef{standard representation of map A1 A2, 3, Hamel basis}.
\end{proof}
\fi

%% file: Schauder.Basis.Eq.tex

\def\LcDA{\mathcal{LC}(D;A_1;A_2)}

\DefEquation
{
n^{-1}d=dn^{-1}
}
{inverse integer in ring}

\DefEq
{
\[
f\in\mathcal{LC}(D;A_1,...,A_n;A)
\]
}
{polymap 5, f}

\DefEq
{
$a_1\in A_1$, ..., $a_{n-1}\in A_{n-1}$
}
{polymap 5, a}

\DefEq
{
\[
h\in\mathcal{LC}(D;A_1,...,A_{n-1};\mathcal{LC}(D;A_n;A))
\]
}
{polymap 5, h}

\DefEquation
{
\lim_{\gi m\rightarrow\infty}\lim_{\gi n\rightarrow\infty}
\sum_{\gi j=\gi 1}^{\gi m}
\sum_{\gi i=\gi 1}^{\gi n}
f^{\gi i}_{\gi j}(e^{\gi j}_1,e_{2\cdot\gi i})
}
{map in Schauder basis}

\DefEq
{
\[
\lim_{\gi j\rightarrow\infty}
\left\|\sum_{\gi i=\gi 1}^{\gi n}
f^{\gi i}_{\gi j}e_{2\cdot\gi i}\right\|_2=0
\]
}
{map in Schauder basis 1}

\DefEq
{
$(e^{\gi j}_1,e_{2\cdot\gi i})$.
}
{L(A1,A2) Schauder basis}

\DefEquation
{
n^{-1}d=n^{-1}dnn^{-1}=dn^{-1}
}
{inverse integer in ring, 2}

\DefEq
{
\[
|a_n-a|<\frac{\epsilon}{|d|}
\]
}
{an-a d}

\DefEq
{
\begin{align*}
|a_nd-ad|=|(a_n-a)d|&=|a_n-a||d|<\frac{\epsilon}{|d|}|d|=\epsilon
\\
|da_n-da|=|d(a_n-a)|&=|d||a_n-a|<|d|\frac{\epsilon}{|d|}=\epsilon
\end{align*}
}
{an-a d1}

\DefEq
{
\[
pd=\lim_{n\rightarrow\infty}(p_nd)=\lim_{n\rightarrow\infty}(dp_n)=dp
\]
}
{complete ring and real number}

\DefEq
{
\begin{align*}
(a,v)\in D\times A&\rightarrow av\in A
\\
(v,w)\in A\times A&\rightarrow vw\in A
\end{align*}
}
{topological D algebra}

\DefEq
{
\symb{\|a\|}0{norm on D module}
\[a\in A\rightarrow \ShowSymbol{norm on D module}\in R\]
}
{norm on D module}

\DefEq
{
\item $\|a\|\ge 0$
\label{item: norm on D module 1}
}
{norm on D module 1}

\DefEq
{
\item $\|a\|=0$
\label{item: norm on D module 2}
}
{norm on D module 2, 1}

\DefEq
{
$a=0$
}
{norm on D module 2, 2}

\DefEquation
{
\left\|\frac 1{\|x\|_1}x\right\|_1
=\frac 1{\|x\|_1}\|x\|_1=1
}
{x/|x|}

\DefEquation
{
\begin{matrix}
\displaystyle
\left\|\frac 1{\|x_1\|_1}x_1\right\|_1
=\frac 1{\|x_1\|_1}\|x_1\|_1=1
&...&
\displaystyle
\left\|\frac 1{\|x_n\|_n}x_n\right\|_n
=\frac 1{\|x_n\|_n}\|x_n\|_n=1
\end{matrix}
}
{x/|x|1n}

\DefEquation
{
\left\|\frac 1{\|x\|_1}f\circ x\right\|_2
=\left\|f\circ\left(\frac 1{\|x\|_1} x\right)\right\|_2
\le\|f\|
}
{f x/|x|}

\DefEquation
{
\left\|\frac 1{\|x_1\|_1...\|x_n\|_n}f\circ(x_1,...,x_n)\right\|
=\left\|f\circ\left(\frac 1{\|x_1\|_1} x_1,...,
\frac 1{\|x_n\|_n} x_n\right)\right\|
\le\|f\|
}
{f x/|x|1n}

\DefEq
{
$\|f\|_{2\cdot 3}$,
}
{|f|23}

\DefEquation
{
\frac 1{\|x\|_1}\|f\circ x\|_2
\le\|f\|
}
{f x/|x| 1}

\DefEquation
{
\frac 1{\|x_1\|_1...\|x_n\|_n}\|f\circ(x_1,...,x_n)\|
\le\|f\|
}
{f x/|x|1n 1}

\DefEq
{
\item $\|a+b\|\le \|a\|+\|b\|$
\label{item: norm on D module 3, 1}
\item $\|da\|=|d|\,\|a\|$, $d\in D$, $a\in A$
\label{item: norm on D module 3, 2}
}
{norm on D module 3}

\DefEquation
{
|(g\circ f)^{\gi i}_{\gi j}|
=
\left|\sum_{\gi k=\gi 1}^{\gi\infty}
g^{\gi i}_{\gi k}f^{\gi k}_{\gi j}\right|
\le
\sum_{\gi k=\gi 1}^{\gi\infty}
|g^{\gi i}_{\gi k}f^{\gi k}_{\gi j}|
}
{gf+A1->A3 1}

\DefEquation
{
|(g\circ f)^{\gi i}_{\gi j}|
\le
\sum_{\gi k=\gi 1}^{\gi\infty}
|g^{\gi i}_{\gi k}|
\sum_{\gi k=\gi 1}^{\gi\infty}
|f^{\gi k}_{\gi j}|
}
{gf+A1->A3 2}

\DefEquation
{
\begin{split}
\sum_{\gi i=\gi 1}^{\gi\infty}
\sum_{\gi j=\gi 1}^{\gi\infty}
|(g\circ f)^{\gi i}_{\gi j}|
&\le
\sum_{\gi i=\gi 1}^{\gi\infty}
\sum_{\gi j=\gi 1}^{\gi\infty}
\left(
\sum_{\gi k=\gi 1}^{\gi\infty}
|g^{\gi i}_{\gi k}|
\sum_{\gi k=\gi 1}^{\gi\infty}
|f^{\gi k}_{\gi j}|
\right)
\\
&=
\left(
\sum_{\gi i=\gi 1}^{\gi\infty}
\sum_{\gi k=\gi 1}^{\gi\infty}
|g^{\gi i}_{\gi k}|
\right)
\left(
\sum_{\gi j=\gi 1}^{\gi\infty}
\sum_{\gi k=\gi 1}^{\gi\infty}
|f^{\gi k}_{\gi j}|
\right)
<\infty
\end{split}
}
{gf+A1->A3 3}

\DefEq
{$a_i\ge 0$, $b_i\ge 0$, $i=1$, ...,$n$.
}
{a,b>0}

\DefEq
{$a_i\ge 0$, $b_i\ge 0$, $i=1$, ...,$\infty$.
}
{a,b>0 infty}

\DefEquation
{
\sum_{i=1}^na_ib_i<\sum_{i=1}^na_i\sum_{i=1}^nb_i
}
{s|ab|<s|a|s|b|}

\DefEquation
{
\sum_{i=1}^{\infty}a_ib_i<\sum_{i=1}^{\infty}a_i\sum_{i=1}^{\infty}b_i
}
{s|ab|<s|a|s|b| infty}

\DefEquation
{
a_1b_1+a_2b_2\le(a_1+a_2)(b_1+b_2)
=a_1b_1+\underline{a_1b_2+a_2b_1}+a_2b_2
}
{s|ab|<s|a|s|b| n=2}

\DefEquation
{
\sum_{i=1}^{k-1}a_ib_i+a_kb_k
\le\left(\sum_{i=1}^{k-1}a_i+a_k\right)\left(\sum_{i=1}^{k-1}b_i+b_k\right)
}
{s|ab|<s|a|s|b| n=k}

\DefEquation
{
\|a\|\le\|a-b\|+\|b\|
}
{normed module, norm of difference, 2}

\DefEq
{
\[
a=(a-b)+b
\]
}
{normed module, norm of difference, 1}

\DefEquation
{
\|a-b\|\ge \|a\|-\|b\|
}
{normed module, norm of difference}

\DefEq
{
$\{a_n\}$
\symb{\lim_{n\rightarrow\infty}a_n}0{limit of sequence}
\[
a=\ShowSymbol{limit of sequence}
\]
}
{limit of sequence, normed module}

\DefEq
{
$\|a_n-a\|<\epsilon$
}
{an-a}

\DefEq
{
$\|a_p-a_q\|<\epsilon$
}
{ap-aq}

\DefEq
{
\symb{\|a\|}1{norm on d algebra}
}
{norm on d algebra}

\DefEquation
{
\|ab\|\le\|a\|\,\|b\|
}
{norm on d algebra 3}

\DefEq
{
$a\in A$, $\|a\|=1$,
}
{unit sphere in algebra}

\DefEq
{
\[f:A_1\rightarrow A_2\]
}
{f A1 A2}

\DefEq
{
\[g:A_1\rightarrow A_2\]
}
{g A1 A2}

\DefEq
{
\[f:A_2\rightarrow A_3\]
}
{f A2 A3}

\DefEquation
{
\|(x'+b)-(x+b)\|_1=\|x'-x\|_1<\delta
}
{continuous map, x 1}

\DefEq
{
\|x'-x\|_1<\delta
}
{continuous map, x}

\DefEq
{
\begin{matrix}
\|x'_1-x_1\|_1<\delta
&...&
\|x'_n-x_n\|_n<\delta
\end{matrix}
}
{|xi|<e, 1-n}

\DefEq
{
\[
\begin{matrix}
\|x'_1-x_1\|_1<\delta
&...&
\|x'_{n-1}-x_{n-1}\|_{n-1}<\delta
\end{matrix}
\]
}
{polymap 4, x}

\DefEq
{
$\|x\|_1<\delta_1$
}
{continuous map, x, 1}

\DefEq
{
$\|x'_1-x_1\|_1<\delta_1$, ..., $\|x'_n-x_n\|_n<\delta_1$
}
{continuous multivariable map, x, 1}

\DefEq
{
$\|x'_1-x_1\|_1<\delta_2$, ..., $\|x'_n-x_n\|_n<\delta_2$
}
{continuous multivariable map, x, 2}

\DefEq
{
$\|x\|_1<\delta_2$
}
{continuous map, x, 2}

\DefEq
{
\|f(x')-f(x)\|_2<\epsilon
}
{continuous map, y}

\DefEquation
{
\|f\circ x'-f\circ x\|_2<\epsilon
}
{continuous map, y at x}

\DefEq
{
\|f(x'_1,...,x'_n)-f(x_1,...,x_n)\|<\epsilon
}
{continuous multivariable map, y}

\DefEq
{
\[\|f\circ(x'_1,...,x'_n)-f\circ(x_1,...,x_n)\|<\epsilon\]
}
{polymap 30, y}

\DefEquation
{
\|f\circ(x'_1,...,x'_n)-f\circ(x_1,...,x_n)\|
}
{polymap 5, y}

\DefEquation
{
\|f\circ(x_1,...,x_{n-1},x'_n)-f\circ(x_1,...,x_{n-1},x_n)\|<\epsilon
}
{polymap 3, y}

\DefEquation
{
\|f\circ(x'_1,...,x'_{n-1},x_n)-f\circ(x_1,...,x_{n-1},x_n)\|<\epsilon\|x_n\|_n
}
{polymap 4, y}

\DefEquation
{
\begin{split}
&f\circ(x'_1,...,x'_n)-f\circ(x_1,...,x_n)
\\=
&(h\circ(x'_1,...,x'_{n-1}))\circ x'_n-(h\circ(x_1,...,x_{n-1}))\circ x_n
\\=
&(h\circ(x'_1,...,x'_{n-1}))\circ x'_n-(h\circ(x'_1,...,x'_{n-1}))\circ x_n
\\+
&(h\circ(x'_1,...,x'_{n-1}))\circ x_n-(h\circ(x_1,...,x_{n-1}))\circ x_n
\end{split}
}
{polymap 51, y}

\DefEquation
{
\begin{split}
&\|f\circ(x'_1,...,x'_n)-f\circ(x_1,...,x_n)\|
\\ \le
&\|(h\circ(x'_1,...,x'_{n-1}))\circ x'_n-(h\circ(x'_1,...,x'_{n-1}))\circ x_n\|
\\+
&\|(h\circ(x'_1,...,x'_{n-1}))\circ x_n-(h\circ(x_1,...,x_{n-1}))\circ x_n\|
\end{split}
}
{polymap 52, y}

\DefEq
{
\[\|(h\circ(x_1,...,x_{n-1}))\circ x'_n-(h\circ(x_1,...,x_{n-1}))\circ x_n\|
<\epsilon\]
}
{polymap 31, y}

\DefEquation
{
\|(h\circ(x'_1,...,x'_{n-1}))\circ x'_n-(h\circ(x'_1,...,x'_{n-1}))\circ x_n\|
<\frac{\epsilon}2
}
{h(x)dx<e/2}

\DefEq
{
\[
\|(h\circ(x'_1,...,x'_{n-1}))\circ x_n-(h\circ(x_1,...,x_{n-1}))\circ x_n\|
=0<\frac{\epsilon}2
\]
}
{dh(x)x<e/2 1}

\DefEquation
{
\begin{split}
&\|(h\circ(x'_1,...,x'_{n-1}))\circ x_n-(h\circ(x_1,...,x_{n-1}))\circ x_n\|
\\
=&\|(h\circ(x'_1,...,x'_{n-1})-h\circ(x_1,...,x_{n-1}))\circ x_n\|
\\
<&\frac{\epsilon}{2\|x_n\|_n}\|x_n\|_n=\frac{\epsilon}2
\end{split}
}
{dh(x)x<e/2 2}

\DefEquation
{
\|(h\circ(x'_1,...,x'_{n-1}))\circ x_n-(h\circ(x_1,...,x_{n-1}))\circ x_n\|
<\frac{\epsilon}2
}
{dh(x)x<e/2 3}

\DefEquation
{
\|h\circ(x'_1,...,x'_{n-1})-h\circ(x_1,...,x_{n-1})\|
<\frac{\epsilon}{2\|x_n\|_n}
}
{dh(x)<e/2}

\DefEquation
{
\begin{split}
&\|(h\circ(x'_1,...,x'_{n-1}))\circ x_n-(h\circ(x_1,...,x_{n-1}))\circ x_n\|
\\=&\|(h\circ(x'_1,...,x'_{n-1})-h\circ(x_1,...,x_{n-1}))\circ x_n\|
\\<&\epsilon\|x_n\|_n
\end{split}
}
{polymap 41, y}

\DefEquation
{
\frac
{\|(h\circ(x'_1,...,x'_{n-1})-h\circ(x_1,...,x_{n-1}))\circ x_n\|}{\|x_n\|_n}
<\epsilon
}
{polymap 42, y}

\DefEquation
{
\begin{split}
&\|h\circ(x'_1,...,x'_{n-1})-h\circ(x_1,...,x_{n-1})\|
\\ \le
&\frac
{\|(h\circ(x'_1,...,x'_{n-1})-h\circ(x_1,...,x_{n-1}))\circ x_n\|}{\|x_n\|_n}
\end{split}
}
{polymap 43, y}

\DefEquation
{
\|h\circ(x'_1,...,x'_{n-1})-h\circ(x_1,...,x_{n-1})\|
<\epsilon
}
{polymap 44, y}

\DefEq
{
$h\circ(x_1,...,x_{n-1})$
}
{polymap 32, y}

\DefEquation
{
\|f\circ x\|_2<\frac{\epsilon}d
}
{continuous map, y, d}

\DefEquation
{
\|f(x'_1,...,x'_n)-f(x_1,...,x_n)\|<\frac{\epsilon}d
}
{continuous multivariable map, y, d}

\DefEquation
{
\|f\circ x\|_2<\frac{\epsilon}2
}
{continuous map, y, 1}

\DefEquation
{
\|f(x'_1,...,x'_n)-f(x_1,...,x_n)\|<\frac{\epsilon}2
}
{continuous multivariable map, y, 1}

\DefEquation
{
\|g(x'_1,...,x'_n)-g(x_1,...,x_n)\|<\frac{\epsilon}2
}
{continuous multivariable map, y, 2}

\DefEq
{
\[f\in\LDA\setminus\LcDA\]
}
{LcDA-LDA}

\DefEquation
{
\|g\circ x\|_2<\frac{\epsilon}2
}
{continuous map, y, 2}

\DefEquation
{
\begin{split}
\|f\circ(x'+b)-f\circ(x+b)\|_2
&=\|(f\circ x'+f\circ b)-(f\circ x+f\circ b)\|_2
\\&=\|f\circ x'-f\circ x\|_2<\epsilon
\end{split}
}
{continuous map, y 1}

\DefEq
{
\[\delta= \mathrm{min}(\delta_1, \delta_2)\]
}
{delta 1 2}

\DefEq
{
\[
\|(f+g)\circ x\|_2=\|f\circ x+g\circ x\|_2
\le\|f\circ x\|_2+\|g\circ x\|_2\le\epsilon
\]
}
{continuous map, y, 3}

\DefEq
{
\begin{align*}
&\|(f+g)(x'_1,...,x'_n)-(f+g)(x_1,...,x_n)\|
\\
=&\|f(x'_1,...,x'_n)+g(x'_1,...,x'_n)-f(x_1,...,x_n)-g(x_1,...,x_n)\|
\\
\le&\|f(x'_1,...,x'_n)-f(x_1,...,x_n)\|
+\|g(x'_1,...,x'_n)-g(x_1,...,x_n)\|\le\epsilon
\end{align*}
}
{continuous multivariable map, y, 3}

\DefEq
{
\begin{align*}
\|(d\,f)(x'_1,...,x'_n)-(d\,f)(x_1,...,x_n)\|
=&\|df(x'_1,...,x'_n)-df(x_1,...,x_n)\|
\\
=&|d|\,\|f(x'_1,...,x'_n)-f(x_1,...,x_n)\|\le\epsilon
\end{align*}
}
{continuous multivariable map, y, 4}

\DefEq
{
\[
\|(d\,f)\circ x\|_2=\|d(f\circ x)\|_2
=|d|\,\|f\circ x\|_2\le\epsilon
\]
}
{continuous map, y, 4}

\DefEq
{
$\|x\|_1<\delta$
}
{continuous map, x, 3}

\DefEq
{
$\|x'_1-x_1\|_1<\delta$, ..., $\|x'_n-x_n\|_n<\delta$
}
{continuous multivariable map, x, 3}

\DefEq
{
$\|x'_n-x_n\|_n<\delta$
}
{polymap 3, x}

\DefEq
{
$\|x'_n-x_n\|_n<\delta_1$
}
{|xn|<d1}

\DefEq
{
$\|x'_i-x_i\|_i<\delta_2$, $1\le i<n$
}
{|xi|<d2}

\DefEq
{
$\gi i=\gi j$,
$\|e_{\gi j}\|=\infty$.
}
{Schauder basis, norm is finite, i=j}

\DefEq
{
$\{a^{\gi i}\}_{\gi i=\gi 1}^{\gi\infty}$,
$a^{\gi i}\in D$,
}
{ai in D}

\DefEq
{
\[\sum_{\gi i=\gi 1}^{\gi\infty}|a^{\gi i}|<\infty\]
}
{ai in D 1}

\DefEq
{
$\epsilon\in R$, $\epsilon>0$,
}
{epsilon in R}

\DefEq
{
\begin{matrix}
|f^{\gi i}_{\gi j}|<\epsilon
&\gi i>\gi N&\gi j>\gi M
\end{matrix}
}
{|f|<epsilon}

\DefEq
{
\[
F_1=\max\{|f^{\gi i}_{\gi j}|,\gi 1\le\gi i\le\gi N,
\gi 1\le\gi j\le\gi M\}
\]
}
{F1=max f}

\DefEquation
{
\begin{matrix}
\|f^{\gi i}_{\gi j}(e^{\gi i}_1,e_{2\cdot\gi j})\|<\epsilon
&\gi i>\gi N&\gi j>\gi M
\end{matrix}
}
{|f|<epsilon 1}

\DefEquation
{
y_n=\frac 1{n\,\|x_n\|_1}x_n
}
{xn=}

\DefEq
{
\[f:A_1\times...\times A_k\rightarrow A\]
}
{f A1...Ak->A}

\DefEq
{
\[f:A_1\times...\times A_n\rightarrow A\]
}
{f A1...An->A}

\DefEq
{
\[h:A_1\times...\times A_{n-1}\rightarrow \mathcal L(D;A_n;A)\]
}
{h 1n-1->A}

\DefEq
{
\[h:A_1\times...\times A_{k-1}\rightarrow \mathcal L(D;A_k;A)\]
}
{h 1k-1->A}

\DefEq
{
$\|h\|<\infty$.
}
{|h|<infty}

\DefEq
{
$\|f\|=\|h\|<\infty$.
}
{|f|=|h|<infty}

\DefEq
{
\[h\circ(x_1,...,x_{k-1})\in\mathcal L^+(D;A_k;A)\]
}
{h() in L+}

\DefEq
{
\[\|h\circ(x_1,...,x_{k-1})\|<\infty\]
}
{|h()|<infty}

\DefEq
{
\[
a,b\in A^+(\Basis e)=>ab\in A^+(\Basis e)
\]
}
{ab in A+}

\DefEq
{
\[(h\circ(x_1,...,x_{k-1}))\circ x_k\in A^+(e)\]
}
{h()a in A+}

\DefEq
{
$a_i\in A$, $i\in I$,
}
{ai in I}

\DefEquation
{
\text{span}(a_i,i\in I)\subset A
}
{span in A}

\DefEq
{
$a_1$, $a_2\in A$,
}
{a1 a2 A}

\DefEq
{
$a\in A$, $d\in D$,
}
{aA dD}

\DefEq
{
\symb{A^+(\Basis e)}0{A plus Schauder}
\[
\ShowSymbol{A plus Schauder}=
\{a\in A:
a=a^{\gi i}\,e_{\gi i},
\ \sum_{\gi i=\gi 1}^{\gi\infty}|a^{\gi i}|<\infty\}
\]
}
{A plus Schauder}

\DefEquation
{
\|f\circ e_{1\cdot\gi i}\|_2\le\|f\|
}
{|fe|<|f|}

\DefEquation
{
\begin{matrix}
a_1\in A_1^+(\Basis e_1)&
a_1=a_1^{\gi i}\,e_{1\cdot\gi i}
\end{matrix}
}
{a1 in A1 plus}

\DefEq
{
\[
\|C\|<\infty
\]
}
{|C|<infty}

\DefEq
{
$a_2\in A_2^+(\Basis e_2)$.
}
{a2 in A2 plus}

\DefEquation
{
\sum_{\gi i=\gi 1}^{\gi\infty}|a_2^{\gi i}|
=
\sum_{\gi i=\gi 1}^{\gi\infty}|a_1^{\gi i}|\,\|f\circ e_{1\cdot\gi i}\|_2
<
\|f\|\sum_{\gi i=\gi 1}^{\gi\infty}|a_1^{\gi i}|<\infty
}
{|a2|<infty}

\DefEquation
{
\sum_{\gi i=\gi 1}^{\gi\infty}|a_2^{\gi i}|
=
\sum_{\gi i=\gi 1}^{\gi\infty}|a_1^{\gi i}
f_{\gi i}^{\gi j}|\,\|e_{1\cdot\gi i}\|_2
=
\sum_{\gi i=\gi 1}^{\gi\infty}|a_1^{\gi i}
f_{\gi i}^{\gi j}|
<\infty
}
{|a2|<infty 1}

\DefEq
{
\[
\|f\circ e_{1\cdot\gi i}\|_2\le\|f\|\,\|e_{1\cdot\gi i}\|_1
\]
}
{|fe|<|f|, 1}

\DefEq
{
\[\|e_{\gi i}\|<\infty\]
}
{Schauder basis, norm is finite}

\DefEquation
{
\frac{\|f(x)\|_2}{\|x\|_1}=\left\|f\left(\frac x{\|x\|_1}\right)\right\|_2
}
{norm of linear map, module, 1}

\DefEquation
{
\frac{\|f(x_1,...,x_n)\|}{\|x\|_1...\|x\|_n}
=\left\|f\left(\frac {x_1}{\|x_1\|_1},...,\frac {x_n}{\|x_n\|_n}\right)\right\|
}
{norm of polylinear map, module, 1}

\DefEq
{
\[
\|f\circ x-f\circ y\|_2=\|f\circ (x-y)\|_2\le \|f\|\ \|x-y\|_1
\]
}
{|f|<infty 1}

\DefEq
{
\[\|x-y\|_1<\delta\]
}
{|f|<infty 2}

\DefEq
{
\[\|f\circ x-f\circ y\|_2\le \|f\|\ \delta=\epsilon\]
}
{|f|<infty 3}

\DefEq
{
\[\LDA\rightarrow R\]
}
{LDA->R}

\DefEq
{
\[A'\rightarrow R\]
}
{A'->R}

\DefEq
{
\symb{\|f\|}0{norm of functional}
\[
\ShowSymbol{norm of functional}=
\text{sup}\frac{|f\circ x|}{\|x\|_0}
\]
}
{norm of functional}

\DefEq
{
\[a=(a-a^{\gi i}e_{\gi i})+a^{\gi i}e_{\gi i}\]
}
{a=a-+ai}

\DefEq
{
\[a=(a-a^{\gi i}e_{1\cdot\gi i})+a^{\gi i}e_{1\cdot\gi i}\]
}
{a=a-+ai L(A1,A2)}

\DefEq
{
\[\|a\|_1\le\|a-a^{\gi i}e_{\gi i}\|_1+\|a^{\gi i}e_{\gi i}\|_1\]
}
{a=a-+ai 1}

\DefEq
{
\[\|a\|_1\le\|a-a^{\gi i}e_{1\cdot\gi i}\|_1+\|a^{\gi i}e_{1\cdot\gi i}\|_1\]
}
{a=a-+ai 1 L(A1,A2)}

\DefEquation
{
\|e^{\gi i}\|_2=\text{sup}\frac{|e^{\gi i}\circ a|}{\|a\|_1}
=\frac{|a^{\gi i}|}{|a^{\gi i}|\,\|e_{\gi i}\|_1}
}
{a=a-+ai 4}

\DefEquation
{
\|(e^{\gi i}_1,e_{2\cdot\gi j})\|
=\text{sup}\frac{\|(e^{\gi i}_1,e_{2\cdot\gi j})\circ a\|_2}{\|a\|_1}
=\frac{|a^{\gi i}|\,\|e_{2\cdot\gi j}\|_2}{|a^{\gi i}|\,\|e_{1\cdot\gi i}\|_1}
}
{a=a-+ai 7}

\DefEq
{
$(e^{\gi i}_1,e_{2\cdot\gi j})$
}
{basis of L(A1,A2)}

\DefEq
{
\[\|(e^{\gi i}_1,e_{2\cdot\gi j})\|=1\]
}
{normal basis L(A1,A2)}

\DefEq
{
$a=a^{\gi i}e_{\gi i}$,
}
{a=a-+ai 2}

\DefEq
{
$a=a^{\gi i}e_{1\cdot\gi i}$,
}
{a=a-+ai 2 L(A1,A2)}

\DefEq
{
\[\|a\|_1=\|a^{\gi i}e_{\gi i}\|_1=|a^{\gi i}|\,\|e_{\gi i}\|_1\]
}
{a=a-+ai 3}

\DefEquation
{
\|a\|_1=\|a^{\gi i}e_{1\cdot\gi i}\|_1=|a^{\gi i}|\,\|e_{1\cdot\gi i}\|_1
}
{a=a-+ai 3 L(A1,A2)}

\DefEquation
{
\|(e^{\gi i}_1,e_{2\cdot\gi j})\|
=\frac {\|e_{2\cdot\gi j}\|_2}{\|e_{1\cdot\gi i}\|_1}
}
{norm of basis L(A1,A2)}

\DefEquation
{
\|e^{\gi i}\|_2=\frac 1{\|e_{\gi i}\|_1}
}
{norm of conjugate basis}

\DefEq
{
\[(e^{\gi i}_1,e_{2\cdot\gi j})\circ a=a^{\gi i}e_{2\cdot\gi j}\]
}
{a=a-+ai 5}

\DefEquation
{
\|(e^{\gi i}_1,e_{2\cdot\gi j})\circ a\|_2=
\|a^{\gi i}e_{2\cdot\gi j}\|_2=
|a^{\gi i}|\,\|e_{2\cdot\gi j}\|_2
}
{a=a-+ai 6}

\DefEq
{
\symb{\|f\|}0{norm of map}
}
{norm of map}

\DefEquation
{
\ShowSymbol{norm of map}=
\text{sup}\frac{\|f\circ x\|_2}{\|x\|_1}
}
{norm of map, module}

\DefEquation
{
\|a_3\|_3\le\|h\circ a_1\|_{2\cdot 3}\|a_2\|_2
}
{|a3|=|ha1||a2|}

\DefEquation
{
\|f\circ(x_1,...,x_n)\|
=\|(h\circ(x_1,...,x_{n-1}))\circ x_n\|
\le\|h\circ (x_1,..,x_{n-1})\|\,\|x_n\|_n
}
{|x|=|hx1n-1||xn|}

\DefEquation
{
\|f\|=\text{sup}\{\|f\circ x\|_2:\|x\|_1=1\}
}
{norm of linear map, module}

\DefEquation
{
\|f\|=\text{sup}\{\|f\circ(x_1,...,x_n)\|:\|x_i\|_i=1,1\le i\le n\}
}
{norm of polylinear map, module}

\DefEq
{
$A_i$, $i=1$, ..., $k$,
}
{A1k}

\DefEquation
{
\|h\circ a_1\|_{2\cdot 3}\le\|h\|\,\|a_1\|_1
}
{|ha1|}

\DefEquation
{
\|h\circ(x_1,...,x_{n-1})\|\le\|h\|\,\|x_1\|_1...\|x_{n-1}\|_{n-1}
}
{|hx1n-1|}

\DefEquation
{
\|f\circ(a_1,...,a_k)\|\le\|f\|\,\|a_1\|...\|a_k\|
}
{|ha1k-1|}

\DefEquation
{
\|f\circ x\|_2\le\|f\|\,\|x\|_1
}
{|fx|<|f||x|}

\DefEq
{
\|f\circ(x_1,...,x_n)\|\le\|f\|\,\|x_1\|_1...\|x_n\|_n
}
{|fx|<|f||x|1n}

\DefEq
{
$\mathcal L^+(D;A_1(\Basis e_1);A_2(\Basis e_2))$
}
{L plus A1 A2}

\DefEq
{
\[f\in\mathcal L^+(D;A_1(\Basis e_1);A_2(\Basis e_2))\]
}
{f+A1->A2}

\DefEq
{
\[g\in\mathcal L^+(D;A_2(\Basis e_2);A_3(\Basis e_3))\]
}
{g+A2->A3}

\DefEq
{
\[g\circ f\in\mathcal L^+(D;A_1(\Basis e_1);A_3(\Basis e_3))\]
}
{gf+A1->A3}

\DefEq
{
$F<\infty$
}
{F<infty}

\DefEquation
{
|f^{\gi i}_{\gi j}|\le F
}
{|f|<<F}

\DefEq
{
$(\Basis e_1,\Basis e^2)$
}
{Basis e1 e2}

\DefEquation
{
\frac{\|(h\circ(x_1,...,x_{n-1}))\circ x_n)\|}{\|x_n\|_n}
\le\|f\|\,\|x_1\|_1...\|x_{n-1}\|_{n-1}
}
{|hx|/|xn|}

\DefEquation
{
\|h\circ(x_1,...,x_{n-1})\|
\le\|f\|\,\|x_1\|_1...\|x_{n-1}\|_{n-1}
}
{|hx|<|f||x|}

\DefEq
{
\[\|f\circ x\|_2=0\]
}
{|f|=0 1}

\DefEq
{
$f\circ x=0$
}
{|f|=0 2}

\DefEquation
{
\begin{array}{r@{\,}l}
\displaystyle\text{sup}\frac{\|(f_1+f_2)\circ x\|_2}{\|x\|_1}
&\displaystyle=\text{sup}\frac{\|f_1\circ x+f_2\circ x\|_2}{\|x\|_1}
\\[6pt]
\ &\displaystyle\le\text{sup}\frac{\|f_1\circ x\|_2+\|f_2\circ x\|_2}{\|x\|_1}
\\&\displaystyle\le\text{sup}\frac{\|f_1\circ x\|_2}{\|x\|_1}
+\text{sup}\frac{\|f_2\circ x\|_2}{\|x\|_1}
\end{array}
}
{|f1+f2| 1}

\DefEquation
{
\text{sup}\frac{\|(d\,f)\circ x\|_2}{\|x\|_1}
=\text{sup}\frac{\|d(f\circ x)\|_2}{\|x\|_1}
\le\text{sup}\frac{|d|\,\|f\circ x\|_2}{\|x\|_1}
=|d|\,\text{sup}\frac{\|f\circ x\|_2}{\|x\|_1}
}
{|df| 1}

\DefEquation
{
\begin{array}{r@{\,}l}
\displaystyle\text{sup}\frac{\|(f\circ g)\circ x\|_3}{\|x\|_1}
&\displaystyle=\text{sup}\frac{\|f\circ(g\circ x)\|_3}{\|x\|_1}
=\text{sup}\left(\frac{\|f\circ(g\circ x)\|_3}{\|g\circ x\|_2}
\frac{\|g\circ x\|_2}{\|x\|_1}\right)
\\
&\displaystyle\le\text{sup}\frac{\|f\circ(g\circ x)\|_3}{\|g\circ x\|_2}
\,\text{sup}\frac{\|g\circ x\|_2}{\|x\|_1}
\end{array}
}
{|fg| 1}

\DefEq
{
\[h\circ a_1:A_2\rightarrow A_3\]
}
{h a_1}

\DefEquation
{
\|f\|=\text{sup}\frac{\|f\circ(a_1,a_2)\|_3}{\|a_1\|_1\,\|a_2\|_2}
}
{|f| a1a2}

\DefEq
{
\[f\circ(a_1,a_2)=(h\circ a_1)\circ a_2\]
}
{f a1,a2}

\DefEq
{
f\circ(x_1,...,x_n)=(h\circ(x_1,...,x_{n-1}))\circ x_n
}
{f1n=hx}

\DefEq
{
f\circ(a_1,...,a_k)=(h\circ(a_1,...,a_{k-1}))\circ a_k
}
{f1k=ha}

\DefEq
{
$h\circ(a_1,...,a_{n-1})$
}
{polymap 3, algebra}

\DefEquation
{
h:A_1\rightarrow\mathcal L(D;A_2;A_3)
}
{h A1 A2->A3}

\DefEquation
{
f:A_1\times A_2\rightarrow A_3
}
{f A1,A2->A3}

\DefEquation
{
\|a_3\|_3\le\|h\|\,\|a_1\|_1\,\|a_2\|_2
}
{a3 h a1 a2}

\DefEquation
{
\|f\circ(x_1,...,x_n)\|
\le\|h\|\,\|x_1\|_1\,...\,\|x_n\|_n
}
{x h x1 xn}

\DefEquation
{
\|f\|=\|h\|
}
{|f|=|h|}

\DefEquation
{
\|f\|\le\|h\|
}
{|f|<|h|}

\DefEquation
{
\|h\|\le\|f\|
}
{|h|<|f|}

\DefEq
{
$\mathcal L(D;A_1;\mathcal L(D;A_2;A_3))$
}
{L(A1,L(A2,A3))}

\DefEq
{
\[(f,g)\rightarrow f\circ g\]
}
{f g ->}

\DefEq
{
$g\circ A_1\ne A_2$,
}
{g A1}

\DefEquation
{
\text{sup}\frac{\|(f\circ g)\circ x\|_3}{\|x\|_1}
\le\text{sup}\frac{\|f\circ y\|_3}{\|y\|_2}
\,\text{sup}\frac{\|g\circ x\|_2}{\|x\|_1}
}
{|fg| 2}

\DefEquation
{
\text{sup}\frac{\|f\circ(g\circ x)\|_3}{\|g\circ x\|_2}
\le\text{sup}\frac{\|f\circ y\|_3}{\|y\|_2}
}
{|fg| 0}

\DefEq
{
\[\|f_1+f_2\|\le\|f_1\|+\|f_2\|\]
}
{|f1+f2| 2}

\DefEq
{
\[\|d\,f\|=|d|\,\|f\|\]
}
{|df| 2}

\DefEquation
{
\begin{matrix}
f(rx)=rf(x)
&
r\in R
\end{matrix}
}
{linear map and real number}

\DefEquation
{
\begin{matrix}
f(r_1x_1,...,r_nx_n)=r_1...r_nf(x_1,...,x_n)
&
r_1,...,r_n\in R
\end{matrix}
}
{polylinear map and real number}

\DefEq
{
\[
\frac{\|f(rx)\|_2}{\|rx\|_1}
=\frac{|r|\ \|f(x)\|_2}{|r|\ \|x\|_1}
=\frac{\|f(x)\|_2}{\|x\|_1}
\]
}
{linear map and real number, 1}

\DefEq
{
\[
\frac{\|f(r_1x_1,...,r_nx_n)\|}{\|r_1x_1\|_1...\|r_nx_n\|_n}
=\frac{|r_1|...|r_n|\ \|f(x_1,...,x_n)\|}
{|r_1|\,\|x\|_1...|r_n|\,\|x_n\|_n}
=\frac{\|f(x_1,...,x_n)\|}{\|x_1\|_1...\|x_n\|_n}
\]
}
{polylinear map and real number, 1}

\DefEq
{
$f\circ g$
}
{fg}

\DefEquation
{
\|f\circ g\|\le\|f\|\,\|g\|
}
{|fg|}

\DefEq
{
\[\text{sup}(f_1(x)f_2(x))\le\text{sup}f_1(x)\text{sup}f_2(x)\]
}
{fi A1 R 1}

\DefEq
{
\[
\begin{matrix}
f_1:A_1\rightarrow R&f_2:A_1\rightarrow R
\end{matrix}
\]
}
{fi A1 R}

\DefEq
{
$\|f\|=\infty$.
}
{|f|=infty}

\DefEq
{
\[
\|f\|<\sum_{\gi i=\gi 1}^{\gi\infty}
\sum_{\gi j=\gi 1}^{\gi\infty}|f^{\gi i}_{\gi j}|
\]
}
{|f|<ij|f|}

\DefEq
{
$\|f\|<\infty$.
}
{|f|<infty}

\DefEquation
{
\|a\|<\sum_{\gi i=\gi 1}^{\gi\infty}|a^{\gi i}|
}
{|a|<|ai|}

\DefEquation
{
\|a\|<\sum_{\gi i=\gi 1}^{\gi\infty}|a^{\gi i}|\,\|e_{\gi i}\|
}
{|a|<|ai| 1}

\DefEquation
{
\|f\circ x_n\|_2>n\,\|x_n\|_1
}
{|f(x)|>n}

\DefEq
{
\[
\|f\circ y_n\|_2=
\left\|f\circ\left(\frac 1{n\,|x_n|_1}x_n\right)\right\|_2
=\frac 1{n\,\|x_n\|_1}\|f\circ x_n\|_2>1
\]
}
{|f(y)|=1}

\DefEq
{
$A_i$, $i=1$, ..., $n$,
}
{Ai 1n}

\DefEq
{
\[
f:A_1\times...\times A_n\rightarrow A
\]
}
{A1...An->A}

\DefEq
{
\[
g:A_1\times...\times A_n\rightarrow A
\]
}
{g A1...An->A}

\DefEq
{
\[
f:A_1\times...\times A_{k-1}\rightarrow A
\]
}
{f A1...Ak-1 ->A}

\DefEq
{
$a_i\in A_i^+(\Basis e_i)$.
}
{ai in Ai+}

\DefEq
{
\[
\begin{matrix}
a=f\circ(a_1,...,a_n)&a\in A^+(\Basis e)
\end{matrix}
\]
}
{a=f(a1...an) A+}

\DefEq
{
\[
xy=C(x,y)
\]
}
{xy=Cxy}

\DefEq
{
\[
\sum_{\gi k=\gi 1}^{\gi\infty}C^{\gi k}_{\gi{ij}}\,e_{\gi k}
\]
}
{structural constants of algebra in Schauder basis}

\DefEq
{
\[
\begin{matrix}
a=a^{\gi i}\,e_{\gi i}&b=b^{\gi i}\,e_{\gi i}
&a, b\in A^+(\Basis e)
\end{matrix}
\]
}
{a b in basis of algebra plus}

\DefEq
{
\[
ab\in A^+(\Basis e)
\]
}
{ab in A+}

\DefEq
{
\symb{\lim_{n\rightarrow\infty}a_n}0{limit of sequence}
\[
a=\ShowSymbol{limit of sequence}
\]
}
{limit of sequence, normed ring}

\DefEq
{
$\gi n>\gi j$
\[
\left\|
\sum_{i=1}^na^{\gi i}\,e_{\gi i}
\right\|
<
\sum_{i=1}^n\|a^{\gi i}\,e_{\gi i}\|
=
\sum_{i=1}^n|a^{\gi i}|\,\|e_{\gi i}\|
=\infty
\]
}
{Schauder basis, norm is finite, n>j}

\DefEq
{
\symb{\{e_{\gi i}\}_{\gi i=\gi 1}^{\gi\infty}}0{Schauder basis}
$\Basis e=\ShowSymbol{Schauder basis}$
}
{Schauder basis}

\DefEq
{
\symb{\{a^{\gi i}\}_{\gi i=\gi 1}^{\gi\infty}}1{Schauder basis, coordinates},
$a^{\gi i}\in D$,
}
{Schauder basis, coordinates}

\DefEq
{
$a^{\gi j}\ne 0$.
}
{Schauder basis, norm is finite, j}

\DefEq
{
$\|a\|<\infty$
}
{Schauder basis, norm is finite, 1}

\DefEq
{
\[F=\max(F_1,\epsilon)\]
}
{F F1 epsilon}

\DefEquation
{
\sum_{\gi i=\gi 1}^{\gi n}f^{\gi i}_{\gi j}\,e_{2\cdot\gi i}
}
{linear map, algebra, Schauder basis}

\DefEq
{
$\|a\|=\infty$
}
{Schauder basis, norm is finite, 2}

\DefEq
{
$\|a-b\|$ 
}
{Schauder basis, norm is finite, 3}

\DefEq
{
\[
e'_{\gi i}=\frac 1{\|e_{\gi i}\|}e_{\gi i}
\]
}
{Banach algebra, normal basis, 1}

\DefEq
{
$\|e_{\gi i}\|=1$
}
{Banach algebra, normal basis}

\DefEq
{
a
=a^{\gi i}\,e_{\gi i}
=\lim_{\gi n\rightarrow\gi{\infty}}
\sum_{\gi i=\gi 1}^{\gi n}a^{\gi i}\,e_{\gi i}
}
{Schauder basis, coordinates 0}

\DefEq
{
\[a=a^{\gi i}\,e_{\gi i}\]
}
{a=ae}

\DefEq
{
$\gi n\rightarrow\gi{\infty}$.
}
{Schauder basis, coordinates 3}

\DefEquation
{
\left\|\sum_{\gi i=\gi 1}^{\gi n}a^{\gi i}\,e_{\gi i}\right\|
<\sum_{\gi i=\gi 1}^{\gi n}|a^{\gi i}|\,\|e_{\gi i}\|
=\sum_{\gi i=\gi 1}^{\gi n}|a^{\gi i}|
}
{Schauder basis, coordinates 2}

\DefEq
{
\begin{align}
\left\|\sum_{i=p}^{\infty}a^{\gi i}\,e_{\gi i}\right\|&<\epsilon
\EqLabel{Schauder basis, limit, p}
\\
\left\|\sum_{i=p}^qa^{\gi i}\,e_{\gi i}\right\|&<\epsilon
\EqLabel{Schauder basis, limit, p, q}
\end{align}
}
{Schauder basis, limit}

\DefEquation
{
\begin{matrix}
\displaystyle
\left\|\sum_{i=p}^qa_1^{\gi i}\,e_{\gi i}\right\|
<\frac {\epsilon}2
&
\displaystyle
\left\|\sum_{i=p}^qa_2^{\gi i}\,e_{\gi i}\right\|
<\frac {\epsilon}2
\end{matrix}
}
{Schauder basis, limit p12}

\DefEquation
{
\left\|\sum_{i=p}^qa^{\gi i}\,e_{\gi i}\right\|
<\frac {\epsilon}{|d|}
}
{Schauder basis, limit p22}

\DefEquation
{
\left\|\sum_{i=p}^q(a_1^{\gi i}+a_2^{\gi i})\,e_{\gi i}\right\|
=
\left\|\sum_{i=p}^q(a_1^{\gi i}\,e_{\gi i}
+a_2^{\gi i}\,e_{\gi i})\right\|
<
\left\|\sum_{i=p}^qa_1^{\gi i}\,e_{\gi i}\right\|
+\left\|\sum_{i=p}^qa_2^{\gi i}\,e_{\gi i}\right\|
<\epsilon
}
{Schauder basis, limit p12 1}

\DefEquation
{
\left\|\sum_{i=p}^qd\,a^{\gi i}\,e_{\gi i}\right\|
<
|d|\left\|\sum_{i=p}^qa^{\gi i}\,e_{\gi i}\right\|
<\epsilon
}
{Schauder basis, limit p22 1}

\DefEquation
{
\sum_{i=p}^q(a_1^{\gi i}+a_2^{\gi i})\,e_{\gi i}
}
{Schauder basis, limit p12 2}

\DefEquation
{
\sum_{i=p}^qd\,a^{\gi i}\,e_{\gi i}
}
{Schauder basis, limit p22 2}

\DefEq
{
\symb{\|f\|}0{norm of polymap}
\begin{equation}
\ShowSymbol{norm of polymap}=
\text{sup}\frac{|f(x)|}{\|x\|_1...\|x\|_n}
\EqLabel{norm of polymap, module}
\end{equation}
}
{norm of polymap, module}

\DefEq
{
\symb{\LcDA}1{set continuous linear maps, vector}
}
{set continuous linear maps, module}

\DefEq
{
\symb{\mathcal C(D;A_1,...,A_n;A)}
1{set continuous multivariable maps}
}
{set continuous multivariable maps, module}

\DefEq
{
\symb{\mathcal{LC}(D;A_1,...,A_n;A)}
1{set continuous polylinear maps}
}
{set continuous polylinear maps, module}

\DefEquation
{
n^{-1}dn=nn^{-1}d=d
}
{inverse integer in ring, 1}

%% file: Center.Ring.English.tex
\input{Center.Ring.Eq}

\begin{definition}
\label{definition: center of ring}
Let $D$ be a ring.\footnote{\citeBib{Serge Lang},
page 89.} The set \symb{Z(D)}1{center of ring}
of elements $a\in D$ such that
\ShowEq{center of ring}
for all $x\in D$,
is called
\AddIndex{center of ring $D$}{center of ring}.
\qed
\end{definition}

\begin{theorem}
\label{theorem: center of ring}
The center $Z(D)$ of ring $D$ is subring of ring $D$.
\end{theorem}
\begin{proof}
The statement follows immediately from definition \ref{definition: center of ring}.
\end{proof}

\begin{definition}
Let $D$ be a ring with unit element $e$.\footnote{I
made definition according to definition
from \citeBib{Serge Lang},
pages 89, 90.} The map
\[
l:Z\rightarrow D
\]
such that $l(n)=ne$ is a homomorphism of rings, and its kernel is
an ideal $(n)$, generated by integer $n\ge 0$. We have
canonical injective homomorphism
\[
Z/nZ\rightarrow D
\]
which is an isomorphism between $Z/nZ$ and subring of $D$.
If $nZ$ is prime ideal, then we have two cases.
\begin{itemize}
\item $n=0$. $D$ contains as subring
a ring which isomorphic to $Z$, and which is often identified
with $Z$. In that case, we say that $D$ has
\AddIndex{characteristic $0$}{ring has characteristic 0}.
\item $n=p$ for some prime number $p$. $D$ has
\AddIndex{characteristic $p$}{ring has characteristic p}, and $D$
contains an isomorphic image of $F_p=Z/pZ$.
\end{itemize}
\qed
\end{definition}

\begin{theorem}
\label{theorem: integer subring of ring}
Let $D$ be ring of characteristic $0$ and let $d\in D$.
Then every integer $n\in Z$ commutes with $d$.
\end{theorem}
\begin{proof}
We prove statement by induction.
The statement is evident for $n=0$ and $n=1$.
Let statement be true for $n=k$.
From chain of equation
\[
(k+1)d=kd+d=dk+d=d(k+1)
\]
Evidence of statement for $n=k+1$ follows.
\end{proof}

\begin{theorem}
\label{theorem: integer subring of center}
Let $D$ be ring of characteristic $0$.
Then ring of integers $Z$ is subring of center $Z(D)$ of ring $D$.
\end{theorem}
\begin{proof}
Corollary of theorem \ref{theorem: integer subring of ring}.
\end{proof}

%% file: Center.Ring.Eq.tex

\DefEquation
{
ax=xa
}
{center of ring}

%% file: Biblio.English.tex
\OpenBiblio


\BiblioItem{Einstein: Electrodynamics of Moving Bodies}
{
Albert Einstein,
On the Electrodynamics of Moving Bodies, 1905,
\\
The Principle of Relativity: A Collection of Original
Memoirs on the Special and General Theory of Relativity , 37 - 65,
\\
Courier Dover Publications, 1952; ISBN-13: 978-0486600819
\\
Zur Elektrodynamik der bewegter K\"orper. Ann. Phys., 1905, 17, 891-921. 
}%

\BiblioItem{Einstein: On the Relativity Principle}
{
Albert Einstein,
On the Relativity Principle and the Conclusions Drawn from It, 1907,
\\
The Collected Papers of Albert Einstein, Volume 2:
The Swiss Years: Writings, 1900-1909. English translation. 252 - 311.
\\
Anna Beck, translator; Peter Havas, consultant.
Princeton University Press, 1989; ISBN-13: 9780691085494
\\
\"Uber das Relativit\"atsprinzip und die aus demselben gezogenen Folgerungen. 
Jahrb. d. Radioaktivit\"at u. Elektronik, 1907, 4, 411-462. 
}%

\BiblioItem{Einstein: Foundations of general relativity}
{
Albert Einstein,
Die Grundlage der allgemeinen Relativit\"atstheorie,
Ann. Phys., 1916, {\bf 49}, 769 - 822,\\
Einstein's Annalen Papers: The Complete Collection 1901-1922,
edited by J\"urgen Renn, 517 - 571,\\
Wiley-VCH Verlag GmbH \& Co. KGaA, 2005
}%

\BiblioItem{Einstein: Geometry and Experience}
{
Albert Einstein, Geometry and Experience, (1921)\\
Albert Einstein, Sidelights on Relativity, 25 - 56,\\
Courier Dover Publications, 1983
}%

\BiblioItem{Einstein: Main problems of general relativity}
{
Albert Einstein,
Grundgedanken und Probleme der Relativit\"atstheorie, (1923),\\
Nobelstiftelsen, Les Prix Nobel en 1921 - 1922,
Imprimerie Royale, Stockholm, 1923
}%

\BiblioItem{Einstein: Noneuclidean Geometry and Physics}
{
Albert Einstein,
Nichtenklidische Geometrie in der Physik Neue Rundschan, (1925)
Berlin, S. 16 - 20
}%

\BiblioItem{Einstein: Isaak Newton}
{
Albert Einstein,
Isaak Newton, 1927,
Out of My Later Years, 
Citadel Press, 1995, 219 - 223
}%

\BiblioItem{Einstein: On Science}
{
Albert Einstein,
On Science, 
Cosmic Religion, with Other Opinions and Aphorisms,142 - 146,
New York, 1931, 97 - 103
}%

\BiblioItem{Einstein: Autobiographical Notes}
{
Albert Einstein,
Autobiographical Notes, 1949,\\
Paul A. Schilpp, editor, Albert Einstein: Philosopher-Scientist,
Evanston, 
Illinois, The Library of Living Philosophers, 1949, 1 - 95
}%

\BiblioItem{Feynman 1}
{
Richard Phillips Feynman, Robert B. Leighton, Matthew Linzee Sands.
The Feynman lectures on physics: Volume 1.
Mainly Mechanics, Radiation, and Heat.
Addison\Hyph Wesley, 1965.
}%

\BiblioItem{0538731877}
{
James Shipman, Jerry D. Wilson and Aaron Todd.
Introduction to Physical Science.
Cengage Learning, 2009; ISBN 0538731877.
}%

\BiblioItem{Cite: 104}
{
Cite 104, Source unknown
}%

\BiblioItem{Ghez}
{
Ghez et al.,
The First Measurement of Spectral Lines in a Short-Period Star Bound to the Galaxy's Central Black Hole: A Paradox of Youth,
\href{http://www.journals.uchicago.edu/ApJ/journal/issues/ApJL/v586n2/16990/brief/16990.abstract.html}{ApJL, 586, L127} (2003),
eprint \href{http://arxiv.org/abs/astro-ph/0302299}{arXiv:astro-ph/0302299} (2003)
}%

\BiblioItem{Schodel}
{
R. Sch\"odel et al.,
A star in a 15.2-year orbit around the supermassive black hole at the centre of the Milky Way,
\href{http://www.nature.com/cgi-taf/DynaPage.taf?file=/nature/journal/v419/n6908/abs/nature01121_fs.html}{Nature 419, 694} (2002)
}%

\BiblioItem{Mielke}
{
Eckehard W. Mielke, Affine generalization of the Komar complex of general relativity,
\href{http://prola.aps.org/searchabstract/PRD/v63/i4/e044018}{Phys. Rev. D 63, 044018} (2001)
}%

\BiblioItem{Obukhov}
{
Yu. N. Obukhov and J. G. Pereira, Metric\hyph affine approach to teleparallel gravity,
\href{http://scitation.aip.org/getabs/servlet/GetabsServlet?prog=normal&id=PRVDAQ000067000004044016000001&idtype=cvips&gifs=Yes}
{Phys. Rev. D 67, 044016} (2003),
eprint \href{http://arxiv.org/abs/gr-qc/0212080}{arXiv:gr-qc/0212080} (2002)
}%

\BiblioItem{Sardanashvily}
{
Giovanni Giachetta, Gennadi Sardanashvily, Dirac Equation in Gauge and Affine-Metric Gravitation Theories,
eprint \href{http://arxiv.org/abs/gr-qc/9511035}{arXiv:gr-qc/9511035} (1995)
}%

\BiblioItem{Gauge}
{
Frank Gronwald and Friedrich W. Hehl, On the Gauge Aspects of Gravity, eprint
\href{http://arxiv.org/abs/gr-qc/9602013}{arXiv:gr-qc/9602013} (1996)
}%

\BiblioItem{Neeman}
{
Yuval Neeman, Friedrich W. Hehl, Test Matter in a Spacetime with Nonmetricity, eprint
\href{http://arxiv.org/abs/gr-qc/9604047}{arXiv:gr-qc/9604047} (1996)
}%

\BiblioItem{torsion}
{
F. W. Hehl, P. von der Heyde, G. D. Kerlick, and J. M. Nester,
General relativity with spin and torsion: Foundations and prospects,\\
\href{http://prola.aps.org/abstract/RMP/v48/i3/p393_1}{Rev. Mod. Phys. 48, 393} (1976)
}%

\BiblioItem{Megged}
{
O. Megged, Post-Riemannian Merger of Yang-Mills Interactions with Gravity,
eprint \href{http://arxiv.org/abs/hep-th/0008135}{arXiv:hep-th/0008135} (2001)
}%


\BiblioItem{gr-qc-9604027}
{
Yu.N. Obukhov, E.J. Vlachynsky, W. Esser, R. Tresguerres and F.W. Hehl,
An exact solution of the metric\hyph affine gauge theory with dilation, shear, and spin charges,
eprint \href{http://arxiv.org/abs/gr-qc/9604027}{arXiv:gr-qc/9604027} (1996)
}%

\BiblioItem{4419-7514}
{
Mari\'an Fabian, Petr Habala, Petr H\'ajek, Vicente Montesinos, V\'aclav Zizler.
Banach Space Theory: The Basis for Linear and Nonlinear Analysis.
\\
Springer; New York, 2010; ISBN-13: 978-1441975140
}%

\BiblioItem{Weinberg I}
{
Steven Weinberg.
The Quantum Theory of Fields. Volume I. Foundations.
Cambridge university press, 1995
}%

\BiblioItem{Weinberg II}
{
Steven Weinberg.
The Quantum Theory of Fields. Volume II. Modern applications.
Cambridge university press, 1996
}%

\BiblioItem{Reinhardt}
{
Walter Greiner, Joachim Reinhardt. Field Quantization. Springer.
}%

\BiblioItem{978-3540875604}
{
Walter Greiner, Joachim Reinhardt. Quantum Electrodynamics. Springer, 2009.
}%

\BiblioItem{978-1898563020}
{
H. Robert Mills. Practical Astronomy. Woodhead Publishing, 1994. ISBN-13: 978-1898563020.
}%

\BiblioItem{Landau}
{
L. D. Landau, E. M. Lifshich, The classical theory of fields.
Oxford, New York, Pergamon Press
}%

\BiblioItem{Wheeler}
{
Ignazio Ciufolini, John Wheeler. Gravitation and Inertia.
Princeton university press.
}%

\BiblioItem{Gravitation MTW}
{
Charles W. Misner, Kip S. Thorne, John Archibald Wheeler.
Gravitation.
W. H. Freeman and Company, San Francisco, 1973.
}%

\BiblioItem{Anderson98}
{
J. D. Anderson, P. A. Laing, E. L. Lau, A. S. Liu, M. M. Nieto, and S. G. Turyshev,
Indication, from Pioneer 10/11, Galileo, and Ulysses Data, of an Apparent Anomalous, Weak, Long-Range Acceleration,
\href{http://prola.aps.org/abstract/PRL/v81/i14/p2858_1}{Phys. Rev. Lett. 81, 2858}, (1998),
eprint \href{http://arxiv.org/abs/gr-qc/9808081}{arXiv:gr-qc/9808081} (1998)
}%

\BiblioItem{Anderson02}
{
J. D. Anderson, P. A. Laing, E. L. Lau, A. S. Liu, M. M. Nieto, and S. G. Turyshev,
Study of the anomalous acceleration of Pioneer 10 and 11,
\href{http://prola.aps.org/searchabstract/PRD/v65/i8/e082004}{Phys. Rev. D 65, 082004, 50 pp.}, (2002),
eprint \href{http://arxiv.org/abs/gr-qc/0104064}{arXiv:gr-qc/0104064} (2001)
}%


\BiblioItem{H. Aslaksen}
{
H. Aslaksen.  Quaternionic determinants \textit{Math.
Intelligencer} {\bf 18}(3), pp.57-65, (1996).
}%

\BiblioItem{L. Chen: Definition of determinant}
{
L. Chen, Definition of determinant and Cramer solutions over
quaternion field, \textit{Acta Math. Sinica (N.S.)} {\bf 7},
pp.171-180, (1991).
}%

\BiblioItem{L. Chen: Inverse matrix}
{
L. Chen,
Inverse matrix and properties of double determinant over quaternion
field, \textit{Sci. China, Ser. A} {\bf 34}, pp.528-540, (1991).
}%

\BiblioItem{N. Cohen S. De Leo}
{
N. Cohen, S. De Leo, The quaternionic determinant, \textit{The Electronic Journal Linear
Algebra} {\bf 7}, pp.100-111, (2000).
}%

\BiblioItem{Dyson: Quaternion determinants}
{
F. J. Dyson, Quaternion determinants, \textit{Helvetica Phys.
Acta} {\bf 45}, pp. 289-302, (1972).
}%

\BiblioItem{Melvin Hausner}
{
Melvin Hausner,
A Vector Space Approach to Geometry,
Dover Publications, 1998
}%

\BiblioItem{Serge Lang}
{
Serge Lang,
Algebra, Springer, 2002
}%

\BiblioItem{9780534423230}
{
Charles Lanski.
Concepts In Abstract Algebra.
American Mathematical Soc., 2005, ISBN 9780534423230
}%

\BiblioItem{Burris Sankappanavar}
{
S. Burris, H.P. Sankappanavar,
A Course in Universal Algebra, Springer-Verlag (March, 1982),
\\eprint
\href{http://www.math.uwaterloo.ca/~snburris/htdocs/ualg.html}
{http://www.math.uwaterloo.ca/~snburris/htdocs/ualg.html}
\\(The Millennium Edition)
}%

\BiblioItem{Shilov single 3}
{
G. E. Shilov,
Calculus, Single Variable Functions, Part 3,
Moscow, Nauka, 1970
}%

\BiblioItem{Shilov}
{
G. E. Shilov,
Calculus, Multivariable Functions,
Moscow, Nauka, 1972
}%

\BiblioItem{Kolmogorov Fomin}
{
A. N. Kolmogorov and S. V. Fomin.
Elements of the Theory of Functions and Functional Analysis.
Courier Dover Publication, 1999
}%

\BiblioItem{Lebedev Vorovich}
{
L. P. Lebedev, I. I. Vorovich,
Functional Analysis in Mechanics,
Springer, 2002
}%

\BiblioItem
{Rashevsky}
{
P. K. Rashevsky, Riemann Geometry and Tensor Calculus,\\
Moscow, Nauka, 1967
}%

\BiblioItem
{Kurosh: High Algebra}
{
A. G. Kurosh, High Algebra,
Moscow, Nauka, 1968
}%

\BiblioItem
{Kurosh: General Algebra}
{
A. G. Kurosh, Lectures on General Algebra,
Chelsea Pub Co, 1965 
}%

\BiblioItem
{Sabinin: Smooth Quasigroups}
{
Lev V. Sabinin, Smooth Quasigroups and Loops,
Kluwer Academic Publisher, 1999 
}%

\BiblioItem{Dubrovin Fomenko Novikov part 1}
{
B. A. Dubrovin, A. T. Fomenko, S. P. Novikov,
Modern Geometry - Methods and Applications,\\
Part 1, The Geometry of Surfaces, Transformation Groups, and Fields,\\
Translated by Robert G. Burns,\\
Springer - New York, 1992
}%

\BiblioItem{Korn}
{
Granino A. Korn, Theresa M. Korn,
Mathematical Handbook for Scientists and Engineer,
McGraw-Hill Book Company, New York, San Francisco,
Toronto, London, Sydney, 1968
}%

\BiblioItem{Hocking Young Topology}
{
John G. Hocking, Gail S. Young,
Topology,\\
Courier Dover Publications, 1988
}%

\BiblioItem{Olver: Lie groups to differential equations}
{
Peter J. Olver,
Applications of Lie groups to differential equations,\\
Springer, 2000
}%

\BiblioItem{Tartaglia}
{
Angelo Tartaglia and Matteo Luca Ruggiero,
Angular Momentum Effects in Michelson\Hyph Morley Type Experiments,
Gen.Rel.Grav. 34, 1371-1382 (2002),\\
eprint \href{http://arxiv.org/abs/gr-qc/0110015}{arXiv:gr-qc/0110015} (2001)
}%

\BiblioItem{Tomozawa}
{
Yukio Tomozawa, Speed of Light in Gravitational Fields, eprint
\href{http://arxiv.org/abs/astro-ph/0303047}{arXiv:astro-ph/0303047} (2004)
}%

\BiblioItem{Magueijo}
{
Joao Magueijo,
Covariant and locally Lorentz-invariant varying speed of light theories,
\href{http://prola.aps.org/abstract/PRD/v62/i10/e103521}{Phys. Rev. D 62, 103521} (2000),
eprint \href{http://arxiv.org/abs/gr-qc/0007036}{arXiv:gr-qc/0007036} (2000)
}%

\BiblioItem{Bassett}
{
Bruce A. Bassett, Stefano Liberati, Carmen Molina-Paris, and Matt Visser,
Geometrodynamics of variable-speed-of-light cosmologies,
\href{http://prola.aps.org/abstract/PRD/v62/i10/e103518}{Phys. Rev. D 62}, 103518 (2000),
eprint \href{http://arxiv.org/abs/astro-ph/0001441}{arXiv:astro-ph/0001441} (2000)
}%

\BiblioItem{C.A. Deavours The Quaternion Calculus}
{
C.A. Deavours, The Quaternion Calculus, 
American Mathematical Monthly, {\bf 80} (1973), pp. 995 - 1008
}%

\BiblioItem{Straumann}
{
Lochlainn O'Raifeartaigh and Norbert Straumann,
Gauge theory: Historical origins and some modern developments,
\href{http://prola.aps.org/abstract/RMP/v72/i1/p1_1}{Rev. Mod. Phys. 72, 1} (2000)
}%

\BiblioItem{Lammerzahl}
{
Claus L\"ammerzahl, Mark P. Haugan,
On the interpretation of Michelson\Hyph Morley experiments,
{Phys. Lett. A282 223-229} (2001),\\
eprint \href{http://arxiv.org/abs/gr-qc/0103052}{arXiv:gr-qc/0103052} (2001)
}%

\BiblioItem{0305117}
{
Holger Mueller, Sven Herrmann, Claus Braxmaier, Stephan Schiller, Achim Peters.
Modern Michelson-Morley Experiment using Cryogenic Optical Resonators.
eprint \href{http://arxiv.org/abs/physics/0305117}{arXiv:physics/0305117} (2003)
\\
Phys. Rev. Lett. 91:020401, 2003
}%

\BiblioItem{0706.2031}
{
Holger Mueller, Paul Louis Stanwix, Michael Edmund Tobar,
Eugene Ivanov, Peter Wolf, Sven Herrmann, Alexander Senger,
Evgeny Kovalchuk, Achim Peters.
Relativity tests by complementary rotating Michelson-Morley experiments.
eprint \href{http://arxiv.org/abs/0706.2031}{arXiv:0706.2031 [physics.class-ph]} (2006)
\\
Phys. Rev. Lett. 99:050401, 2007
}%

\BiblioItem{1008.1205}
{
M. Nagel, K. M\"ohle, K. D\"oringshoff, S. Herrmann, A. Senger, E.V. Kovalchuk, A. Peters.
Testing Lorentz Invariance by Comparing Light Propagation in Vacuum and Matter.
eprint \href{http://arxiv.org/abs/1008.1205}{arXiv:1008.1205 [physics.ins-det]} (2010)
}%

\BiblioItem{1109.4897}
{
The OPERA Collaboration.
Measurement of the neutrino velocity with the OPERA detector in the CNGS beam.
eprint \href{http://arxiv.org/abs/1109.4897}{arXiv:1109.4897 [hep-ex]} (2011)
}%

\BiblioItem{Ranada}
{
Antonio F. Ranada,
Pioneer acceleration and variation of light speed: experimental situation,
eprint \href{http://arxiv.org/abs/gr-qc/0402120}{arXiv:gr-qc/0402120} (2004)
}%

\BiblioItem{Gelfand Minlos: rotation and Lorentz groups}
{
Izrail Moiseevich Gelfand, Robert Adolfovich Minlos,
Representations of the rotation and Lorentz groups and their applications;\\
Engl. transl. ed. H. K. Farahat; Transl. by G. Cummins and T. Boddongton;\\
Pergamon Press, 1963
}%

\BiblioItem{math.QA-0208146}
{
I. Gelfand, S. Gelfand, V. Retakh, R. Wilson,
Quasideterminants,\\
eprint \href{http://arxiv.org/abs/math.QA/0208146}{arXiv:math.QA/0208146} (2002)
}%

\BiblioItem{q-alg-9705026}
{
I. Gelfand, V. Retakh,
Quasideterminants, I,\\
eprint \href{http://arxiv.org/abs/q-alg/9705026}{arXiv:q-alg/9705026} (1997)
}%

\BiblioItem{Gelfand Retakh 1991}
{
I. Gelfand and V. Retakh, Determinants of Matrices over Noncommutative Rings, Funct.
Anal. Appl. 25 (1991), no. 2, 91-102
}%

\BiblioItem{Gelfand Retakh 1992}
{
I. Gelfand and V. Retakh, A Theory of Noncommutative Determinants and Characteristic
Functions of Graphs, Funct. Anal. Appl. 26 (1992), no. 4, 1-20
}%

\BiblioItem{hep-th-9407124}
{
I. M. Gelfand, D. Krob, A. Lascoux, B. Leclerc, V.S. Retakh and J.-Y. Thibon,
Noncommutative symmetric functions,\\
eprint \href{http://arxiv.org/abs/hep-th/9407124}{arXiv:hep-th/9407124} (1994)
}%

\BiblioItem{Naimark Shtern: Theory of group representations}
{
Mark Aronovich Naimark, Aleksandr Isaakovich Shtern,
Theory of group representations;\\
Heidelberg, 1982
}%

\BiblioItem{Barut Raczka: Theory of group representations}
{
Asim Orhan Barut; Ryszard R\c{a}czka;
Theory of group representations and applications;\\
World Scientific Publishing Co. Pre. Ltd., 1986
}%

\BiblioItem{Mihalev Pilz: concise handbook of algebra}
{
Aleksandr Vasilevich Mikhalev; G\"{u}nter Pilz;
The concise handbook of algebra;\\
Kluwer Academic Publishers, 2002
}%

\BiblioItem{Shafarevich: Basic notions of algebra}
{
I. R. Shafarevich,
Basic notions of algebra,\\
Translated from the Russian by M. Reid,\\
Springer, 2005
}%

\BiblioItem{Elsgolts: Differential Equations}
{
Lev Elsgolts,
Differential Equations and the Calculus of Variations,\\
University Press of the Pacific, 2003 
}%

\BiblioItem{Baez Huerta: algebra of grand unified theories}
{
John Baez; John Huerta;
The algebra of grand unified theories;\\
Bull. Amer. Math. Soc. {\bf 47} (2010), 483-552
}%

\BiblioItem{J. Fan: Determinants}
{
J. Fan, Determinants and multiplicative functionals
on quaternion matrices, \textit{Linear Algebra and Its
Applications} {\bf 369}, pp. 193-201, (2003).
}%

\BiblioItem{Carl Faith 1}
{
Carl Faith, Algebra: Rings, Modules and Categories I,
Springer - Verlag, Berlin - Heidelberg - New York, 1973
}%

\BiblioItem{Gilson Nimmo Ohta}
{
 C.R.Gilson, J.J.C.Nimmo, Y.Ohta, Quasideterminant solutions of a non-Abelian Hirota-Miwa
 equation, \textit{Journal of Physics A: Mathematical and Theoretical} {\bf 40}(42), pp.
 12607-12617,(2007).
}%

\BiblioItem{Haider Hassan}
{
B. Haider, M. Hassan, Quasideterminant solutions of an integrable chiral model in two
 dimensions, \textit{Journal of Physics A: Mathematical and Theoretical} {\bf 42} (35), art. no.
 355211, (2009).
}%



\BiblioItem{0702447}
{
I.I. Kyrchei, Cramer's rule for quaternion systems of linear equations,
\textit{Journal of Mathematical Sciences} {\bf 155}(6), 839-858, (2008).
 Translated from  \textit{Fundamental and Appl. Math.}
 {\bf 13}(4), pp.67-94, (2007). (in Russian)\\
eprint
\href{http://arxiv.org/abs/math/0702447}{arXiv:math.RA/0702447}
(2007)
}%

\BiblioItem{1004.4380}
{
I.I. Kyrchei, Cramer's rule for some quaternion matrix
    equations,  \textit{Applied Mathematics and Computation} {\bf 217}(5), pp.2024-2030, (2010).\\eprint
\href{http://arxiv.org/abs/1004.4380
}{arXiv:math.RA/arXiv:1004.4380 } (2010)
}%

\BiblioItem{1005.0736}
{
I.I. Kyrchei,Determinantal representations of the Moore-Penrose inverse
 over the quaternion skew field and corresponding Cramer's rules,
 \\
eprint
\href{http://arxiv.org/abs/1005.0736}{arXiv:math.RA/1005.0736}
(2010)
}%

\BiblioItem{0412.391}
{
Aleks Kleyn,
Basis Manifold,
eprint \href{http://arxiv.org/abs/math.DG/0412391}{arXiv:math.DG/0412391} (2007)
}%

\BiblioItem{0405.027}
{
Aleks Kleyn,
Reference Frame in General Relativity,\\
eprint \href{http://arxiv.org/abs/gr-qc/0405027}{arXiv:gr-qc/0405027} (2008)
}%

\BiblioItem{0405.028}
{
Aleks Kleyn, Metric\hyph Affine Manifold,\\
eprint \href{http://arxiv.org/abs/gr-qc/0405028}{arXiv:gr-qc/0405028} (2008)
}%

\BiblioItem{0612.111}
{
Aleks Kleyn,
Biring of Matrices,\\
eprint \href{http://arxiv.org/abs/math.OA/0612111}{arXiv:math.OA/0612111} (2007)
}%

\BiblioItem{0701.238}
{
Aleks Kleyn,
Lectures on Linear Algebra over Division Ring,\\
eprint \href{http://arxiv.org/abs/math.GM/0701238}{arXiv:math.GM/0701238} (2010)
}%

\BiblioItem{0702.561}
{
Aleks Kleyn,
Fibered $\mathfrak{F}$\Hyph Algebra,\\
eprint \href{http://arxiv.org/abs/math.DG/0702561}{arXiv:math.DG/0702561} (2007)
}%

\BiblioItem{math.RA-0501237}
{
Aleks Kleyn,
Vector Space Over Division Ring,\\
eprint \href{http://arxiv.org/abs/math.RA/0412391}{arXiv:math.RA/0501237} (2007)
}%

\BiblioItem{math.RA-0501237v1}
{
Aleks Kleyn,
Module Over Division Ring, version 1,\\
eprint \href{http://arxiv.org/abs/math/0501237v1}{arXiv:math.RA/0501237v1} (2005)
}%

\BiblioItem{0707.2246}
{
Aleks Kleyn,
Fibered Correspondence,\\
eprint \href{http://arxiv.org/abs/0707.2246}{arXiv:0707.2246} (2007)
}%

\BiblioItem{0803.2620}
{
Aleks Kleyn,
Morphism of \Ts Representations,\\
eprint \href{http://arxiv.org/abs/0803.2620}{arXiv:0803.2620} (2008)
}%

\BiblioItem{0803.3276}
{
Aleks Kleyn,
Lorentz Transformation and General Covariance Principle,\\
eprint \href{http://arxiv.org/abs/0803.3276}{arXiv:0803.3276} (2009)
}%

\BiblioItem{0812.4763}
{
Aleks Kleyn,
Introduction into Calculus over Division Ring,\\
eprint \href{http://arxiv.org/abs/0812.4763}{arXiv:0812.4763} (2010)
}%

\BiblioItem{0906.0135}
{
Aleks Kleyn,
Introduction into Geometry over Division Ring,\\
eprint \href{http://arxiv.org/abs/0906.0135}{arXiv:0906.0135} (2010)
}%

\BiblioItem{0909.0855}
{
Aleks Kleyn,
Quaternion Rhapsody,\\
eprint \href{http://arxiv.org/abs/0909.0855}{arXiv:0909.0855} (2010)
}%

\BiblioItem{0912.3315}
{
Aleks Kleyn,
Representation of Universal Algebra,\\
eprint \href{http://arxiv.org/abs/0912.3315}{arXiv:0912.3315} (2009)
}%

\BiblioItem{0912.4061}
{
Aleks Kleyn,
Linear Equation in Finite Dimensional Algebra,\\
eprint \href{http://arxiv.org/abs/0912.4061}{arXiv:0912.4061} (2010)
}%

\BiblioItem{1001.4852}
{
Aleks Kleyn,
The Matrix of Linear Maps,\\
eprint \href{http://arxiv.org/abs/1001.4852}{arXiv:1001.4852} (2010)
}%

\BiblioItem{1003.1544}
{
Aleks Kleyn,
Linear Maps of Free Algebra,\\
eprint \href{http://arxiv.org/abs/1003.1544}{arXiv:1003.1544} (2010)
}%

\BiblioItem{1006.2597}
{
Aleks Kleyn,
The G\^ateaux Derivative and Integral over Banach Algebra,\\
eprint \href{http://arxiv.org/abs/1006.2597}{arXiv:1006.2597} (2010)
}%

\BiblioItem{1011.3102}
{
Aleks Kleyn,
Polylinear Map of Free Algebra,\\
eprint \href{http://arxiv.org/abs/1011.3102}{arXiv:1011.3102} (2010)
}%

\BiblioItem{1104.5197}
{
Aleks Kleyn,
$C^*$-Rhapsody,\\
eprint \href{http://arxiv.org/abs/1104.5197}{arXiv:1104.5197} (2011)
}%

\BiblioItem{1105.4307}
{
Aleks Kleyn,
Algebra with Conjugation,\\
eprint \href{http://arxiv.org/abs/1105.4307}{arXiv:1105.4307} (2011)
}%

\BiblioItem{1107.1139}
{
Aleks Kleyn,
Linear Maps of Quaternion Algebra,\\
eprint \href{http://arxiv.org/abs/1107.1139}{arXiv:1107.1139} (2011)
}%

\BiblioItem{1107.5037}
{
Aleks Kleyn,
Orthogonal Basis and Motion in Finsler Geometry,\\
eprint \href{http://arxiv.org/abs/1107.5037}{arXiv:1107.5037} (2011)
}%

\BiblioItem{1202.6021}
{
Aleks Kleyn,
Maps of Conjugation of Quaternion Algebra,\\
eprint \href{http://arxiv.org/abs/1202.6021}{arXiv:1202.6021} (2012)
}%

\BiblioItem{8433-5163}
{
Aleks Kleyn,
Linear Maps of Free Algebra: First Steps in Noncommutative Linear Algebra,\\
Lambert Academic Publishing, 2010
}%

\BiblioItem{8443-0072}
{
Aleks Kleyn,
Representation Theory: Representation of Universal Algebra,\\
Lambert Academic Publishing, 2011
}%

\BiblioItem{Lauve: Quantum coordinates}
{
A. Lauve, Quantum- and quasi-Plucker coordinates,
\textit{Journal of Algebra} {\bf 296}(2), pp.440-461,
(2006).
}%

\BiblioItem{Lewis D. W. Quaternion algebras}
{
Lewis D. W. Quaternion algebras and the algebraic legacy
of Hamilton's quaternions, \textit{Irish Math. Soc. Bulletin} {\bf
57}, pp. 41-64, (2006).
}%

\BiblioItem{0812.2865}
{
Jos\'e Miguel Figueroa-O'Farrill,
Three lectures on 3-algebras,
eprint \href{http://arxiv.org/abs/0812.2865}{arXiv:0812.2865} (2008)
}%

\BiblioItem{1202.4546}
{
Ming-Liang Hu,
Disentanglement, Bell-nonlocality violation
and teleportation capacity of the decaying tripartite states,
eprint \href{http://arxiv.org/abs/1202.4546}{arXiv:1202.4546} (2012)
}%

\BiblioItem{1203.1629}
{
Borivoje Dakic, Yannick Ole Lipp, Xiaosong Ma, Martin Ringbauer,
Sebastian Kropatschek, Stefanie Barz, Tomasz Paterek, Vlatko Vedral,
Anton Zeilinger, Caslav Brukner, Philip Walther,
Quantum Discord as Optimal Resource for Quantum Communication,
eprint \href{http://arxiv.org/abs/1203.1629}{arXiv:1203.1629} (2012)
}%

\BiblioItem{Li Nimmo: Darboux transformations}
{
C.X.Li, J.J.C. Nimmo, Darboux transformations for a twisted
derivation and quasideterminant solutions to the super KdV
equation, \textit{Proceedings of the Royal Society A:
Mathematical, Physical and Engineering Sciences} {\bf 466} (2120),
pp. 2471-2493, (2010)
}%

\BiblioItem{Schiebold: Cauchy-type determinants}
{
C. Schiebold, Cauchy-type determinants and integrable
systems, \textit{Linear Algebra and Its Applications} {\bf 433}
(2), pp. 447-475, (2010)
}%

\BiblioItem{Suzuki: Noncommutative spectral decomposition}
{
T. Suzuki, Noncommutative
spectral decomposition with qua\-si\-de\-ter\-mi\-nant, \textit{Advances in
Mathematics} {\bf 217}(5), pp. 2141-2158, (2008)
}%

\BiblioItem{1105.3456}
{
C. W. F. Everitt, D. B. DeBra, B. W. Parkinson, J. P. Turneaure, J. W. Conklin,
M. I. Heifetz, G. M. Keiser, A. S. Silbergleit, T. Holmes, J. Kolodziejczak,
M. Al-Meshari, J. C. Mester, B. Muhlfelder, V. Solomonik, K. Stahl, P. Worden,
W. Bencze, S. Buchman, B. Clarke, A. Al-Jadaan, H. Al-Jibreen, J. Li, J. A. Lipa,
J. M. Lockhart, B. Al-Suwaidan, M. Taber, S. Wang,\\
Gravity Probe B: Final Results of a Space Experiment to Test General Relativity,\\
eprint \href{http://arxiv.org/abs/1105.3456}{arXiv:1105.3456[gr-qc]} (2011)
}%

\BiblioItem{0009305}
{
G. S. Asanov.
Can Neutrinos and High-Energy Particles Test Finsler Metric of Space-Time?\\
eprint \href{http://arxiv.org/abs/hep-ph/0009305}{arXiv:hep-ph/0009305} (2000)
}%

\BiblioItem{Asanov 2004}
{
G. S. Asanov.
Finsleroid - space supplemented by angle and scalar product.\\
Hypercomplex Numbers in Geometry and Physics, {\bf 1}, 2004, p. 40 - 62
}%

\BiblioItem{1004.3007}
{
Sergiu I. Vacaru,
Principles of Einstein-Finsler Gravity and Perspectives in Modern Cosmology,\\
eprint \href{http://arxiv.org/abs/1004.3007}{arXiv:1004.3007[math-ph]} (2010)
}%

\BiblioItem{1012.4148}
{
Sergiu I. Vacaru.
Principles of Einstein-Finsler Gravity and Cosmology.\\
eprint \href{http://arxiv.org/abs/1012.4148}{arXiv:1012.4148[physics.gen-ph]} (2010)
}%

\BiblioItem{1112.5641}
{
Christian Pfeifer, Mattias N.R. Wohlfarth.
Finsler geometric extension of Einstein gravity.\\
eprint \href{http://arxiv.org/abs/1112.5641}{arXiv:1112.5641[gr-qc]} (2011)
}%

\BiblioItem{0711.0056}
{
Zhe Chang, Xin Li.
Lorentz Invariance Violation and Symmetry in Randers\Hyph Finsler Spaces.\\
eprint \href{http://arxiv.org/abs/0711.0056}{arXiv:0711.0056[hep-th]} (2011)
}%

\BiblioItem{Rund Finsler geometry}
{
Hanno Rund,
The differential geometry of Finsler spaces.
\\
Springer - Verlag, Berlin - G\"ottingen - Heidelberg, 1959
}%

\BiblioItem{Smirnov vol 1}
{
V. I. Smirnov,
A Course of Higher Mathematics, volume I.
\\
Translated by D. E. Brown.
\\
Translation, edited and additions made by I. N. Sneddon.
\\
Pergamon Press, Addison-Wesley Publishing Company, 1964
}%

\BiblioItem{Beem Dostoglou Ehrlich}
{
John K. Beem, Stamatis A. Dostoglou, Paul E. Ehrlich,
Advances in differential geometry and general relativity.
\\
American Mathematical Society, 2004
}%

\BiblioItem{978-0719033414}
{
Malcolm Pemberton, Nicholas Rau,
Mathematics for economists: an introductory textbook.
\\
Manchester University Press, November 2001; ISBN-13: 978-0719033414
}%

\BiblioItem{0 521 59180 5}
{
Cyrus D. Cantrell,
Modern mathematical methods for physicists and engineers.
\\
Cambridge University Press, 2000
}%

\BiblioItem{Arveson spectral theory}
{
William Arveson,
A short course on spectral theory.
\\
Springer - Verlag, New York, 2002
}%

\BiblioItem{Robert Hermann}
{
Robert Hermann,
Topics in the mathematics of quantum mechanics.
\\
Math Sci Press, 1973
}%

\BiblioItem{9705.009}
{
John C. Baez,
An Introduction to n-Categories,\\
eprint \href{http://arxiv.org/abs/q-alg/9705009}{arXiv:q-alg/9705009} (1997)
}%

\BiblioItem{0105.155}
{
John C. Baez,
The Octonions,\\
eprint \href{http://arxiv.org/abs/math.RA/0105155}{arXiv:math.RA/0105155} (2002)
}%

\BiblioItem{John Baez: Math Blogs}
{
John C. Baez,
What do mathematicians need to know about blogging?,\\
Notices of the American Mathematical Society,
(2010), 3, {\bf 57}, 333,\\
\url{http://www.ams.org/notices/201003/rtx100300333p.pdf}
}%

\BiblioItem{Tolstoi about Anna Karenina}
{
Tolstoi about Anna Karenina,
in book A Karenina Companion, by C. J. G. Turner,
published by Wilfrid Laurier University Press (August 1993)
}%

\BiblioItem
{Cohn: Universal Algebra}
{
Paul M. Cohn,
Universal Algebra,
Springer, 1981
}%

\BiblioItem
{Maunder: Algebraic Topology}
{
C. R. F. Maunder,
Algebraic Topology,
Dover Publications, Inc, Mineola, New York, 1996
}%

\BiblioItem{Pommaret: Partial Differential Equations}
{
J.-F. Pommaret,
Partial Differential Equations and Group Theory,
Springer, 1994
}%

\BiblioItem{Bourbaki: Set Theory}
{
N. Bourbaki,
Theory of sets,
Springer, 2004
}%

\BiblioItem{Bourbaki: Algebra 1}
{
N. Bourbaki,
Algebra 1,
Springer, 2004
}%

\BiblioItem
{Bourbaki: General Topology 1}
{
N. Bourbaki,
General Topology, Chapters 1 - 4,
Springer, 1989
}

\BiblioItem{Bourbaki: General Topology: Chapter 5 - 10}
{
N. Bourbaki,
General Topology, Chapters 5 - 10,
Springer, 1989
}

\BiblioItem{Bourbaki: Topological Vector Space}
{
N. Bourbaki,
Topological Vector Spaces, Chapters 1 - 5,
Transl. by H. G. Eggleston $\&$ S. Madan,
Springer, 2003
}

\BiblioItem{Bourbaki: Coxeter Group Lie}
{
N. Bourbaki,
Lie Groups and Lie Algebras, Chapters 4 - 6,
Translator Andrew Pressley,
Springer, 2002
}

\BiblioItem{Bourbaki: Real Group Lie}
{
N. Bourbaki,
Lie Groups and Lie Algebras, Chapters 7 - 9,
Translator Andrew Pressley,
Springer, 2005
}

\BiblioItem{Shabat: Complex Analysis}
{
Shabat B. V.,
Introduction to Complex Analysis,
\\ \url{http://www.math.uchicago.edu/~ryzhik/shabat-all.pdf},
\\Translated from Russian by L.Ryzhik, 2003
(Moscow, Nauka, 1969)
}

\BiblioItem{Pontryagin: Topological Group}
{
L. S. Pontryagin,
Selected Works, Volume Two, Topological Groups,
Gordon and Breach Science Publishers, 1986
}

\BiblioItem
{Eisenhart: Riemannian Geometry}
{
Eisenhart,
Riemannian Geometry,
Princeton University Press, Princeton, 1949
}

\BiblioItem
{Eisenhart: Continuous Groups of Transformations}
{
Eisenhart,
Continuous Groups of Transformations,
Dover Publications, New York, 1961
}

\BiblioItem
{Condon Odabasi}
{
Edward Uhler Condon, Halis Odabasi,
Atomic Structure,
CUP Archive, 1980
}

\BiblioItem{Postnikov: Differential Geometry}
{
Postnikov M. M.,
Geometry IV: Differential geometry,
Moscow, Nauka, 1983
}

\BiblioItem{Fihtengolts: Calculus volume 1}
{
Fihtengolts G. M.,
Differential and Integral Calculus Course, volume 1,
Moscow, Nauka, 1969
}

\BiblioItem{Fihtengolts: Calculus volume 2}
{
Fihtengolts G. M.,
Differential and Integral Calculus Course, volume 2,
Moscow, Nauka, 1969
}

\BiblioItem{Hatcher: Algebraic Topology}
{
Allen Hatcher,
Algebraic Topology,
Cambridge University Press, 2002
}

\BiblioItem{geometry of differential equations}
{
Vinogradov, A. M., Krasil'shchik, I. S., and Lychagin, V. V.,
Introduction to geometry of nonlinear differential equations,
Nauka, Moscow, 1986
}

\BiblioItem{cohomological analysis}
{
A. M. Vinogradov,
Cohomological Analysis of Partial Differential Equations
and Secondary Calculus,
American Mathematical Society, 2001
}

\BiblioItem{0801.1734}
{
Brandon S. DiNunno, Richard A. Matzner,
The Volume Inside a Black Hole,\\
eprint \href{http://arxiv.org/abs/0801.1734v1}{arXiv:0801.1734v1} (2008)
}

\BiblioItem{0702.447}
{
Ivan Kyrchei,
Cramer's rule for some quaternion matrix equations,\\
eprint \href{http://arxiv.org/abs/math/0702447}{arXiv:math.RA/0702447} (2007)
}

\BiblioItem{Izrail M. Gelfand: Quaternion Groups}
{
I. M. Gelfand, M. I. Graev,
Representation of Quaternion Groups over Localy Compact and
Functional Fields,\\
Funct. Anal. Appl. {\bf 2} (1968) 19 - 33;\\
Izrail Moiseevich Gelfand, Semen Grigorevich Gindikin,\\
Izrail M. Gelfand: Collected Papers, volume II, 435 - 449,\\
Springer, 1989
}

\BiblioItem{Richard D. Schafer}
{
Richard D. Schafer,
An Introduction to Nonassociative Algebras,
Dover Publications, Inc., New York, 1995
}

\BiblioItem{Bamberg Sternberg}
{
Paul Bamberg, Shlomo Sternberg,
A course in mathematics for students of physics,
Cambridge University Press, 1991
}

\BiblioItem{Conway Smith}
{
John Horton Conway, Derek Alan Smith,
On quaternions and octonions: their geometry, arithmetic, and symmetry,
A K Peters, Natick, Massachussets, 2003
}

\BiblioItem{Fueter}
{
Fueter, R.
Die Funktionentheorie der Differentialgleichungen $\Delta u = 0$ und
$\Delta \Delta u = 0$ mit vier reellen Variablen.
Comment. Math. Helv. {\bf 7} (1935), 307-330
}

\BiblioItem{Sudbery Quaternionic Analysis}
{
A. Sudbery,
Quaternionic Analysis,
Math. Proc. Camb. Phil. Soc. (1979), {\bf 85}, 199 - 225
}

\BiblioItem{0902.4771}
{
Fabrizio Colombo, Graziano Gentili, Irene Sabadini,
A Cauchy kernel for slice regular functions,\\
eprint \href{http://arxiv.org/abs/0902.4771v1}{arXiv:0902.4771v1} (2009)
}

\BiblioItem{Vadim Komkov}
{
Vadim Komkov,
Variational Principles of Continuum Mechanics with Engineering Applications: Critical Points Theory,\\
Springer, 1986
}

\BiblioItem{Alain Connes 1994}
{
Alain Connes,
Noncommutative Geometry,\\
Academic Press, 1994
}

\BiblioItem{Hamilton papers 3}
{
Sir William Rowan Hamilton,
The Mathematical Papers, Vol. III, Algebra,\\
Cambridge at the University Press, 1967
}

\BiblioItem{Hamilton Elements of Quaternions 1}
{
Sir William Rowan Hamilton,
Elements of Quaternions, Volume I,\\
Longmans, Green, and Co., London, New York, and Bombay, 1899
}

\BiblioItem{Cartan geometry in reper}
{
Elie Cartan, Vladislav V. Goldberg, Serge\u{i} Pavlovich Finikov,\\
Riemannian geometry in an orthogonal frame:
from lectures delivered by Elie Cartan at the Sorbonne in 1926-1927,\\
translated by Vladislav V. Goldberg,\\
World Scientific, 2001
}

\BiblioItem{Cartan differential form}
{
Henri Cartan.
Differential calculus. Differential forms.\\
Moscow. Mir, 1971
}

\BiblioItem{Arnautov Glavatsky Mikhalev}
{
V. I. Arnautov, S. T. Glavatsky, A. V. Mikhalev,\\
Introduction to the theory of topological rings and modules,
Volume 1995,\\
Marcel Dekker, Inc, 1996
}

\BiblioItem{Moore Yaqub}
{
Hal G. Moore, Adil Yaqub,
A first course in linear algebra with applications,
Edition 3, Academic Press, 1998 
}

\BiblioItem{math.CV-0405471}
{
S. V. Ludkovsky,
Differentiable functions of Cayley-Dickson numbers,\\
eprint \href{http://arxiv.org/abs/math.CV/0405471}{arXiv:math.CV/0405471} (2004)
}%

\BiblioItem{W.Bertram H.Glockner K.Neeb}
{
W.Bertram, H.Glockner, K.Neeb,
Differential Calculus over General Base Fields and Rings,
Expositiones Mathematicae (2004), Volume 22, Issue 3, Pages 213-282
}

\CloseBiblio

%% file: Index.English.tex
\OpenIndex
\SetIndexSpace%
\Index
   {$1$-\rcd form}%
   {1-rcd form, vector spaces}%
\SetIndexSpace%
\Index
   {$2$\Hyph ary fibered relation}%
   {2 ary fibered relation}%
\SetIndexSpace%
\Index
   {$A\CRcirc$\Hyph basis for module}%
   {A CRcirc basis, module over algebra}%
\Index
   {$A\CRcirc$\Hyph linearly dependent set of vectors}%
   {CRcirc linearly dependent, Astar module over algebra}%
\Index
   {$A\CRcirc$\Hyph linearly independent set of vectors}%
   {CRcirc linearly independent, Astar module over algebra}%
\Index
   {$A$\Hyph linear mapping of modules}%
   {A linear map of modules}%
\Index
   {$\mathcal A(A)$\Hyph mapping}%
   {A(A) mapping}%
\Index
   {$A$\Hyph module}%
   {module over algebra}%
\Index
   {$A$\Hyph valued function}%
   {A valued function}%
\Index
   {Abelian $\Omega$\Hyph group}%
   {Abelian Omega group}%
\Index
   {absolute value on division ring}%
   {absolute value on division ring}%
\Index
   {\Acr linear mapping of modules}%
   {Acr linear map of modules}%
\Index
   {$A\CRcirc$\Hyph linear combination}%
   {ACRcirc linear combination}%
\Index
   {active representation}%
   {active representation}%
\Index
   {active representation of group $G(f)$ in basis manifold of representation}%
   {active representation in basis manifold}%
\Index
   {active representation of group $G(\Vector f)$ in basis manifold of tower of representations}%
   {active representation in basis manifold, tower of representations}%
\Index
   {active \sT representation}%
   {active representation, vector space}%
\Index
   {active transformation of basis manifold of representation}%
   {active transformation of basis, representation}%
\Index
   {active transformation of basis manifold of tower of representations}%
   {active transformation of basis, tower of representations}%
\Index
   {active transformation on basis manifold}%
   {active transformation}%
\Index
   {active transformation on the set of \rcd bases}%
   {active transformation, vector space}%
\Index
   {affine basis}%
   {Affine Basis}%
\Index
   {affine structure on set}%
   {affine structure on set}%
\Index
   {affine transformation group}%
   {drc affine transformation group}%
\Index
   {affine transformation group}%
   {affine transformation group}%
\Index
   {affine transformation on basis manifold}%
   {affine transformation}%
\Index
   {algebra of fractions of algebra with conjugation}%
   {algebra of fractions of algebra with conjugation}%
\Index
   {algebra of polynomials over $D$\Hyph algebra}%
   {algebra of polynomials over D algebra}%
\Index
   {algebra of rational mappings of algebra}%
   {algebra of rational mappings of algebra}%
\Index
   {algebra over ring}%
   {algebra over ring}%
\Index
   {algebra with conjugation}%
   {algebra with conjugation}%
\Index
   {alternative representation of matrix}%
   {Alternative representation}%
\Index
   {anholonomic coordinate}%
   {anholonomic coordinate}%
\Index
   {anholonomic coordinates of connection}%
   {anholonomic coordinates of connection}%
\Index
   {anholonomic coordinates of vector}%
   {vector anholonomic coordinates}%
\Index
   {anholonomic coordinates on manifold}%
   {anholonomic coordinates on manifold}%
\Index
   {anholonomity object}%
   {anholonomity object}%
\Index
   {antihomomorphism of fibered groups}%
   {antihomomorphism of fibered groups}%
\Index
   {antisymmetric $2$\Hyph ary fibered relation}%
   {antisymmetric 2 ary fibered relation}%
\Index
   {$A\RCstar$\Hyph basis for vector space}%
   {Arc basis, vector space}%
\Index
   {arity of operation}%
   {arity of operation}%
\Index
   {associative $D$\Hyph algebra}%
   {associative D algebra}%
\Index
   {associative law for $A\star$\Hyph linear mappings of vector spaces}%
   {associative law for Astar linear mappings of vector spaces}%
\Index
   {associative law for $A\star$\Hyph module}%
   {associative law, Astar module over algebra}%
\Index
   {associative law for $A\star$\Hyph vector space}%
   {associative law, Astar vector space}%
\Index
   {associative law for covariant \sT representation}%
   {associative law for covariant starT representation}%
\Index
   {associative law for covariant \Ts representation}%
   {associative law for covariant Tstar representation}%
\Index
   {associative law for $D$\Hyph module}%
   {associative law, D module}%
\Index
   {associative law for \Drc linear maps of vector bundles}%
   {associative law for drc linear maps of vector bundles}%
\Index
   {associative law for $\mathcal D\star$\Hyph vector fields}%
   {associative law, Dstar vector fields}%
\Index
   {associative law for $D\star$\Hyph vector space}%
   {associative law, Dstar vector space}%
\Index
   {associative law for \rcd linear maps of vector spaces}%
   {associative law for rcd linear maps of vector spaces}%
\Index
   {associative law for $\star A$\Hyph module}%
   {associative law, starA module over algebra}%
\Index
   {associative law for $\star D$\Hyph vector space}%
   {associative law, starD vector space}%
\Index
   {associative law for twin representations}%
   {associative law for twin representations}%
\Index
   {associative law of composition of fibered correspondences}%
   {associative law, composition of fibered correspondences}%
\Index
   {associator of $D$\Hyph algebra}%
   {associator of algebra}%
\Index
   {$A\star$\Hyph antilinear mapping of algebra with conjugation}%
   {antilinear mapping of algebra with conjugation}%
\Index
   {$A\star$\Hyph linear map of vector spaces}%
   {Astar linear map of vector spaces}%
\Index
   {$A\star$\Hyph vector space}%
   {Astar vector space}%
\Index
   {$A\star$\Hyph module}%
   {Astar-module}%
\Index
   {$A\star$\Hyph product of vector over scalar}%
   {Astar product of vector over scalar, Astar module}%
\Index
   {$A\star$\hyph product of vector over scalar}%
   {Astar product of vector over scalar, vector space}%
\Index
   {auto parallel line}%
   {auto parallel line}%
\Index
   {automorphism of representation of $\Omega$\Hyph algebra}%
   {automorphism of representation}%
\Index
   {automorphism of tower of representations}%
   {automorphism of tower of representations}%
\Index
   {$(^j_i)$\hyph \CR quasideterminant}%
   {j i cr-quasideterminant}%
\Index
   {norm of quaternion}%
   {norm of quaternion}%
\SetIndexSpace%
\Index
   {Banach $D$\Hyph algebra}%
   {Banach algebra}%
\Index
   {Banach $D$\Hyph module}%
   {Banach module}%
\Index
   {base of fibered correspondence}%
   {base of fibered correspondence}%
\Index
   {base of mapping}%
   {base of map}%
\Index
   {basis dual to basis}%
   {dual basis}%
\Index
   {basis for \sups rows \rcd vector space}%
   {basis, c rows rcd vector space}%
\Index
   {basis for $D$\Hyph vector space}%
   {basis, D vector space}%
\Index
   {basis for \rcd vector space}%
   {basis, rcd vector space}%
\Index
   {basis manifold of affine space}%
   {Basis Manifold, Affine Space}%
\Index
   {basis manifold of central affine space}%
   {Basis Manifold, Central Affine Space, division ring}%
\Index
   {basis manifold of central affine space}%
   {Basis Manifold, Central Affine Space}%
\Index
   {basis manifold of Euclid space}%
   {Basis Manifold, Euclid Space}%
\Index
   {basis manifold of Euclid space}%
   {Basis Manifold, Euclid Space, division ring}%
\Index
   {basis manifold of \rcd affine space}%
   {Basis Manifold, rcd Affine Space, division ring}%
\Index
   {basis manifold of \rcd vector space}%
   {basis manifold of rcd vector space}%
\Index
   {basis manifold of representation}%
   {basis manifold representation F algebra}%
\Index
   {basis manifold of tower of representations}%
   {basis manifold tower of representations}%
\Index
   {basis manifold of vector space}%
   {basis manifold of vector space}%
\Index
   {basis of $A$\Hyph module}%
   {basis of A module}%
\Index
   {basis of $A\CRcirc$\Hyph module}%
   {basis of ACRcirc module}%
\Index
   {basis of algebra $\mathcal L(A;A)$}%
   {basis of algebra L(A,A)}%
\Index
   {basis of $\RCcirc A$\Hyph module}%
   {basis of CRcircA module}%
\Index
   {basis of \subs rows \drc vector space}%
   {basis, r rows drc vector space}%
\Index
   {basis of representation}%
   {basis of representation}%
\Index
   {basis of tower of representations}%
   {basis of tower of representations}%
\Index
   {basis of vector space}%
   {Basis}%
\Index
   {basis vector of representation of Lie group}%
   {basis vector of representation of Lie group}%
\Index
   {basis vector of representation of Lie group over algebra $A$}%
   {basis vector of representation of Lie group over algebra A}%
\Index
   {biring}%
   {biring}%
\Index
   {bundle of level $2$}%
   {bundle of level 2}%
\Index
   {bundle of level $n$}%
   {bundle of level n}%
\SetIndexSpace%
\Index
   {\subs row of matrix}%
   {c row}%
\Index
   {\sups rows \rcd vector space}%
   {sups rows rcd vector space}%
\Index
   {$c$\hyph row of matrix}%
   {c-row}%
\Index
   {Cartan connection}%
   {Cartan connection}%
\Index
   {Cartan curvature}%
   {Cartan curvature}%
\Index
   {Cartan derivative}%
   {Cartan derivative}%
\Index
   {Cartan symbol}%
   {Cartan symbol}%
\Index
   {Cartan transport}%
   {Cartan transport}%
\Index
   {Cartesian power $\Bundle A$ of bundle $\Bundle B$}%
   {Cartesian power A of bundle B}%
\Index
   {Cartesian power $A$ of set $B$}%
   {Cartesian power of set}%
\Index
   {Cartesian power $n$ of bundle $\Bundle E$}%
   {Cartesian power n of bundle E}%
\Index
   {Cartesian power $n$ of $\mathfrak{H}$\Hyph algebra}%
   {Cartesian power of algebra}%
\Index
   {category of \drc vector spaces}%
   {category of drc vector spaces}%
\Index
   {category of fibered correspondences over diagonal}%
   {category of fibered correspondences over diagonal}%
\Index
   {category of reduced fibered correspondences}%
   {category of reduced fibered correspondences}%
\Index
   {category of \Ts representations of $\Omega_1$\Hyph algebra $A$}%
   {category of Tstar representations of Omega1 algebra}%
\Index
   {category of \Ts representations of $\Omega_1$\Hyph algebra from category $\mathcal A$}%
   {category of Tstar representations of Omega1 algebra from category}%
\Index
   {Cauchy sequence}%
   {Cauchy sequence}%
\Index
   {center of $D$\Hyph algebra $A$}%
   {center of algebra}%
\Index
   {center of ring $D$}%
   {center of ring}%
\Index
   {central affine basis}%
   {Central Affine Basis, division ring}%
\Index
   {central affine basis}%
   {Central Affine Basis}%
\Index
   {column determinant}%
   {column determinant}%
\Index
   {column vector}%
   {column vector}%
\Index
   {commutative $D$\Hyph algebra}%
   {commutative D algebra}%
\Index
   {commutative diagram of correspondences}%
   {commutative diagram of correspondences}%
\Index
   {commutator of $D$\Hyph algebra}%
   {commutator of algebra}%
\Index
   {compact\hyph open topology}%
   {compact open topology}%
\Index
   {complete division ring}%
   {complete division ring}%
\Index
   {complete ring}%
   {complete ring}%
\Index
   {complete system of linear partial differential equations}%
   {Complete System of Linear Partial Differential Equations}%
\Index
   {completely integrable system}%
   {completely integrable system}%
\Index
   {component of linear map}%
   {component of linear map}%
\Index
   {component of polylinear map}%
   {component of polylinear map}%
\Index
   {component of the G\^ateaux derivative}%
   {component of Gateaux derivative}%
\Index
   {component of the G\^ateaux derivative of second order}%
   {component of Gateaux derivative of Second Order}%
\Index
   {composition of fibered correspondences}%
   {composition of fibered correspondences}%
\Index
   {composition of reduced fibered correspondences}%
   {composition of reduced fibered correspondences}%
\Index
   {condition of reducibility of products}%
   {condition of reducibility of products}%
\Index
   {conjugate of quaternion $x$}%
   {conjugate of quaternion}%
\Index
   {conjugated $D$\Hyph  module}%
   {conjugated D module}%
\Index
   {conjugation in algebra}%
   {conjugation in algebra}%
\Index
   {conjugation in ring}%
   {conjugation in ring}%
\Index
   {connection coefficients in $D$\Hyph affine space}%
   {connection coefficients, D affine space}%
\Index
   {continuous correspondence}%
   {continuous correspondence}%
\Index
   {continuous map}%
   {continuous map}%
\Index
   {continuous multivariable map}%
   {continuous multivariable map}%
\Index
   {contravariant \sT representation of fibered group}%
   {contravariant starT representation of fibered group}%
\Index
   {contravariant \sT representation of group}%
   {contravariant starT representation of group}%
\Index
   {contravariant \Ts representation of fibered group}%
   {contravariant Tstar representation of fibered group}%
\Index
   {contravariant \Ts representation of group}%
   {contravariant Tstar representation of group}%
\Index
   {coordinate \Drc vector bundle}%
   {coordinate drc vector bundle}%
\Index
   {coordinate isomorphism}%
   {coordinate isomorphism}%
\Index
   {coordinate matrix of set of vectors}%
   {coordinate matrix of set of vectors}%
\Index
   {coordinate matrix of vector}%
   {coordinate matrix of vector}%
\Index
   {coordinate matrix of vector field in \rcD basis}%
   {coordinate matrix of vector field in drc basis}%
\Index
   {coordinate \rcd isomorphism}%
   {coordinate rcd isomorphism}%
\Index
   {coordinate \rcd vector space}%
   {coordinate rcd vector space}%
\Index
   {coordinate reference frame}%
   {coordinate reference frame}%
\Index
   {coordinate representation in $\Omega_2$\Hyph algebra}%
   {coordinate representation, Omega_2 algebra}%
\Index
   {coordinate representation in \rcd vector space}%
   {coordinate representation, rcd vector space}%
\Index
   {coordinate representation in tuple of $\VX\Omega$\Hyph algebras}%
   {coordinate tower of representations, Omega algebra}%
\Index
   {coordinate representation of group in vector space}%
   {coordinate representation, vector space}%
\Index
   {coordinate vector space}%
   {coordinate vector space}%
\Index
   {coordinates of a geometric object in $\Omega_2$\Hyph algebra $M$}%
   {coordinates of geometric object, representation g}%
\Index
   {coordinates of a geometric object in tuple of $\VX\Omega$\Hyph algebras}%
   {coordinates of geometric object, tower of representations g}%
\Index
   {coordinates of basis of representation}%
   {coordinates of basis relative to basis, representation}%
\Index
   {coordinates of element $m$ of representation $f$ relative to set $X$}%
   {coordinates of element relative to set, representation}%
\Index
   {coordinates of endomorphism of representation}%
   {coordinates of endomorphism, representation}%
\Index
   {coordinates of endomorphism of tower of representations}%
   {coordinates of endomorphism, tower of representations}%
\Index
   {coordinates of geometric object}%
   {coordinates of geometric object, vector space}%
\Index
   {coordinates of geometric object in coordinate \rcd vector space}%
   {coordinates of geometric object, coordinate rcd vector space}%
\Index
   {coordinates of geometric object in coordinate representation}%
   {coordinates of geometric object, coordinate vector space}%
\Index
   {coordinates of geometric object in coordinate space of representation}%
   {coordinates of geometric object, coordinate representation}%
\Index
   {coordinates of geometric object in coordinate space of tower of representations}%
   {coordinates of geometric object, coordinate tower of representations}%
\Index
   {coordinates of geometric object in \rcd vector space}%
   {coordinates of geometric object, rcd vector space}%
\Index
   {coordinates of point $A$ of affine space $\overset{\circ}{A}$ relative to basis $(O,\Basis e)$}%
   {coordinates in affine space}%
\Index
   {coordinates of point of \rcd affine space relative to basis}%
   {coordinates in rcd affine space}%
\Index
   {coordinates of representation}%
   {coordinates of representation, drc vector space}%
\Index
   {coordinates of representation}%
   {coordinates of representation}%
\Index
   {coordinates of set of vectors}%
   {coordinates of set of vectors}%
\Index
   {coordinates of vector}%
   {coordinates of vector}%
\Index
   {coordinates of vector field in \Drc basis}%
   {coordinates of vector field in drc basis}%
\Index
   {coordinates of vector relative to Hamel basis}%
   {coordinates of vector, Hamel basis}%
\Index
   {coordinates of vector relative to Schauder basis}%
   {coordinates of vector, Schauder basis}%
\Index
   {correspondence continuous on the set}%
   {correspondence continuous on the set}%
\Index
   {correspondence of homomorphism}%
   {correspondence of homomorphism}%
\Index
   {covariant \sT representation of fibered group}%
   {covariant starT representation of fibered group}%
\Index
   {covariant \sT representation of group}%
   {covariant starT representation of group}%
\Index
   {covariant \Ts representation of fibered group}%
   {covariant Tstar representation of fibered group}%
\Index
   {covariant \Ts representation of group}%
   {covariant Tstar representation of group}%
\Index
   {\CR inverse element of biring}%
   {cr-inverse element}%
\Index
   {\CR matrix group}%
   {cr-matrix group}%
\Index
   {\CR power}%
   {cr power}%
\Index
   {\CR product of matrices}%
   {cr-product of matrices}%
\Index
   {$\CRcirc$\Hyph product of matrices of mappings}%
   {cr product of matrices of mappings}%
\Index
   {\crd vector space}%
   {crd vector space}%
\Index
   {$C^*$\Hyph algebra}%
   {Cstar-algebra}%
\Index
   {curvilinear coordinates of point in affine space}%
   {curvilinear coordinates of point in affine space}%
\Index
   {\subs rows \drc vector space}%
   {subs rows drc vector space}%
\SetIndexSpace%
\Index
   {$D$\Hyph affine space}%
   {d affine space}%
\Index
   {$D$\Hyph linear functional}%
   {D linear functional}%
\Index
   {$D$\Hyph linearly dependent vectors of $D$\Hyph module}%
   {D linearly dependent, module}%
\Index
   {$D$\Hyph linearly independent vectors of $D$\Hyph module}%
   {D linearly independent, module}%
\Index
   {$D$\Hyph affine connection on manifold with affine connections}%
   {D affine connection, affine manifold}%
\Index
   {$D$\Hyph algebra}%
   {D algebra}%
\Index
   {$D$\Hyph basis for module}%
   {D basis, module}%
\Index
   {$D$\Hyph module}%
   {D-module}%
\Index
   {$D$\Hyph module}%
   {D module}%
\Index
   {$D$\Hyph valued variable}%
   {D valued variable}%
\Index
   {$D$\Hyph vector function}%
   {d vector function}%
\Index
   {$D$\Hyph affine connection coefficients on manifold}%
   {D affine connection coefficients, manifold}%
\Index
   {$D$\hyph vector space}%
   {D vector space}%
\Index
   {\dcr vector space}%
   {dcr vector space}%
\Index
   {determinant of matrix}%
   {determinant}%
\Index
   {deviation of trajectories}%
   {deviation of trajectories}%
\Index
   {diagonal in bundle}%
   {diagonal in bundle}%
\Index
   {diagram of correspondences}%
   {diagram of correspondences}%
\Index
   {diagram of representations}%
   {diagram of representations}%
\Index
   {dimension of \rcd vector space}%
   {dimension of vector space}%
\Index
   {direct product of bundles}%
   {Cartesian product of bundles}%
\Index
   {direct product of $D$\Hyph vector spaces}%
   {direct product of D vector spaces}%
\Index
   {direct product of division rings}%
   {direct product of division rings}%
\Index
   {direct product of \rcd vector spaces}%
   {direct product, rcd vector space}%
\Index
   {direct product of representations of fibered group}%
   {direct product of representations of fibered group}%
\Index
   {direct product of representations of group}%
   {direct product of representations of group}%
\Index
   {direct product of total spaces}%
   {Cartesian product of total spaces}%
\Index
   {direct product of \Ts representations of group}%
   {direct product of Tstar representations of group}%
\Index
   {direct sum of representations}%
   {direct sum of representations}%
\Index
   {distributive law for $A\star$\Hyph module}%
   {distributive law, Astar module over algebra}%
\Index
   {distributive law for $A\star$\Hyph vector space}%
   {distributive law, Astar vector space}%
\Index
   {distributive law for $D$\Hyph module}%
   {distributive law, D module}%
\Index
   {distributive law for $\mathcal D\star$\Hyph vector fields}%
   {distributive law, Dstar vector fields}%
\Index
   {distributive law for $D\star$\Hyph vector space}%
   {distributive law, Dstar vector space}%
\Index
   {distributive law for $\star A$\Hyph module}%
   {distributive law, starA module over algebra}%
\Index
   {distributive law for $\star D$\Hyph vector space}%
   {distributive law, starD vector space}%
\Index
   {double determinant}%
   {double determinant}%
\Index
   {\Drc basis for vector  bundle}%
   {drc basis, vector bundle}%
\Index
   {\Drc linear map of vector bundles}%
   {drc linear map of vector bundles}%
\Index
   {\Drc linearly dependent vector fields}%
   {linearly dependent vector fields}%
\Index
   {\Drc linearly independent vector fields}%
   {linearly independent vector fields}%
\Index
   {\drc representation of group}%
   {drc linear representation of group}%
\Index
   {\drc vector}%
   {drc vector}%
\Index
   {\drc vector space}%
   {drc vector space}%
\Index
   {$D\star$\Hyph affine space}%
   {Dstar affine space}%
\Index
   {$D\star$\Hyph antilinear homomorphism}%
   {Dstar antilinear homomorphism}%
\Index
   {$D\star$\Hyph antilinear mapping of ring with conjugation}%
   {antilinear mapping of ring with conjugation}%
\Index
   {\Ds component of coordinates of vector $\Vector r$}%
   {Dstar component of coordinates of vector, D vector space}%
\Index
   {$D\star$\Hyph linear homomorphism}%
   {Dstar linear homomorphism}%
\Index
   {$\mathcal D\star$\Hyph vector bundle}%
   {Dstar vector bundle}%
\Index
   {$\mathcal D\star$\Hyph vector field}%
   {Dstar vector field}%
\Index
   {$D\star$\hyph  vector space}%
   {Dstar vector space}%
\Index
   {$\mathcal D\star$\hyph linear composition of vector fields}%
   {linear composition of vector fields}%
\Index
   {$D\star$\Hyph module}%
   {Dstar-module}%
\Index
   {$\mathcal D\star$\hyph product of vector field over scalar}%
   {Dstar product of vector field over scalar, vector space}%
\Index
   {$D\star$\hyph product of vector over scalar}%
   {Dstar product of vector over scalar, vector space}%
\Index
   {dual space of \rcd vector space}%
   {dual space of rcd vector space}%
\Index
   {duality principle for biring}%
   {duality principle for biring}%
\Index
   {duality principle for biring of matrices}%
   {duality principle for biring of matrices}%
\Index
   {\rcd linear span in vector space}%
   {linear span, vector space}%
\SetIndexSpace%
\Index
   {effective representation of division ring}%
   {effective representation of division ring}%
\Index
   {effective representation of fibered $\Omega$\Hyph algebra}%
   {effective representation of fibered Omega-algebra}%
\Index
   {effective representation of group}%
   {effective representation of group}%
\Index
   {effective representation of $\Omega$\Hyph algebra $A$}%
   {effective representation of algebra}%
\Index
   {effective representation of ring}%
   {effective representation of ring}%
\Index
   {effective \Ts representation of $\Omega$\Hyph algebra $A$}%
   {effective Tstar representation of algebra}%
\Index
   {effective \Ts representation of fibered division ring}%
   {effective representation of fibered division ring}%
\Index
   {effective \Ts representation of fibered group}%
   {effective representation of fibered group}%
\Index
   {effective \Ts representation of group}%
   {effective Tstar representation of group}%
\Index
   {endomorphism of representation of $\Omega$\Hyph algebra}%
   {endomorphism of representation}%
\Index
   {endomorphism of representation regular on generating set $X$}%
   {endomorphism of representation, regular on set}%
\Index
   {endomorphism of representation singular on generating set $X$}%
   {endomorphism of representation, singular on set}%
\Index
   {endomorphism of tower of representations}%
   {endomorphism of tower of representations}%
\Index
   {endomorphism of tower of representations regular on tuple of generating sets}%
   {endomorphism of representation, regular on tuple}%
\Index
   {endomorphism of tower of representations singular on tuple of generating sets}%
   {endomorphism of representation, singular on tuple}%
\Index
   {enhanced Lie group}%
   {enhanced Lie group}%
\Index
   {equivalence generated by representation $f$}%
   {equivalence of representation}%
\Index
   {essential parameters in a set of functions}%
   {essential parameters}%
\Index
   {Euclidean metric on division ring}%
   {Euclidean metric on division ring}%
\Index
   {Euclidean scalar product in $D$\Hyph vector space}%
   {Euclidean scalar product, vector space}%
\Index
   {Euclidean scalar product on division ring}%
   {Euclidean scalar product on division ring}%
\Index
   {expansion of vector relative to basis converges}%
   {expansion converges}%
\Index
   {expansion of vector relative to basis converges normally}%
   {expansion converges normally}%
\Index
   {extended matrix of \drc linear equations}%
   {extended matrix, system of drc linear equations}%
\Index
   {extended matrix of \rcd linear equations}%
   {extended matrix, system of rcd linear equations}%
\Index
   {extension of correspondence}%
   {extension of correspondence}%
\Index
   {extreme line}%
   {extreme line}%
\SetIndexSpace%
\Index
   {fibered coordinate \Drc isomorphism}%
   {fibered coordinate drc isomorphism}%
\Index
   {fibered correspondence from $\Bundle A$ to $\Bundle B$}%
   {fibered correspondence from A to B}%
\Index
   {fibered correspondence in $\Bundle{A}$}%
   {fibered correspondence in A}%
\Index
   {fibered correspondence of homomorphism}%
   {fibered correspondence of homomorphism}%
\Index
   {fibered equivalence}%
   {fibered equivalence}%
\Index
   {fibered group}%
   {fibered group}%
\Index
   {fibered identification morphism}%
   {fibered identification morphism}%
\Index
   {fibered little group}%
   {fibered little group}%
\Index
   {fibered morphism from bundle $\Bundle A$ into $\Bundle B$}%
   {fibered morphism from A into B}%
\Index
   {fibered natural morphism}%
   {fibered natural morphism}%
\Index
   {fibered $\Omega$\Hyph algebra}%
   {fibered Omega-algebra}%
\Index
   {fibered $\Omega$\Hyph subalgebra}%
   {fibered Omega-subalgebra}%
\Index
   {fibered ordering}%
   {fibered ordering}%
\Index
   {fibered preordering}%
   {fibered preordering}%
\Index
   {fibered ring}%
   {fibered ring}%
\Index
   {fibered stability group}%
   {fibered stability group}%
\Index
   {fibered subset}%
   {fibered subset}%
\Index
   {field-strength tensor}%
   {field-strength tensor}%
\Index
   {filter $\mathfrak{F}$ converges to $A$}%
   {filter converges}%
\Index
   {Finsler metric}%
   {Finsler metric}%
\Index
   {Finsler space}%
   {Finsler space}%
\Index
   {Finsler structure}%
   {Finsler structure}%
\Index
   {first Newton law}%
   {First Newton law}%
\Index
   {free $A$\Hyph module}%
   {free A module}%
\Index
   {free algebra over ring}%
   {free algebra over ring}%
\Index
   {free module over ring}%
   {free module over ring}%
\Index
   {free \Ts representation of fibered group}%
   {free representation of fibered group}%
\Index
   {free \Ts representation of group}%
   {free representation of group}%
\Index
   {Frenet transport}%
   {Frenet transport}%
\Index
   {function homogeneous of degree $k$}%
   {function homogeneous}%
\Index
   {function of division ring \Ds differentiable in the Fr\'echet sense}%
   {function Dstar differentiable in Frechet sense, division ring}%
\Index
   {fundamental sequence}%
   {fundamental sequence}%
\SetIndexSpace%
\Index
   {$G$\Hyph reference frame}%
   {G reference frame}%
\Index
   {$G$\Hyph basis of vector space}%
   {G-basis}%
\Index
   {$G$\Hyph coordinates of basis}%
   {G-coordinates}%
\Index
   {$G$\Hyph space}%
   {GSpace}%
\Index
   {the G\^ateaux \crd derivative of map $\Vector f$ of $D$\hyph vector space $\Vector V$ to $D$\hyph vector space $\Vector W$}%
   {Gateaux crd derivative of map, D vector space}%
\Index
   {the G\^ateaux derivative of map}%
   {Gateaux derivative of map}%
\Index
   {the G\^ateaux derivative of order $n$}%
   {Gateaux derivative of Order n}%
\Index
   {the G\^ateaux derivative of second order}%
   {Gateaux derivative of Second Order}%
\Index
   {the G\^ateaux differential of map}%
   {Gateaux differential of map}%
\Index
   {the G\^ateaux \drc derivative of map $\Vector f$ of $D$\Hyph vector space $\Vector V$ to $D$\Hyph vector space $\Vector W$}%
   {Gateaux drc derivative of map, D vector space}%
\Index
   {the G\^ateaux \Ds derivative of map $f$ of division ring $D$}%
   {Gateaux Dstar derivative of map, division ring}%
\Index
   {the G\^ateaux mixed partial derivative}%
   {Gateaux partial derivative of Second Order}%
\Index
   {the G\^ateaux partial derivative}%
   {Gateaux partial derivative}%
\Index
   {the G\^ateaux partial \drc derivative of map $f^b$ with respect to variable $x^a$}%
   {Gateaux partial drc derivative of map with respect to variable, D vector space}%
\Index
   {the G\^ateaux partial \drc derivative of map $f^b$ with respect to variable $x^a$}%
   {Gateaux partial crd derivative of map with respect to variable, D vector space}%
\Index
   {the G\^ateaux \sD derivative of map $f$ of division ring $D$}%
   {Gateaux starD derivative of map, division ring}%
\Index
   {generating set of representation}%
   {generating set of representation}%
\Index
   {generating set of subrepresentation}%
   {generating set of subrepresentation}%
\Index
   {generator of linear map}%
   {generator of linear map, division ring}%
\Index
   {geometric object defined in $\Omega_2$\Hyph algebra $M$}%
   {geometric object, representation g}%
\Index
   {geometric object defined in \rcd vector space}%
   {geometric object, rcd vector space}%
\Index
   {geometric object defined in tuple of $\VX\Omega$\Hyph algebras $\VX A$}%
   {geometric object, tower of representations g}%
\Index
   {geometric object in coordinate representation}%
   {geometric object, coordinate vector space}%
\Index
   {geometric object in coordinate representation defined in $\Omega_2$\Hyph algebra $M$}%
   {geometric object, coordinate representation g}%
\Index
   {geometric object in coordinate representation defined in \rcd vector space}%
   {geometric object, coordinate rcd vector space}%
\Index
   {geometric object in coordinate representation defined in tuple of $\VX\Omega$\Hyph algebras $\VX A$}%
   {geometric object, coordinate tower of representations g}%
\Index
   {geometric object in vector space}%
   {geometric object, vector space}%
\Index
   {geometric object of type $H$}%
   {geometric object of type H, representation g}%
\Index
   {geometric object of type $A$ in vector space}%
   {geometric object of type A, vector space}%
\Index
   {group algebra}%
   {group algebra}%
\Index
   {group of automorphisms of representation}%
   {group of automorphisms of representation}%
\SetIndexSpace%
\Index
   {Hadamard inverse of matrix}%
   {Hadamard inverse of matrix}%
\Index
   {Hamel basis}%
   {Hamel basis}%
\Index
   {hermitian conjugated vector}%
   {hermitian conjugated vector}%
\Index
   {hermitian conjugation in division ring}%
   {hermitian conjugation, division ring}%
\Index
   {hermitian matrix}%
   {hermitian matrix}%
\Index
   {hermitian metric on division ring}%
   {hermitian metric on division ring}%
\Index
   {hermitian scalar product in $D$\Hyph vector space}%
   {hermitian scalar product, vector space}%
\Index
   {hermitian scalar product on division ring}%
   {hermitian scalar product on division ring}%
\Index
   {holonomic coordinates of connection}%
   {holonomic coordinates of connection}%
\Index
   {holonomic coordinates of vector}%
   {vector holonomic coordinates}%
\Index
   {homogeneous bundle of fibered group}%
   {homogeneous bundle of fibered group}%
\Index
   {homogeneous linear geometric object}%
   {homogeneous linear geometric object}%
\Index
   {homogeneous linear representation of Lie group}%
   {homogeneous Linear Representation of Lie Group}%
\Index
   {homogeneous map of degree $k$ over field $F$}%
   {homogeneous map of degree over field, D vector space}%
\Index
   {homogeneous space of group}%
   {homogeneous space of group}%
\Index
   {homomorphism of fibered groups}%
   {homomorphism of fibered groups}%
\Index
   {homomorphism of fibered universal algebras}%
   {homomorphism of fibered universal algebras}%
\SetIndexSpace%
\Index
   {ideal of algebra}%
   {ideal of algebra}%
\Index
   {infinitesimal generator}%
   {infinitesimal generator}%
\Index
   {infinitesimal generators of group Lie}%
   {infinitesimal generators of group Lie}%
\Index
   {invariance principle in \drc vector space}%
   {invariance principle}%
\Index
   {invariance principle in representation of universal algebra}%
   {invariance principle, representation g}%
\Index
   {invariance principle in tower of representations of universal algebras}%
   {invariance principle, tower of representations g}%
\Index
   {invariance principle in vector space}%
   {invariance principle, vector space}%
\Index
   {inverse fibered correspondence}%
   {inverse fibered correspondence}%
\Index
   {inverse reduced fibered correspondence}%
   {inverse reduced fibered correspondence}%
\Index
   {involution in quaternion algebra}%
   {involution, quaternion algebra}%
\Index
   {isomorphism of fibered $\Omega$\Hyph algebras}%
   {isomorphism of fibered Omega-algebras}%
\Index
   {isomorphism of repesentations of $\Omega$\Hyph algebra}%
   {isomorphism of repesentations of Omega algebra}%
\Index
   {isotropic vector}%
   {isotropic vector}%
\SetIndexSpace%
\Index
   {$(^j_i)$\hyph $\RCcirc$\Hyph quasideterminant}%
   {j i RCcirc-quasideterminant}%
\Index
   {Jacobian complete system of differential equations}%
   {Jacobian complete system of differential equations}%
\Index
   {Jacobian complete system of \drv differential equations}%
   {Jacobian complete system of drc differential equations}%
\Index
   {$(ji)$\hyph quasideterminant}%
   {j i quasideterminant}%
\Index
   {the Jacobi\Hyph G\^ateaux matrix of map}%
   {Jacobi Gateaux matrix of map}%
\SetIndexSpace%
\Index
   {kernel of inefficiency of representation of fibered group}%
   {kernel of inefficiency of representation of fibered group}%
\Index
   {kernel of inefficiency of representation of group}%
   {kernel of inefficiency of representation of group}%
\Index
   {kernel of inefficiency of \Ts representation of group $G$}%
   {kernel of inefficiency of Tstar representation of group}%
\Index
   {kernel of linear map}%
   {kernel of linear map}%
\Index
   {Killing equation}%
   {Killing equation}%
\Index
   {Killing equation of second type}%
   {Killing equation second type}%
\Index
   {Killing vector of second type}%
   {Killing vector second type}%
\Index
   {Kronecker symbol}%
   {Kronecker symbol}%
\SetIndexSpace%
\Index
   {left cofactor of entry of matrix}%
   {left cofactor, matrix}%
\Index
   {left defined Lie algebra of Lie group}%
   {left defined Lie algebra}%
\Index
   {left double cofactor of entry of matrix}%
   {left double cofactor}%
\Index
   {left fraction}%
   {left fraction}%
\Index
   {left ideal of algebra}%
   {left ideal of algebra}%
\Index
   {left invariant vector field}%
   {left invariant vector}%
\Index
   {left module over $D$\Hyph algebra $A$}%
   {left module over algebra}%
\Index
   {left module over a ring $D$}%
   {left module over ring}%
\Index
   {left principal ideal}%
   {left principal ideal}%
\Index
   {left shift of $R$\Hyph module}%
   {left shift of module}%
\Index
   {left shift on fibered group}%
   {Tstar shift, fibered group}%
\Index
   {left shift on group}%
   {left shift}%
\Index
   {left shift on group}%
   {left shift, group}%
\Index
   {left structural constant of Lie algebra}%
   {left structural constant of Lie algebra}%
\Index
   {left vector space}%
   {left vector space}%
\Index
   {left-ordered cycle notation of permutation}%
   {left-ordered cycle notation of permutation}%
\Index
   {left-side contravariant representation of group}%
   {left-side contravariant representation of group}%
\Index
   {left-side covariant representation of group}%
   {left-side covariant representation of group}%
\Index
   {left-side representation of fibered $\Omega$\Hyph algebra}%
   {left-side representation of fibered Omega-algebra}%
\Index
   {left-side representation of $\Omega_1$\Hyph algebra $A$ in $\Omega_2$\Hyph algebra $M$}%
   {left-side representation of algebra}%
\Index
   {left-side transformation}%
   {left-side transformation}%
\Index
   {left-side transformation on bundle}%
   {left-side transformation of bundle}%
\Index
   {Lie algebra of Lie group}%
   {algebra Lie group Lie}%
\Index
   {Lie derivative}%
   {Lie derivative}%
\Index
   {Lie derivative of connection}%
   {Lie derivative of connection}%
\Index
   {Lie derivative of metric}%
   {Lie derivative of metric}%
\Index
   {Lie group basic operators}%
   {Lie group basic operators}%
\Index
   {lift of correspondence}%
   {lift of correspondence}%
\Index
   {lift of mapping}%
   {lift of map}%
\Index
   {limit of correspondence with respect to the filter}%
   {limit of correspondence with respect to the filter}%
\Index
   {limit of filter}%
   {limit of filter}%
\Index
   {limit of sequence}%
   {limit of sequence}%
\Index
   {limit set of filter}%
   {limit set of filter}%
\Index
   {linear combination of vectors of $A$\Hyph module}%
   {linear combination in A module}%
\Index
   {linear composition of  vectors}%
   {linear composition of  vectors}%
\Index
   {linear dependent vectors of $A$\Hyph module}%
   {linear dependent vectors, module}%
\Index
   {linear functional on $D$\Hyph module}%
   {linear functional, D module}%
\Index
   {linear geometric object}%
   {linear geometric object}%
\Index
   {linear independent vectors of $A$\Hyph module}%
   {linear independent vectors, module}%
\Index
   {linear map}%
   {linear map}%
\Index
   {linear map of $D$\Hyph vector spaces}%
   {linear map of D vector spaces}%
\Index
   {linear map of $D$\Hyph vector spaces over field $F$}%
   {linear map over field, vector space}%
\Index
   {linear map of division ring generated by map $G$}%
   {linear map generated by map, division ring}%
\Index
   {linear map of \rcd vector spaces}%
   {linear map of rcd vector spaces}%
\Index
   {linear mapping of division ring}%
   {linear mapping of division ring}%
\Index
   {linear representation of group}%
   {linear representation of group}%
\Index
   {linear representation of Lie group}%
   {Linear Representation of Lie Group}%
\Index
   {linear transformation of \rcd affine space}%
   {linear transformation, rcd affine space}%
\Index
   {linearly dependent vectors}%
   {linearly dependent vectors}%
\Index
   {linearly independent vectors}%
   {linearly independent vectors}%
\Index
   {little group}%
   {little group}%
\Index
   {local reference frame}%
   {local reference frame}%
\Index
   {locally compact at point $p$ space}%
   {locally compact at point space}%
\Index
   {locally compact space}%
   {locally compact space}%
\Index
   {Lorentz transformation}%
   {Lorentz transformation}%
\SetIndexSpace%
\Index
   {$m$\Hyph dimensional parallelepiped}%
   {m dimensional parallelepiped}%
\Index
   {$m$\Hyph vector}%
   {m-vector}%
\Index
   {manifold with $D$\Hyph affine connections}%
   {manifold with D- affine connections}%
\Index
   {map continuous with respect to set of arguments}%
   {map continuous with respect to set of arguments}%
\Index
   {map differentiable in the G\^ateaux sense}%
   {map differentiable in Gateaux sense}%
\Index
   {map of $\gi n$ $D$\Hyph valued variables}%
   {map of n D valued variables}%
\Index
   {map of type $G$ on manifold}%
   {map of type G on manifold}%
\Index
   {map polylinear over finite dimensional algebras}%
   {map polylinear over finite dimensional algebras}%
\Index
   {mapping of rings polylinear over commutative ring}%
   {map polylinear over commutative ring, ring}%
\Index
   {mapping space}%
   {mapping space}%
\Index
   {matrix of antilinear homomorphism}%
   {matrix of antilinear homomorphism}%
\Index
   {matrix of bilinear function}%
   {matrix of bilinear function}%
\Index
   {matrix of endomorphisms of $\Omega$\Hyph algebra}%
   {matrix of endomorphisms of Omega algebra}%
\Index
   {matrix of fibered \Drc linear map}%
   {matrix of fibered drc linear map}%
\Index
   {matrix of linear homomorphism}%
   {matrix of linear homomorphism}%
\Index
   {matrix of linear map}%
   {matrix of linear map}%
\Index
   {matrix of linear mappings}%
   {matrix of linear mappings}%
\Index
   {matrix of mappings}%
   {matrix of mappings}%
\Index
   {matrix of quadratic map}%
   {matrix of quadratic map, division ring}%
\Index
   {metric tensor in Minkowski space}%
   {metric tensor, Minkowski space}%
\Index
   {metric-affine manifold}%
   {metric-affine manifold}%
\Index
   {Minkowski space}%
   {Minkowski space, Finsler}%
\Index
   {minor matrix}%
   {minor matrix}%
\Index
   {module over ring}%
   {module over ring}%
\Index
   {morphism from tower of \Ts representations into tower of \Ts representations}%
   {morphism from tower of representations into tower of representations}%
\Index
   {morphism of fibered \Ts representations from $\Bundle F$ into $\Bundle G$}%
   {morphism of fibered representations from f into g}%
\Index
   {morphism of representation $f$}%
   {morphism of representation f}%
\Index
   {morphism of representations from $f$ into $g$}%
   {morphism of representations from f into g}%
\Index
   {morphism of representations of $\Omega_1$\Hyph algebra in $\Omega_2$\Hyph algebra}%
   {morphism of representations of Omega1 algebra in Omega2 algebra}%
\Index
   {morphism of \Ts representations of fibered $\Omega$\Hyph algebra}%
   {morphism of representations of fibered Omega algebra}%
\Index
   {motion of Minkowski space}%
   {motion, Minkowski space}%
\Index
   {movement on basis manifold}%
   {movement transformation}%
\SetIndexSpace%
\Index
   {$n$\Hyph algebra over the ring}%
   {n algebra over ring}%
\Index
   {$n$\Hyph ary fibered relation}%
   {fibered relation}%
\Index
   {nonmetricity}%
   {nonmetricity}%
\Index
   {nonsingular bilinear function}%
   {nonsingular bilinear function}%
\Index
   {nonsingular system of \rcd linear equations}%
   {nonsingular system of linear equations}%
\Index
   {nonsingular tensor}%
   {nonsingular tensor, algebra}%
\Index
   {nonsingular transformation}%
   {nonsingular transformation}%
\Index
   {norm in quaternion algebra}%
   {norm, quaternion algebra}%
\Index
   {norm of functional}%
   {norm of functional}%
\Index
   {norm of map}%
   {norm of map}%
\Index
   {norm of polylinear map}%
   {norm of polymap}%
\Index
   {norm on $D$\Hyph algebra}%
   {norm on D algebra}%
\Index
   {norm on $D$\Hyph vector space}%
   {norm on D vector space}%
\Index
   {norm on $D$\Hyph module}%
   {norm on D module}%
\Index
   {norm on ring}%
   {norm on ring}%
\Index
   {normal basis}%
   {normal basis}%
\Index
   {normed $D$\Hyph algebra}%
   {normed D algebra}%
\Index
   {normed $D$\Hyph module}%
   {normed D module}%
\Index
   {normed $D$\Hyph vector space}%
   {normed D vector space}%
\Index
   {normed ring}%
   {normed ring}%
\Index
   {not complete group}%
   {not complete group}%
\Index
   {not complete $\Omega$\Hyph algebra}%
   {not complete Omega algebra}%
\Index
   {nucleus of $D$\Hyph algebra $A$}%
   {nucleus of algebra}%
\SetIndexSpace%
\Index
   {octonion algebra}%
   {octonion algebra}%
\Index
   {$\Omega$\Hyph group}%
   {Omega group}%
\Index
   {$\Omega$\Hyph linear mapping}%
   {Omega linear map}%
\Index
   {$\Omega_2$\Hyph word of element of representation relative to generating set}%
   {word of element relative to generating set, representation}%
\Index
   {operation on bundle}%
   {operation on bundle}%
\Index
   {opposite algebra to algebra $P$}%
   {opposite algebra}%
\Index
   {opposite fibered preordering}%
   {opposite fibered preordering}%
\Index
   {orbit of linear mapping}%
   {orbit of linear mapping}%
\Index
   {orbit of representation of fibered group}%
   {orbit of representation of fibered group}%
\Index
   {orbit of representation of group}%
   {orbit of representation of group}%
\Index
   {orbit of \Ts representation of group}%
   {orbit of Tstar  representation of group}%
\Index
   {origin of coordinate system of affine space}%
   {origin of coordinate system of affine space}%
\Index
   {origin of coordinate system of $\star D$\Hyph affine space}%
   {origin of coordinate system of starD affine space}%
\Index
   {orthogonal basis in Minkowski space}%
   {orthogonal basis, Minkowski space}%
\Index
   {orthogonality in Minkowski space}%
   {Minkowski orthogonality}%
\Index
   {orthonormal basis in Minkowski space}%
   {orthonormal basis, Minkowski space}%
\Index
   {orthonornal basis}%
   {Orthonornal Basis}%
\Index
   {orthonornal basis}%
   {Orthonornal Basis, division ring}%
\SetIndexSpace%
\Index
   {passive representation of group $G(f)$ in basis manifold of representation}%
   {passive representation in basis manifold}%
\Index
   {parallel shift of \rcd affine space}%
   {parallel shift, rcd affine space}%
\Index
   {parallelogram}%
   {parallelogram}%
\Index
   {partial linear map}%
   {partial linear map}%
\Index
   {passive representation}%
   {passive representation}%
\Index
   {passive representation of group $G(\Vector f)$ in basis manifold of tower of representations}%
   {passive representation in basis manifold, tower of representations}%
\Index
   {passive \sT representation}%
   {passive starT representation}%
\Index
   {passive transformation of the basis manifold of representation}%
   {passive transformation of basis, representation}%
\Index
   {passive transformation of the basis manifold of tower of representations}%
   {passive transformation of basis, tower of representations}%
\Index
   {passive transformation on basis manifold}%
   {passive transformation}%
\Index
   {passive transformation on the set of \rcd bases}%
   {passive transformation, vector space}%
\Index
   {permutability property of trace}%
   {permutability property of trace}%
\Index
   {pfaffian derivative}%
   {pfaffian derivative}%
\Index
   {polylinear map}%
   {polylinear map}%
\Index
   {polylinear skew symmetric map}%
   {polylinear map skew symmetric}%
\Index
   {polylinear symmetric map}%
   {polylinear map symmetric}%
\Index
   {polymorphism of representations}%
   {polymorphism of representations}%
\Index
   {polyvector}%
   {polyvector}%
\Index
   {potential energy}%
   {potential energy}%
\Index
   {principal ideal}%
   {principal ideal}%
\Index
   {product of geometric object and constant}%
   {product of geometric object and constant}%
\Index
   {product of geometric object and constant in vector space}%
   {product of geometric object and constant, vector space}%
\Index
   {product of groups}%
   {product of groups}%
\Index
   {product of map over scalar}%
   {product of map over scalar}%
\Index
   {product of morphisms of representations of universal algebra}%
   {product of morphisms of representations of universal algebra}%
\Index
   {product of morphisms of tower of representations}%
   {product of morphisms of tower of representations}%
\Index
   {product of morphisms of \Ts representations of fibered $\Omega$\Hyph algebra}%
   {product of morphisms of representations of fibered Omega algebra}%
\Index
   {product of objects in category}%
   {product of objects in category}%
\Index
   {projection of bundle $\Bundle E$ along fiber $E$}%
   {projection of bundle along fiber}%
\Index
   {pseudo\Hyph Euclidean metric on division ring}%
   {pseudo-Euclidean metric on division ring}%
\Index
   {pseudo\Hyph Euclidean scalar product in $D$\Hyph vector space}%
   {pseudo-Euclidean scalar product, vector space}%
\Index
   {pseudo-Euclidean scalar product on division ring}%
   {pseudo-Euclidean scalar product on division ring}%
\SetIndexSpace%
\Index
   {quadratic form in division ring}%
   {quadratic form, division ring}%
\Index
   {quadratic map of division ring}%
   {Quadratic Map of Division Ring}%
\Index
   {quasi affine transformation on basis manifold}%
   {quasi affine transformation}%
\Index
   {quasi affine transformation on basis manifold}%
   {quasi affine drc transformation}%
\Index
   {quasi movement on basis manifold}%
   {quasi movement, division ring}%
\Index
   {quasi movement on basis manifold}%
   {quasi movement}%
\Index
   {quasiclosed ring of mappings}%
   {quasiclosed ring of mappings}%
\Index
   {quasideterminant}%
   {quasideterminant definition}%
\Index
   {quasimotion of Minkowski space}%
   {Quasimotion, Minkowski space}%
\Index
   {quaternion algebra}%
   {quaternion algebra}%
\Index
   {quaternion algebra $E$ over the field $F$}%
   {quaternion algebra over the field}%
\Index
   {quotient bundle}%
   {quotient bundle}%
\SetIndexSpace%
\Index
   {\drc basis dual to \rcd basis of vector space}%
   {basis dual to basis, rcd vector space}%
\Index
   {$(^j_i)$\hyph \RC quasideterminant}%
   {j i RC-quasideterminant}%
\Index
   {\sups row of matrix}%
   {r row}%
\Index
   {$R$\Hyph module}%
   {R- module}%
\Index
   {$r$\hyph row of matrix}%
   {r-row}%
\Index
   {rank of Hermitian matrix by principal minors}%
   {rank of Hermitian matrix by principal minors}%
\Index
   {rank of quadratic map of division ring}%
   {rank of quadratic map, division ring}%
\Index
   {\RC inverse element of biring}%
   {rc-inverse element}%
\Index
   {\RC major minor matrix}%
   {RC-major minor}%
\Index
   {\RC matrix group}%
   {rc-matrix group}%
\Index
   {\RC nonsingular matrix}%
   {RC nonsingular matrix}%
\Index
   {\RC power}%
   {rc power}%
\Index
   {\RC product of matrices}%
   {rc-product of matrices}%
\Index
   {$\RCcirc$\Hyph product of matrices of mappings}%
   {rc product of matrices of mappings}%
\Index
   {\RC quasideterminant}%
   {RC-quasideterminant}%
\Index
   {\RC rank of matrix}%
   {rc-rank of matrix}%
\Index
   {\RC singular matrix}%
   {RC singular matrix}%
\Index
   {$\RCcirc A$\Hyph basis for module}%
   {RCcirc A basis, module over algebra}%
\Index
   {$\RCcirc A$\Hyph linearly dependent set of vectors}%
   {RCcirc linearly dependent, starA module over algebra}%
\Index
   {$\RCcirc A$\Hyph linearly independent set of vectors}%
   {RCcirc linearly independent, starA module over algebra}%
\Index
   {$\RCcirc$\Hyph nonsingular matrix of $\mathcal A(A)$\Hyph mappings}%
   {RCcirc nonsingular matrix of A(A) mappings}%
\Index
   {$\RCcirc$\Hyph nonsingular matrix of endomorphisms}%
   {RCcirc nonsingular matrix of endomorphisms}%
\Index
   {$\RCcirc$\Hyph nonsingular system of additive equations}%
   {RCcirc nonsingular system of additive equations}%
\Index
   {$\RCcirc$\Hyph quasideterminant}%
   {RCcirc-quasideterminant definition}%
\Index
   {$\RCcirc$\Hyph singular matrix of $\mathcal A(A)$\Hyph mappings}%
   {RCcirc singular matrix of A(A) mappings}%
\Index
   {$\RCcirc$\Hyph singular matrix of endomorphisms}%
   {RCcirc singular matrix of endomorphisms}%
\Index
   {$\RCcirc A$\Hyph linear combination}%
   {RCcircA linear combination}%
\Index
   {\rcd affine basis}%
   {rcd affine basis, division ring}%
\Index
   {\rcd affine plane}%
   {rcd affine plane}%
\Index
   {\rcd affine space}%
   {rcd affine space}%
\Index
   {\rcd affine transformation}%
   {rcd affine transformation}%
\Index
   {\rcd automorphism of vector space}%
   {automorphism of vector space}%
\Index
   {\rcd isomorphism of vector spaces}%
   {isomorphism of vector spaces}%
\Index
   {\rcd linear \Ts representation of group}%
   {rcd linear Tstar representation of group}%
\Index
   {\rcd vector}%
   {rcd vector}%
\Index
   {\rcd vector space}%
   {rcd vector space}%
\Index
   {reduced Cartesian product of bundles}%
   {reduced Cartesian product of bundles}%
\Index
   {reduced Cartesian product of total spaces}%
   {reduced Cartesian product of total spaces}%
\Index
   {reduced fibered correspondence from $\Bundle{A}$ to $\Bundle B$}%
   {reduced fibered correspondence from A to B}%
\Index
   {reduced fibered correspondence in $\Bundle{A}$}%
   {reduced fibered correspondence in A}%
\Index
   {reduced polymorphism of representations}%
   {reduced polymorphism of representations}%
\Index
   {reducible biring}%
   {reducible biring}%
\Index
   {reference frame in event space}%
   {reference frame in event space}%
\Index
   {reference frame manifold}%
   {reference frame manifold}%
\Index
   {reflexive $2$\Hyph ary fibered relation}%
   {reflexive 2 ary fibered relation}%
\Index
   {regular endomorphism of representation}%
   {regular endomorphism of representation}%
\Index
   {regular endomorphism of tower of representations}%
   {regular endomorphism of tower of representations}%
\Index
   {regular quadratic map in division ring}%
   {regular quadratic map, division ring}%
\Index
   {representation of group}%
   {representation of group}%
\Index
   {representation of $\Omega$\Hyph algebra in representation}%
   {representation of Omega algebra in representation}%
\Index
   {representation of $\Omega$\Hyph algebra in tower of representations}%
   {representation of Omega algebra in tower of representations}%
\Index
   {representation of $\Omega$\Hyph algebra $A$ in category $\mathcal B$}%
   {representation of Omega algebra in category}%
\Index
   {representation of $\Omega_1$\Hyph algebra $A$ in $\Omega_2$\Hyph algebra $M$}%
   {representation of algebra}%
\Index
   {representative of geometric object in \drc vector space}%
   {representative of geometric object, drc vector space}%
\Index
   {representative of geometric object in $\Omega_2$\Hyph algebra}%
   {representative of geometric object, representation g}%
\Index
   {representative of geometric object in \rcd vector space}%
   {representative of geometric object, rcd vector space}%
\Index
   {representative of geometric object in tuple of $\VX\Omega$\Hyph algebras}%
   {representative of geometric object, tower of representations g}%
\Index
   {representative of geometric object in vector space}%
   {representative of geometric object, vector space}%
\Index
   {restriction of correspondence $\Phi$ to set $C$}%
   {restriction of correspondence}%
\Index
   {right cofactor of entry of matrix}%
   {right cofactor, matrix}%
\Index
   {right defined Lie algebra of Lie group}%
   {right defined Lie algebra}%
\Index
   {right double cofactor of entry of matrix}%
   {right double cofactor}%
\Index
   {right fraction}%
   {right fraction}%
\Index
   {right ideal of algebra}%
   {right ideal of algebra}%
\Index
   {right invariant vector field}%
   {right invariant vector}%
\Index
   {right module over $D$\Hyph algebra $A$}%
   {right module over algebra}%
\Index
   {right module over a ring $D$}%
   {right module over ring}%
\Index
   {right principal ideal}%
   {right principal ideal}%
\Index
   {right shift on group}%
   {right shift}%
\Index
   {right shift on group}%
   {right shift, group}%
\Index
   {right structural constant of Lie algebra}%
   {right structural constant of Lie algebra}%
\Index
   {right vector space}%
   {right vector space}%
\Index
   {right-ordered cycle notation of permutation}%
   {right-ordered cycle notation of permutation}%
\Index
   {right-side contravariant representation of group}%
   {right-side contravariant representation of group}%
\Index
   {right-side covariant representation of group}%
   {right-side covariant representation of group}%
\Index
   {right-side representation of fibered $\Omega$\Hyph algebra}%
   {right-side representation of fibered Omega-algebra}%
\Index
   {right-side representation of $\Omega_1$\Hyph algebra $A$ in $\Omega_2$\Hyph algebra $M$}%
   {right-side representation of algebra}%
\Index
   {right-side transformation}%
   {right-side transformation}%
\Index
   {ring has characteristic $0$}%
   {ring has characteristic 0}%
\Index
   {ring has characteristic $p$}%
   {ring has characteristic p}%
\Index
   {ring with conjugation}%
   {ring with conjugation}%
\Index
   {root of polynomial}%
   {root of polynomial}%
\Index
   {row determinant}%
   {row determinant}%
\Index
   {row vector}%
   {row vector}%
\SetIndexSpace%
\Index
   {$\star A$\Hyph module}%
   {starA-module}%
\Index
   {$D\star$\Hyph product of \rcd linear map $A$ over scalar}%
   {Dstar product of rcd linear map over scalar}%
\Index
   {$(\RCstar S,\RCstar T)$\Hyph linear map of vector spaces}%
   {rcs rct linear map of vector spaces}%
\Index
   {scalar algebra of algebra}%
   {scalar algebra of algebra}%
\Index
   {scalar algebra of ring}%
   {scalar algebra of ring}%
\Index
   {scalar of element of algebra}%
   {scalar of algebra}%
\Index
   {scalar of element of ring}%
   {scalar of ring}%
\Index
   {scalar of mapping}%
   {scalar of mapping}%
\Index
   {scalar potential}%
   {scalar potential}%
\Index
   {Schauder basis}%
   {Schauder basis}%
\Index
   {second Newton law}%
   {Second Newton law}%
\Index
   {section of bundle}%
   {section of bundle}%
\Index
   {set of coordinates of representation}%
   {coordinate set of representation}%
\Index
   {set of invertible elements of algebra}%
   {set of invertible elements of algebra}%
\Index
   {set of $\Omega_2$\Hyph words of representation}%
   {word set of representation}%
\Index
   {set of tuples of coordinates of tower of representations}%
   {coordinate set of tower of representations}%
\Index
   {set of tuples of $\Vector\Omega$\Hyph words of tower of representations}%
   {word set of tower of representations}%
\Index
   {set of zeros of algebra}%
   {set of zeros of algebra}%
\Index
   {simple polyvector}%
   {simple polyvector}%
\Index
   {single transitive representation of fibered $\Omega$\Hyph algebra}%
   {single transitive representation of fibered Omega-algebra}%
\Index
   {single transitive representation of group}%
   {single transitive representation of group}%
\Index
   {single transitive representation of $\Omega$\Hyph algebra $A$}%
   {single transitive representation of algebra}%
\Index
   {singular linear map}%
   {singular linear map}%
\Index
   {skew product of vectors}%
   {skew product of vectors}%
\Index
   {space of orbits of \Ts representation}%
   {space of orbits of Ts representation}%
\Index
   {spacelike vector}%
   {spacelike vector}%
\Index
   {speed of deviation}%
   {speed of deviation}%
\Index
   {$(\mathcal S\RCstar,\mathcal T\RCstar)$\Hyph linear map of vector bundles}%
   {src trc linear map of vector bundles}%
\Index
   {($S\star$, $\star T$)\hyph bimodule}%
   {(Sstar,starT)-bimodule}%
\Index
   {stability group}%
   {stability group}%
\Index
   {stable set of representation}%
   {stable set of representation}%
\Index
   {standard component of the G\^ateaux derivative}%
   {standard component of Gateaux derivative}%
\Index
   {standard component of linear map}%
   {standard component of linear map}%
\Index
   {standard component of polylinear map}%
   {standard component of polylinear map}%
\Index
   {standard component of quadratic map $f$ over field $F$}%
   {standard component of quadratic map, division ring}%
\Index
   {standard component of tensor}%
   {standard component of tensor}%
\Index
   {standard component over field $F$ of bilitnear map $f$}%
   {standard component of bilinear map, division ring}%
\Index
   {standard coordinates of basis}%
   {standard coordinates of basis}%
\Index
   {standard coordinates of \rcd basis}%
   {standard coordinates of rcd basis}%
\Index
   {standard representation of the G\^ateaux derivative}%
   {Gateaux derivative, standard representation}%
\Index
   {standard representation of linear map}%
   {linear map, standard representation}%
\Index
   {standard representation of matrix}%
   {Standard representation}%
\Index
   {standard representation of polylinear map}%
   {polylinear map, standard representation}%
\Index
   {standard representation of quadratic map of division ring over field $F$}%
   {quadratic map, standard representation, division ring}%
\Index
   {standard representation over field $F$ of bilinear map of division ring}%
   {bilinear map, standard representation, division ring}%
\Index
   {$\star A$\Hyph product of \Acr linear map $\Vector f$ over scalar}%
   {starA product of Acr linear map over scalar}%
\Index
   {$\star A$\Hyph product of $A\star$\Hyph linear map over scalar}%
   {starA product of Astar linear map over scalar}%
\Index
   {$\star D$\hyph vector space}%
   {starD-vector space}%
\Index
   {$\star R$\hyph module}%
   {starR-module}%
\Index
   {$\star A$\Hyph product of vector over scalar}%
   {starA product of vector over scalar, starA module}%
\Index
   {$\star D$\Hyph affine space}%
   {starD affine space}%
\Index
   {\sD component of coordinates of vector $\Vector r$}%
   {starD component of coordinates of vector, D vector space}%
\Index
   {$\star D$\Hyph vector space}%
   {starD vector space}%
\Index
   {$\star D$\hyph product of vector over scalar}%
   {starD product of vector over scalar, vector space}%
\Index
   {\sT representation of fibered group}%
   {starT representation of fibered group}%
\Index
   {\sT representation of fibered $\Omega$\Hyph algebra}%
   {starT representation of fibered Omega-algebra}%
\Index
   {\sT representation of $\Omega_1$\Hyph algebra $A$ in $\Omega_2$\Hyph algebra $M$}%
   {starT representation of algebra}%
\Index
   {\sT shift}%
   {starT shift}%
\Index
   {\sT shift on fibered group}%
   {starT shift, fibered group}%
\Index
   {\sT transformation}%
   {starT transformation}%
\Index
   {\sT transformation on bundle}%
   {starT transformation of bundle}%
\Index
   {structural constants}%
   {structural constants}%
\Index
   {subbundle}%
   {subbundle}%
\Index
   {subbundle of $\mathcal D\star$\hyph vector space}%
   {subbundle of Dstar vector bundle}%
\Index
   {subrepresentation generated by set $X$}%
   {subrepresentation generated by set}%
\Index
   {subrepresentation of representation}%
   {subrepresentation of representation}%
\Index
   {sum of geometric objects in vector space}%
   {sum of geometric objects, vector space}%
\Index
   {sum of geometric objects}%
   {sum of geometric objects}%
\Index
   {sum of linear maps}%
   {sum of linear maps}%
\Index
   {sum of \rcd linear maps}%
   {sum of rcd linear maps, rcd vector spaces}%
\Index
   {superposition of coordinates of the representation $f$ and the element $m$}%
   {superposition of coordinates, representation}%
\Index
   {superposition of coordinates of the tower of representations $\Vector f$ and the element $\VX a$}%
   {superposition of coordinates, tower of representations}%
\Index
   {symmetric $2$\Hyph ary fibered relation}%
   {symmetric 2 ary fibered relation}%
\Index
   {symmetric bilinear map of $D$\Hyph vector space to division ring}%
   {symmetric bilinear map, vector space to division ring}%
\Index
   {symmetric polylinear mapping into associative algebra}%
   {polylinear map symmetric, associative algebra}%
\Index
   {symmetry group}%
   {symmetry group}%
\Index
   {symmetry group}%
   {SymmetryGroup}%
\Index
   {synchronization of reference frame}%
   {synchronization of reference frame}%
\Index
   {system of additive equations}%
   {system of additive equations}%
\Index
   {system of \drc linear equations}%
   {system of drc linear equations}%
\Index
   {system of \rcd linear equations}%
   {system of rcd linear equations}%
\SetIndexSpace%
\Index
   {Taylor polynomial}%
   {Taylor polynomial, division ring}%
\Index
   {Taylor series}%
   {Taylor series, division ring}%
\Index
   {tensor power of algebra}%
   {tensor power of algebra}%
\Index
   {tensor power of division ring}%
   {tensor power of division ring}%
\Index
   {tensor power of representation}%
   {tensor power of representation}%
\Index
   {tensor product of algebras}%
   {tensor product of algebras}%
\Index
   {tensor product of $D$\Hyph vector spaces}%
   {tensor product of D vector spaces}%
\Index
   {tensor product of division rings}%
   {tensor product of division rings}%
\Index
   {tensor product of \Ds vector spaces}%
   {tensor product of Dstar vector spaces}%
\Index
   {tensor product of representations}%
   {tensor product of representations}%
\Index
   {tensor product of representations}%
   {tensor product of representations}%
\Index
   {tensor product of rings over commutative ring}%
   {tensor product of rings}%
\Index
   {tensor product of vector spaces}%
   {tensor product of vector spaces}%
\Index
   {the Fr\'echet \Ds derivative of map $f$ of division ring $D$ at point $x$}%
   {Frechet Dstar derivative of map, division ring}%
\Index
   {timelike vector}%
   {timelike vector}%
\Index
   {topological $D$\Hyph vector space}%
   {topological D vector space}%
\Index
   {topological $D$\Hyph algebra}%
   {topological D algebra}%
\Index
   {topological division ring}%
   {topological division ring}%
\Index
   {topological ring}%
   {topological ring}%
\Index
   {torsion form}%
   {torsion form}%
\Index
   {torsion tensor}%
   {torsion tensor}%
\Index
   {tower of bundles}%
   {tower of bundles}%
\Index
   {tower of effective representations}%
   {tower of effective representations}%
\Index
   {tower of representations of $\Vector{\Omega}$\Hyph algebras}%
   {tower of representations of algebras}%
\Index
   {tower of subrepresentations}%
   {tower of subrepresentations}%
\Index
   {tower of subrepresentations of tower of representations $\Vector f$ generated by tuple of sets $\VX X$}%
   {subrepresentation generated by tuple of sets}%
\Index
   {trace of quaternion}%
   {trace, quaternion algebra}%
\Index
   {transformation coordinated with equivalence}%
   {transformation coordinated with equivalence}%
\Index
   {transformation of universal algebra}%
   {transformation of universal algebra}%
\Index
   {transformation on bundle}%
   {transformation of bundle}%
\Index
   {transitive $2$\Hyph ary fibered relation}%
   {transitive 2 ary fibered relation}%
\Index
   {transitive representation of fibered $\Omega$\Hyph algebra}%
   {transitive representation of fibered Omega-algebra}%
\Index
   {transitive representation of group}%
   {transitive representation of group}%
\Index
   {transitive representation of $\Omega$\Hyph algebra $A$}%
   {transitive representation of algebra}%
\Index
   {\Ts matrices vector space}%
   {matrices vector space}%
\Index
   {\Ts representation of fibered $\Omega$\Hyph algebra}%
   {Tstar representation of fibered Omega-algebra}%
\Index
   {\Ts representation of $\Omega_1$\Hyph algebra $A$ in $\Omega_2$\Hyph algebra $M$}%
   {Tstar representation of algebra}%
\Index
   {\Ts shift}%
   {Tstar shift}%
\Index
   {\Ts transformation}%
   {Tstar transformation}%
\Index
   {\Ts transformation on bundle}%
   {Tstar transformation of bundle}%
\Index
   {tuple of coordinates of element $\Vector a$ relative to tuple of sets $\VX X$}%
   {coordinates of element, tower of representations}%
\Index
   {tuple of equivalence generated by tower of representations $\Vector f$}%
   {tuple of equivalence of tower of representations}%
\Index
   {tuple of generating sets of tower of representations}%
   {tuple of generating sets of tower of representations}%
\Index
   {tuple of generating sets of tower subrepresentations}%
   {tuple of generating sets of tower subrepresentations}%
\Index
   {tuple of $\Vector{\Omega}$\Hyph words of element of tower of representations relative to tuple of generating sets}%
   {tuple of words relative to tuple of generating sets, tower of representations}%
\Index
   {tuple of stable sets of tower of representation}%
   {tuple of stable sets of tower of representations}%
\Index
   {twin representations of associative algebra}%
   {twin representations of associative algebra}%
\Index
   {twin representations of $D$\Hyph algebra}%
   {twin representations of D algebra}%
\Index
   {twin representations of division ring}%
   {twin representations of division ring}%
\Index
   {twin representations of fibered group}%
   {twin representations of fibered group}%
\Index
   {twin representations of group}%
   {twin representations of group}%
\SetIndexSpace%
\Index
   {unit sphere in $D$\Hyph algebra}%
   {unit sphere in algebra}%
\Index
   {unit sphere in division ring}%
   {unit sphere in division ring}%
\Index
   {unit vector}%
   {unit vector}%
\Index
   {unitarity law for  $\mathcal D\star$\Hyph vector fields}%
   {unitarity law, Dstar vector fields}%
\Index
   {unitarity law for $A\star$\Hyph module}%
   {unitarity law, Astar module over algebra}%
\Index
   {unitarity law for $A\star$\Hyph vector space}%
   {unitarity law, Astar vector space}%
\Index
   {unitarity law for $D$\Hyph module}%
   {unitarity law, D module}%
\Index
   {unitarity law for $D\star$\Hyph vector space}%
   {unitarity law, Dstar vector space}%
\Index
   {unitarity law for $\star A$\Hyph module}%
   {unitarity law, starA module over algebra}%
\Index
   {unitarity law for $\star D$\Hyph vector space}%
   {unitarity law, starD vector space}%
\SetIndexSpace%
\Index
   {valued division ring}%
   {valued division ring}%
\Index
   {vector bundle}%
   {vector bundle}%
\Index
   {vector module of algebra}%
   {vector module of algebra}%
\Index
   {vector module of ring}%
   {vector module of ring}%
\Index
   {vector of element of algebra}%
   {vector of algebra}%
\Index
   {vector of element of ring}%
   {vector of ring}%
\Index
   {vector of mapping}%
   {vector of mapping}%
\Index
   {vector potential}%
   {vector potential}%
\Index
   {vector space over field}%
   {vector space over field}%
\Index
   {vector space type}%
   {vector space type}%

\CloseIndex

%% file: Symbol.English.tex
\def\indexname{Special Symbols and Notations}
\OpenIndex

\SetIndexSpace
\Symb
   {minor matrix}%
   {A from b a}%
\Symb
   {minor matrix}%
   {A from columns T}%
\Symb
   {minor matrix}%
   {A from rows S}%
\Symb
   {set of vectors whose expansion relative to the basis $\Basis e$ converges normally}%
   {A plus Schauder}%
\Symb
   {minor matrix}%
   {A without column a}%
\Symb
   {minor matrix}%
   {A without columns T}%
\Symb
   {minor matrix}%
   {A without row b}%
\Symb
   {minor matrix}%
   {A without rows S}%
\Symb
   {$A\CRcirc$\Hyph linear combination}%
   {ACRcirc linear combination 1}%
\Symb
   {$A\CRcirc$\Hyph linear combination}%
   {ACRcirc linear combination 2}%
\Symb
   {active representation of group $G(f)$ in basis manifold $\mathcal B(f)$}%
   {active representation in basis manifold}%
\Symb
   {active representation of group $G(\Vector f)$ in basis manifold $\mathcal B(\Vector f)$}%
   {active representation in basis manifold, tower of representations}%
\Symb
   {affine space}%
   {affine space, division ring}%
\Symb
   {algebra of polynomials over $D$\Hyph algebra $A$}%
   {algebra of polynomials over D algebra}%
\Symb
   {algebra of rational mappings of algebra $A$}%
   {algebra of rational mappings of algebra}%
\Symb
   {affine space}%
   {An}%
\Symb
   {associator of $D$\Hyph algebra}%
   {associator of algebra}%
\Symb
   {\subs row ($c$\hyph row) of matrix}%
   {c row}%
\Symb
   {commutator of $D$\Hyph algebra}%
   {commutator of algebra}%
\Symb
   {component of linear map $A$ of $D$\Hyph vector space}%
   {component of linear map, D vector space}%
\Symb
   {component $p$ of polylinear mapping $\Vector A$}%
   {component of polyadditive map, D vector space}%
\Symb
   {conjugated $D$\Hyph  module}%
   {conjugated D module}%
\Symb
   {linear combination of vectors of $A$\Hyph module}%
   {CR linear combination in A module}%
\Symb
   {\CR power of element $A$ of biring}%
   {cr power}%
\Symb
   {\CR inverse element of biring}%
   {cr-inverse element}%
\Symb
   {\CR product of matrices}%
   {cr-product of matrices}%
\Symb
   {derivative of left shift}%
   {derivative of left shift}%
\Symb
   {derivative of left shift in $1$\Hyph parameter Lie group}%
   {derivative of left shift, 1-Parameter Group}%
\Symb
   {derivative of left shift in $1$\Hyph parameter Lie D group}%
   {derivative of left shift, 1-Parameter Group, algebra}%
\Symb
   {}%
   {derivative of right shift}%
\Symb
   {}%
   {derivative of right shift}%
\Symb
   {derivative of right shift in $1$\Hyph parameter Lie group}%
   {derivative of right shift, 1-Parameter Group}%
\Symb
   {derivative of right shift in $1$\Hyph parameter Lie D group}%
   {derivative of right shift, 1-Parameter Group, algebra}%
\Symb
   {derivative of left shift}%
   {derivative of Tstar shift}%
\Symb
   {\drc vector}%
   {drc vector}%
\Symb
   {coordinates of vector $a$ relative to Hamel basis}%
   {Hamel basis, coordinates}%
\Symb
   {hermitian conjugation in division ring}%
   {hermitian conjugation, division ring}%
\Symb
   {$(^j_i)$\hyph\CR quasideterminant}%
   {j i CR quasideterminant definition}%
\Symb
   {$(ji)$\hyph quasideterminant of matrix $\bfA$}%
   {j i quasideterminant definition}%
\Symb
   {$(^j_i)$\hyph \RC quasideterminant}%
   {j i RC-quasideterminant definition}%
\Symb
   {$(^j_i)$\hyph $\RCcirc$\Hyph quasideterminant}%
   {j i RCcirc-quasideterminant definition}%
\Symb
   {left fraction}%
   {left fraction}%
\Symb
   {left principal ideal}%
   {left principal ideal}%
\Symb
   {left shift in $D$\Hyph algebra}%
   {left shift, D algebra}%
\Symb
   {linear combination of vectors of $A$\Hyph module}%
   {linear combination in A module}%
\Symb
   {transformation of matrix}%
   {matrix, replacing its column}%
\Symb
   {transformation of matrix}%
   {matrix, replacing its row}%
\Symb
   {norm on $D$\Hyph module}%
   {norm on D module}%
\Symb
   {opposite algebra to algebra $A$}%
   {opposite algebra}%
\Symb
   {orbit of linear mapping}%
   {orbit of linear mapping}%
\Symb
   {derivative}%
   {overline nabla_l, definition 2}%
\Symb
   {partial linear map}%
   {partial linear map}%
\Symb
   {principal ideal}%
   {principal ideal}%
\Symb
   {quasideterminant of matrix $\bfA$}%
   {quasideterminant definition}%
\Symb
   {\sups row ($r$\hyph row) of matrix}%
   {r row}%
\Symb
   {\RC power of element $A$ of biring}%
   {rc power}%
\Symb
   {\RC inverse element of biring}%
   {rc-inverse element}%
\Symb
   {\RC product of matrices}%
   {rc-product of matrices}%
\Symb
   {\RC quasideterminant}%
   {RC-quasideterminant definition}%
\Symb
   {$\RCcirc$\Hyph quasideterminant}%
   {RCcirc-quasideterminant definition}%
\Symb
   {\rcd vector}%
   {rcd vector}%
\Symb
   {right principal ideal}%
   {right principal ideal}%
\Symb
   {right shift in $D$\Hyph algebra}%
   {right shift, D algebra}%
\Symb
   {coordinates of vector $a$ relative to Schauder basis}%
   {Schauder basis, coordinates}%
\Symb
   {set of invertible elements of algebra $A$}%
   {set of invertible elements of algebra}%
\Symb
   {set of zeros of algebra $A$}%
   {set of zeros of algebra}%
\Symb
   {set of polylinear maps of rings $R_1$, ..., $R_n$ into module $S$}%
   {set polylinear maps, ring}%
\Symb
   {skew product of vectors $\Vector a_1$, ..., $\Vector a_m$}%
   {skew product of vectors}%
\Symb
   {standard component of tensor in tensor product of algebras}%
   {standard component of tensor, algebra}%
\Symb
   {right shift}%
   {starT shift}%
\Symb
   {\sT shift}%
   {starT shift, fibered group}%
\Symb
   {tensor power of algebra $A$}%
   {tensor power of algebra}%
\Symb
   {tensor product of algebras}%
   {tensor product of algebras}%
\Symb
   {left shift}%
   {Tstar shift}%
\Symb
   {\Ts shift}%
   {Tstar shift, fibered group}%
\Symb
   {anholonomic coordinates of vector}%
   {vector anholonomic coordinates}%
\Symb
   {holonomic coordinates of vector}%
   {vector holonomic coordinates}%

\SetIndexSpace
\Symb
   {basis manifold of \rcd vector space $\Vector V$}%
   {basis manifold of rcd vector space}%
\Symb
   {basis manifold of vector space}%
   {basis manifold of vector space}%
\Symb
   {basis manifold of representation $f$}%
   {basis manifold representation F algebra}%
\Symb
   {basis manifold of tower of representations $\Vector f$}%
   {basis manifold tower of representations}%
\Symb
   {basis manifold of affine space}%
   {Basis Manifold, Affine Space}%
\Symb
   {basis manifold of \rcd affine space}%
   {Basis Manifold, rcd Affine Space, division ring}%
\Symb
   {basis manifold of central affine space}%
   {BCAn}%
\Symb
   {basis manifold of Euclid space}%
   {BEn}%
\Symb
   {Cartesian power $\Bundle A$ of bundle $\Bundle B$}%
   {Cartesian power A of bundle B}%
\Symb
   {Cartesian power $A$ of set $B$}%
   {Cartesian power of set}%
\Symb
   {basis manifold of central affine space}%
   {FCAn}%
\Symb
   {basis manifold of Euclid space}%
   {FEn}%
\Symb
   {lattice of subrepresentations of representation $f$}%
   {lattice of subrepresentations}%
\Symb
   {lattice of towers of subrepresentations of tower of representations $\Vector f$}%
   {lattice of subrepresentations, tower of representations}%
\Symb
   {product of objects $B_1$, ..., $B_n$ in category $\mathcal A$}%
   {product of objects in category, 1 n}%
\Symb
   {right fraction}%
   {right fraction}%
\Symb
   {tensor power of representation}%
   {tensor power of representation}%
\Symb
   {tensor product of representations}%
   {tensor product of representations}%

\SetIndexSpace
\Symb
   {central affine space}%
   {CAn}%
\Symb
   {central affine space}%
   {central affine space}%
\Symb
   {$j$th column determinant of matrix $\bfA$}%
   {column determinant}%
\Symb
   {$\CRcirc$\Hyph product of matrices of maps}%
   {cr product of matrices of mappings}%
\Symb
   {left structural constant of Lie algebra}%
   {left structural constant of Lie algebra}%
\Symb
   {right structural constant of Lie algebra}%
   {right structural constant of Lie algebra}%
\Symb
   {set of continuous multivariable maps}%
   {set continuous multivariable maps}%
\Symb
   {structural constants}%
   {structural constants}%

\SetIndexSpace
\Symb
   {basis vector of representation of Lie group}%
   {basis vector of representation of Lie group}%
\Symb
   {basis vector of representation of Lie group over algebra $A$}%
   {basis vector of representation of Lie group over algebra A}%
\Symb
   {coordinates of basis vector of representation of Lie group over algebra $A$}%
   {basis vector of representation of Lie group over algebra A, coordinates}%
\Symb
   {coordinates of basis vector of representation of Lie group}%
   {basis vector of representation of Lie group, coordinates}%
\Symb
   {component of the G\^ateaux derivative of map $f(x)$}%
   {component of Gateaux derivative}%
\Symb
   {component of the G\^ateaux derivative of map $f(x)$}%
   {component of Gateaux derivative of map, D vector space, short}%
\Symb
   {component of the G\^ateaux derivative of second order of map $f(x)$}%
   {component of Gateaux derivative of Second Order}%
\Symb
   {component of the G\^ateaux derivative of second order of map $f(x)$}%
   {component of Gateaux derivative of Second Order, D vector space}%
\Symb
   {component of the G\^ateaux derivative of map $f(x)$}%
   {component of Gateaux derivative, vector space}%
\Symb
   {conjugation in algebra}%
   {conjugation in algebra}%
\Symb
   {conjugation in ring}%
   {conjugation in ring}%
\Symb
   {coordinate \Drc vector bundle}%
   {coordinate drc vector bundle}%
\Symb
   {coordinate \rcd vector space}%
   {coordinate rcd vector space}%
\Symb
   {coordinate reference frame}%
   {coordinate reference frame, extensive definition}%
\Symb
   {diagonal in bundle $\Bundle A$}%
   {diagonal in bundle, 1}%
\Symb
   {direct product of division rings $D_1$, ..., $D_n$}%
   {direct product of division rings, 1 n}%
\Symb
   {double determinant of matrix $\bfA$}%
   {double determinant}%
\Symb
   {the Fr\'echet \Ds derivative of map $f$ of division ring}%
   {Frechet Dstar derivative of map, division ring}%
\Symb
   {the G\^ateaux \crd derivative of map $\Vector f$ of $D$\hyph vector space $\Vector V$ to $D$\hyph vector space $\Vector W$}%
   {Gateaux crd derivative of map, D vector space}%
\Symb
   {the G\^ateaux derivative of map $f$}%
   {Gateaux derivative of map}%
\Symb
   {the G\^ateaux derivative of map $f$}%
   {Gateaux derivative of map, fraction}%
\Symb
   {the G\^ateaux derivative of order $n$}%
   {Gateaux derivative of Order n}%
\Symb
   {the G\^ateaux derivative of order $n$ of map $f$ of division ring}%
   {Gateaux derivative of Order n, division ring}%
\Symb
   {the G\^ateaux derivative of order $n$ of map $f$ of algebra}%
   {Gateaux derivative of Order n, fraction, algebra}%
\Symb
   {the G\^ateaux derivative of order $n$ of map $f$ of division ring}%
   {Gateaux derivative of Order n, fraction, division ring}%
\Symb
   {the G\^ateaux derivative of second order}%
   {Gateaux derivative of Second Order}%
\Symb
   {the G\^ateaux derivative of second order of mapping $f$ of algebra}%
   {Gateaux derivative of Second Order, fraction, algebra}%
\Symb
   {the G\^ateaux derivative of second order of map $f$ of division ring}%
   {Gateaux derivative of Second Order, fraction, division ring}%
\Symb
   {the G\^ateaux differential of map $f$}%
   {Gateaux differential of map, scalar}%
\Symb
   {the G\^ateaux differential of map $f$}%
   {Gateaux differential of map, vector}%
\Symb
   {the G\^ateaux \drc derivative of map $\Vector f$ of $D$\Hyph vector space $\Vector V$ to $D$\Hyph vector space $\Vector W$}%
   {Gateaux drc derivative of map, D vector space}%
\Symb
   {the G\^ateaux \Ds derivative of map $f$ of division ring $D$}%
   {Gateaux Dstar derivative of map, division ring}%
\Symb
   {the G\^ateaux Jacobian of map of $D$\Hyph vector space}%
   {Gateaux Jacobian of map, D vector space}%
\Symb
   {the G\^ateaux partial \drc derivative of map $f^b$ with respect to variable $v^a$}%
   {Gateaux partial crd derivative of map, 1, D vector space}%
\Symb
   {the G\^ateaux partial \drc derivative of map $f^b$ with respect to variable $v^a$}%
   {Gateaux partial crd derivative of map, 2, D vector space}%
\Symb
   {the G\^ateaux partial \drc derivative of map $f^b$ with respect to variable $v^a$}%
   {Gateaux partial crd derivative of map, 3, D vector space}%
\Symb
   {the G\^ateaux partial derivative}%
   {Gateaux partial derivative}%
\Symb
   {the G\^ateaux mixed partial derivative}%
   {Gateaux partial derivative of Second Order}%
\Symb
   {the G\^ateaux partial \drc derivative of map $f^b$ with respect to variable $v^a$}%
   {Gateaux partial drc derivative of map, 1, D vector space}%
\Symb
   {the G\^ateaux partial \drc derivative of map $f^b$ with respect to variable $v^a$}%
   {Gateaux partial drc derivative of map, 2, D vector space}%
\Symb
   {the G\^ateaux partial \drc derivative of map $f^b$ with respect to variable $v^a$}%
   {Gateaux partial drc derivative of map, 3, D vector space}%
\Symb
   {the G\^ateaux \sD derivative of map $f$ of division ring $D$}%
   {Gateaux starD derivative of map, division ring}%
\Symb
   {matrices vector space}%
   {matrices vector space}%
\Symb
   {Cartan derivative}%
   {overbrace D}%
\Symb
   {derivative}%
   {overline D}%
\Symb
   {derivative $e_{(k)}$}%
   {partial(k)}%
\Symb
   {product of map over scalar}%
   {product of map over scalar}%
\Symb
   {\subs rows \drc vector space}%
   {r rows drc vector space}%
\Symb
   {speed of deviation}%
   {speed of deviation}%
\Symb
   {standard component of the G\^ateaux derivative}%
   {standard component of Gateaux derivative}%
\Symb
   {tensor power of division ring $D$}%
   {tensor power of division ring}%
\Symb
   {tensor product of division rings}%
   {tensor product of division rings}%
\Symb
   {tensor product of rings}%
   {tensor product of rings}%
\Symb
   {vector space type}%
   {vector space type}%

\SetIndexSpace
\Symb
   {$A\CRcirc$\Hyph basis for module}%
   {A CRcirc basis, module over algebra}%
\Symb
   {Jacobian matrix of left shift}%
   {aE, quaternion, Jacobian matrix}%
\Symb
   {affine basis}%
   {Affine Basis}%
\Symb
   {basis of vector space}%
   {Basis e}%
\Symb
   {basis in vector space $\Vector V$}%
   {basis in V}%
\Symb
   {basis of $D$\Hyph module $\mathcal L(D;A_1;A_2)$}%
   {basis L(A1,A2)}%
\Symb
   {basis of vector space}%
   {basis, vector space}%
\Symb
   {basis of $(n)$\hyph vector space}%
   {basis,n vector space}%
\Symb
   {Cartesian power of total spaces}%
   {Cartesian power of total spaces}%
\Symb
   {Cartesian product of total spaces}%
   {Cartesian product of total spaces, definition 1}%
\Symb
   {central affine basis}%
   {Central Affine Basis}%
\Symb
   {basis for \Drc vector bundle}%
   {drc basis, vector bundle}%
\Symb
   {form of reference frame}%
   {dual forms, reference frame}%
\Symb
   {Euclid space}%
   {Euclid space}%
\Symb
   {Euclid space}%
   {Euclid space, division ring}%
\Symb
   {Hamel basis}%
   {Hamel basis}%
\Symb
   {identical transformation of bundle}%
   {identical transformation of bundle}%
\Symb
   {linear automorphism of quaternioin algebra}%
   {mapping E, quaternion}%
\Symb
   {linear automorphism of quaternioin algebra}%
   {mapping E_1, quaternion}%
\Symb
   {linear automorphism of quaternioin algebra}%
   {mapping E_2, quaternion}%
\Symb
   {orthonornal basis}%
   {Orthonornal Basis}%
\Symb
   {pseudo Euclid space}%
   {pseudo Euclid space}%
\Symb
   {pseudo Euclid space}%
   {pseudo Euclid space, division ring}%
\Symb
   {quaternion algebra over the field $F$}%
   {quaternion algebra over the field}%
\Symb
   {quaternion division algebra over the field}%
   {quaternion division algebra over the fieldL}%
\Symb
   {$\RCcirc A$\Hyph linear combination}%
   {RCcircA linear combination 1}%
\Symb
   {$\RCcirc A$\Hyph linear combination}%
   {RCcircA linear combination 2}%
\Symb
   {\rcd affine basis}%
   {rcd affine basis, division ring}%
\Symb
   {reduced Cartesian product of total spaces}%
   {reduced Cartesian product of total spaces, definition 1}%
\Symb
   {Schauder basis}%
   {Schauder basis}%
\Symb
   {set of nonsingular \sT transformations of bundle $\Bundle E$}%
   {set of starT nonsingular transformations of bundle}%
\Symb
   {set of nonsingular \Ts transformations of bundle $\Bundle E$}%
   {set of Tstar nonsingular transformations of bundle}%
\Symb
   {standard coordinates of basis}%
   {standard coordinates of basis}%
\Symb
   {standard coordinates of reference frame}%
   {standard coordinates of reference frame}%
\Symb
   {vector field of reference frame}%
   {vector field of reference frame}%
\Symb
   {vector of basis}%
   {vector of basis}%

\SetIndexSpace
\Symb
   {coordinates of basis in \sups rows \rcd vector space}%
   {basis coordinates, c rows rcd vector space}%
\Symb
   {coordinates of basis in \subs rows \drc vector space}%
   {basis coordinates, r rows drc vector space}%
\Symb
   {basis for \sups rows \rcd vector space}%
   {basis, c rows rcd vector space}%
\Symb
   {basis for \subs rows \drc vector space}%
   {basis, r rows drc vector space}%
\Symb
   {central affine basis}%
   {Central Affine Basis, division ring}%
\Symb
   {component of linear map $f$ of division ring}%
   {component of linear map, division ring}%
\Symb
   {component of polylinear map}%
   {component of polylinear map}%
\Symb
   {fibered morphism from bundle $\Bundle A$ into $\Bundle B$}%
   {fibered morphism from A into B}%
\Symb
   {filter $\mathfrak{F}$ converges to set $A$}%
   {filter converges}%
\Symb
   {homomorphism of fibered universal algebras}%
   {homomorphism of fibered universal algebras}%
\Symb
   {inverse fibered correspondence}%
   {inverse fibered correspondence, 1}%
\Symb
   {inverse reduced fibered correspondence}%
   {inverse reduced fibered correspondence, 1}%
\Symb
   {map to Cartesian product}%
   {map to Cartesian product}%
\Symb
   {norm of functional}%
   {norm of functional}%
\Symb
   {norm of map}%
   {norm of map}%
\Symb
   {norm of polylinear map}%
   {norm of polymap}%
\Symb
   {representation orbit of group $G$}%
   {orbit of Tstar representation of group}%
\Symb
   {orthonornal basis}%
   {Orthonornal Basis, division ring}%
\Symb
   {quaternion algebra  over field ${\rm {\mathbb{F}}}$}%
   {quaternion algebra F a b}%
\Symb
   {reference frame}%
   {reference frame}%
\Symb
   {reference frame, extensive definition}%
   {reference frame, extensive definition}%
\Symb
   {standard component of biadditive map $f$ over field $F$}%
   {standard component of biadditive map, division ring}%
\Symb
   {standard component of linear map}%
   {standard component of linear map, G}%
\Symb
   {standard component of polylinear map}%
   {standard component of polylinear map}%
\Symb
   {standard component of quadratic map $f$ over field $F$}%
   {standard component of quadratic map, division ring}%
\Symb
   {standard component of tensor}%
   {standard component of tensor}%
\Symb
   {sum of linear maps}%
   {sum of linear maps}%

\SetIndexSpace
\Symb
   {affine transformation group}%
   {affine transformation group}%
\Symb
   {\CR matrix group}%
   {cr-matrix group}%
\Symb
   {affine transformation group}%
   {drc affine transformation group}%
\Symb
   {fibered little group of section $h$}%
   {fibered little group}%
\Symb
   {fibered stability group of section $h$}%
   {fibered stability group}%
\Symb
   {algebra Lie of group Lie}%
   {g}%
\Symb
   {left defined algebra Lie of group Lie}%
   {gl}%
\Symb
   {right defined algebra Lie of group Lie}%
   {gr}%
\Symb
   {group of automorphisms of representation $f$}%
   {group of automorphisms of representation}%
\Symb
   {group of homomorphisms of vector space $\Vector V$}%
   {GV}%
\Symb
   {little group of $x$}%
   {little group}%
\Symb
   {orbit of effective covariant \sT representation of fibered group}%
   {orbit of effective covariant starT representation of fibered group}%
\Symb
   {orbit of effective covariant \sT representation of group}%
   {orbit of effective covariant starT representation of group}%
\Symb
   {orbit of effective covariant \Ts representation of fibered group}%
   {orbit of effective covariant Tstar representation of fibered group}%
\Symb
   {orbit of effective covariant \Ts representation of group}%
   {orbit of effective covariant Tstar representation of group}%
\Symb
   {product of groups $G_1$, ..., $G_n$}%
   {product of groups, 1 n}%
\Symb
   {\RC matrix group}%
   {rc-matrix group}%
\Symb
   {stability group of $x$}%
   {stability group}%

\SetIndexSpace
\Symb
   {Hadamard inverse of matrix}%
   {Hadamard inverse of matrix}%
\Symb
   {quaternion algebra}%
   {quaternion algebra H a b}%
\Symb
   {quaternion algebra over real field}%
   {quaternion algebra over real field}%

\SetIndexSpace
\Symb
   {infinitesimal generator of representation}%
   {infinitesimal generator I of representation}%
\Symb
   {infinitesimal generator of representation}%
   {infinitesimal generator i of representation}%
\Symb
   {Lie group infinitesimal generators}%
   {Lie group infinitesimal generators}%
\Symb
   {vector module of algebra $A$}%
   {vector module of algebra}%
\Symb
   {vector module of ring $D$}%
   {vector module of ring}%
\Symb
   {vector of element $d$ of algebra}%
   {vector of algebra}%
\Symb
   {vector of mapping $f$}%
   {vector of mapping}%
\Symb
   {vector of element $d$ of ring}%
   {vector of ring}%

\SetIndexSpace
\Symb
   {closure operator of representation $f$}%
   {closure operator, representation}%
\Symb
   {closure operator of tower of representations $\Vector f$}%
   {closure operator, tower of representations}%
\Symb
   {Jacobian matrix of right shift}%
   {Ea, quaternion, Jacobian matrix}%
\Symb
   {tower of subrepresentations of tower of representations $\Vector f$ generated by tuple of sets $\VX X$}%
   {subrepresentation generated by tuple of sets}%

\SetIndexSpace
\Symb
   {kernel of linear map}%
   {kernel of linear map}%

\SetIndexSpace
\Symb
   {left $ij$th cofactor of entry of matrix}%
   {left cofactor, matrix}%
\Symb
   {left double $ij$th cofactor of entry of matrix}%
   {left double cofactor}%
\Symb
   {left shift}%
   {left shift}%
\Symb
   {Lie derivative of connection}%
   {Lie derivative of connection}%
\Symb
   {Lie derivative of metric}%
   {Lie derivative of metric}%
\Symb
   {limit of correspondence $\Phi$ with respect to the filter $\mathfrak{F}$}%
   {limit of correspondence with respect to the filter}%
\Symb
   {limit of sequence}%
   {limit of sequence}%
\Symb
   {passive transformation}%
   {passive transformation}%
\Symb
   {set of \Acr linear maps of module $\Vector V$ into module $\Vector W$}%
   {set Acr linear maps, module}%
\Symb
   {$D$\Hyph module of continuous linear mappings of normed $D$\Hyph module $A_1$ into normed $D$\Hyph module $A_2$}%
   {set continuous linear mappings, module}%
\Symb
   {set of continuous linear maps}%
   {set continuous linear maps, vector}%
\Symb
   {set of continuous polylinear maps}%
   {set continuous polylinear maps}%
\Symb
   {\rcd vector space of linear maps of \drc vector space $V$ into \drc vector space $W$}%
   {set linear maps, drc vector space}%
\Symb
   {\drc vector space of linear maps of \rcd vector space $V$ into \rcd vector space $W$}%
   {set linear maps, rcd vector space}%
\Symb
   {set of linear maps}%
   {set linear maps, scalar}%
\Symb
   {set of linear maps}%
   {set linear maps, vector}%
\Symb
   {set of left-side nonsingular transformations of set $M$}%
   {set of left-side nonsingular transformations}%
\Symb
   {set of polylinear maps}%
   {set polylinear maps}%
\Symb
   {set of $n$\hyph linear maps}%
   {set polylinear maps An}%
\Symb
   {set of polylinear maps}%
   {set polylinear maps, D vector space}%
\Symb
   {set of polylinear maps of algebras $A_1$, ..., $A_n$ into algebra $A$}%
   {set polylinear maps, finite dimensional algebra}%
\Symb
   {set of \sT representations of division ring $S$ in additive group of division ring $R$}%
   {set sT representations, division ring}%
\Symb
   {set of \Ts representations of division ring $S$ in additive group of division ring $R$}%
   {set Ts representations, division ring}%

\SetIndexSpace
\Symb
   {set of \sT transformations of set $M$}%
   {set of starT transformations}%
\Symb
   {set of transformations of set $M$}%
   {set of transformations}%
\Symb
   {set of \Ts transformations of set $M$}%
   {set of Tstar transformations}%
\Symb
   {space of orbits of effective \sT covariant representation of the group}%
   {space of orbits of effective sT representation}%
\Symb
   {space of orbits of effective \Ts covariant representation of the group}%
   {space of orbits of effective Ts representation}%
\Symb
   {space of orbits of \Ts representation $f$ of group $G$ in set $M$}%
   {space of orbits of Ts representation}%

\SetIndexSpace
\Symb
   {norm of quaternion $x$}%
   {norm, quaternion algebra}%
\Symb
   {nucleus of $D$\Hyph algebra $A$}%
   {nucleus of algebra}%

\SetIndexSpace
\Symb
   {geometric object in coordinate representation defined in \rcd vector space}%
   {geometric object, coordinate rcd vector space}%
\Symb
   {geometric object in coordinate representation}%
   {geometric object, coordinate vector space}%
\Symb
   {geometric object defined in \rcd vector space}%
   {geometric object, rcd vector space}%
\Symb
   {octonion algebra}%
   {octonion algebra}%
\Symb
   {orbit of representation of fibered group $\Bundle G$}%
   {orbit of representation of fibered group}%
\Symb
   {orbit of representation of the group $G$}%
   {orbit of representation of group}%

\SetIndexSpace
\Symb
   {bundle}%
   {bundle}%
\Symb
   {bundle of level $2$}%
   {bundle of level 2}%
\Symb
   {bundle of level $n$}%
   {bundle of level n}%
\Symb
   {Cartesian power $n$ of bundle $\bundle{}{p}{E}{}$}%
   {Cartesian power of bundle}%
\Symb
   {Cartesian product of bundles}%
   {Cartesian product of bundles, definition 1}%
\Symb
   {passive representation of group $G(f)$ in basis manifold $\mathcal B(f)$}%
   {passive representation in basis manifold}%
\Symb
   {passive representation of group $G(\Vector f)$ in basis manifold $\mathcal B(\Vector f)$}%
   {passive representation in basis manifold, tower of representations}%
\Symb
   {reduced Cartesian product of bundles}%
   {reduced Cartesian product of bundles, definition 1}%
\Symb
   {set of nonsingular \sT transformations of bundle $\bundle{}pE{}$}%
   {set of starT nonsingular transformations of bundle, projection}%
\Symb
   {set of nonsingular \Ts transformations of bundle $\bundle{}pE{}$}%
   {set of Tstar nonsingular transformations of bundle, projection}%

\SetIndexSpace
\Symb
   {active transformation}%
   {active transformation}%
\Symb
   {\sups rows \rcd vector space}%
   {c rows rcd vector space}%
\Symb
   {Cartan curvature}%
   {Cartan curvature}%
\Symb
   {\CR rank of matrix}%
   {cr-rank of matrix}%
\Symb
   {diagonal in bundle  $\bundle{}pA{}$}%
   {diagonal in bundle, 2}%
\Symb
   {diagonal in bundle $\Bundle A$}%
   {diagonal in reduced bundle, 2}%
\Symb
   {\Ds component of coordinates of vector $\Vector r$}%
   {Dstar component of coordinates of vector, D vector space}%
\Symb
   {image of $m$ under endomorphism $R$ of effective representation}%
   {endomorphism image, effective representation}%
\Symb
   {image of tuple $\VX a$ under endomorphism $\VX r$ of tower of effective representations}%
   {endomorphism image, tower of effective representations}%
\Symb
   {curvature}%
   {GLn curvature_overline}%
\Symb
   {$\RCcirc$\Hyph product of matrices of maps}%
   {rc product of matrices of mappings}%
\Symb
   {\RC rank of matrix}%
   {rc-rank of matrix}%
\Symb
   {right $ij$th cofactor of entry of matrix}%
   {right cofactor, matrix}%
\Symb
   {right double $ij$th cofactor of entry of matrix}%
   {right double cofactor}%
\Symb
   {right shift}%
   {right shift}%
\Symb
   {$i$th row determinant of matrix $\bfA$}%
   {row determinant}%
\Symb
   {scalar algebra of algebra $A$}%
   {scalar algebra of algebra}%
\Symb
   {scalar algebra of ring $D$}%
   {scalar algebra of ring}%
\Symb
   {scalar of element $d$ of algebra}%
   {scalar of algebra}%
\Symb
   {scalar of mapping $f$}%
   {scalar of mapping}%
\Symb
   {scalar of element $d$ of ring}%
   {scalar of ring}%
\Symb
   {set of right-side nonsingular transformations of set $M$}%
   {set of right-side nonsingular transformations}%
\Symb
   {\sD component of coordinates of vector $\Vector r$}%
   {starD component of coordinates of vector, D vector space}%

\SetIndexSpace
\Symb
   {composition of fibered correspondences}%
   {composition of fibered correspondences}%
\Symb
   {inverse fibered correspondence}%
   {inverse fibered correspondence, 2}%
\Symb
   {inverse reduced fibered correspondence}%
   {inverse reduced fibered correspondence, 2}%
\Symb
   {linear span in vector space}%
   {linear span, vector space}%
\Symb
   {image of basis $X$ under passive transformation $S$}%
   {passive transformation of basis, representation}%
\Symb
   {image of basis $\VX  X$ under passive transformation $\VX s$}%
   {passive transformation of basis, tower of representations}%
\Symb
   {symmetric group}%
   {symmetric group}%

\SetIndexSpace
\Symb
   {category of \Ts representations of $\Omega_1$\Hyph algebra $A$}%
   {category of Tstar representations of Omega1 algebra}%
\Symb
   {category of \Ts representations of $\Omega_1$\Hyph algebra from category $\mathcal A$}%
   {category of Tstar representations of Omega1 algebra from category}%
\Symb
   {tangent plane to group $G$}%
   {TaG}%
\Symb
   {trace of quaternion $x$}%
   {trace, quaternion algebra}%

\SetIndexSpace
\Symb
   {coordinate vector space}%
   {coordinate vector space}%
\Symb
   {coordinates in vector space}%
   {coordinates in vector space}%
\Symb
   {direct product of $\RCstar D_i$\hyph vector spaces $\Vector V_1$, ..., $\Vector V_n$}%
   {direct product, rcd vector space, 1 n}%
\Symb
   {dual space of \rcd vector space $\Vector V$}%
   {dual space of rcd vector space}%
\Symb
   {hermitian conjugated vector}%
   {hermitian conjugated vector}%
\Symb
   {tensor product of $D$\Hyph vector spaces}%
   {tensor product of D vector spaces}%
\Symb
   {tensor product of \Ds vector spaces}%
   {tensor product of Dstar vector spaces}%
\Symb
   {vector space}%
   {V}%

\SetIndexSpace
\Symb
   {set of coordinates of representation $J(f,X)$}%
   {coordinate set of representation}%
\Symb
   {set of tuples of coordinates of tower of representations $\Vector J(\Vector f,\VX X)$}%
   {coordinate set of tower of representations}%
\Symb
   {coordinates of basis $X'$ relative to basis $X$ of representation}%
   {coordinates of basis relative to basis, representation}%
\Symb
   {coordinates of element $m$ of representation $f$ relative to set $X$}%
   {coordinates of element relative to generating set, representation}%
\Symb
   {coordinates of element $m$ relative to set $X$}%
   {coordinates of element relative to set, representation}%
\Symb
   {tuple of coordinates of element $\Vector a*$ relative to tuple of sets $\VX X$}%
   {coordinates of element, tower of representations}%
\Symb
   {geometric object in coordinate representation defined in $\Omega_2$\Hyph algebra $M$}%
   {geometric object, coordinate representation g}%
\Symb
   {geometric object in coordinate representation defined in tuple of $\VX\Omega$\Hyph algebras $\VX A$}%
   {geometric object, coordinate tower of representations g}%
\Symb
   {geometric object defined in $\Omega_2$\Hyph algebra $M$}%
   {geometric object, representation g}%
\Symb
   {geometric object defined in tuple of $\VX\Omega$\Hyph algebras $\VX A$}%
   {geometric object, tower of representations g}%
\Symb
   {geometric object in vector space}%
   {geometric object, vector space}%
\Symb
   {set of coordinates of set $B\subset J(f,X)$}%
   {subset of coordinates of representation}%
\Symb
   {coordinates of tuple of sets $\VX B$ relative to tuple of sets $\VX X$}%
   {subset of coordinates of tower of representations}%
\Symb
   {coordinates of set $B_k$ relative to tuple of sets $\VX X$}%
   {subset of coordinates of tower of representations, k}%
\Symb
   {set of $\Omega_2$\Hyph words representing set $B\subset J(f,X)$}%
   {subset of words of representation}%
\Symb
   {superposition of coordinates of the representation $f$ and the element $m$}%
   {superposition of coordinates, representation}%
\Symb
   {superposition of coordinates of the tower of representations $\Vector f$ and the element $\VX a$}%
   {superposition of coordinates, tower of representations}%
\Symb
   {$\Omega_2$\Hyph word representing element $m\in J(f,X)$}%
   {word of element relative to generating set, representation}%
\Symb
   {set of $\Omega_2$\Hyph words of representation $J(f,X)$}%
   {word set of representation}%
\Symb
   {set of tuples of $\VX{\Omega}$\Hyph words of tower of representations $\Vector J(\Vector f,\VX X)$}%
   {word set of tower of representations}%
\Symb
   {tuple of words of element $\Vector a*$ relative to tuple of sets $\VX X$}%
   {words of element, tower of representations}%

\SetIndexSpace
\Symb
   {conjugate of quaternion $x$}%
   {conjugate of quaternion}%
\Symb
   {local basis of affine space}%
   {local basis of affine space}%
\Symb
   {anholonomic coordinate}%
   {x(k)}%

\SetIndexSpace
\Symb
   {center of $D$\Hyph algebra $A$}%
   {center of algebra}%
\Symb
   {center of ring $D$}%
   {center of ring}%

\SetIndexSpace
\Symb
   {deviation of trajectories}%
   {deviation of trajectories}%
\Symb
   {identical transformation}%
   {identical transformation}%
\Symb
   {image of vector $\Vector e_k\in\Basis e$ under isomorphism to coordinate vector space}%
   {image of vector e_k, coordinate vector space}%
\Symb
   {Kronecker symbol}%
   {Kronecker symbol}%

\SetIndexSpace
\Symb
   {anholonomic coordinates of connection}%
   {anholonomic coordinates of connection}%
\Symb
   {Cartan symbol}%
   {Cartan symbol}%
\Symb
   {connection}%
   {conection overline}%
\Symb
   {connection coefficients in $D$\Hyph affine space}%
   {connection coefficients, D affine space}%
\Symb
   {connection in $D$\Hyph affine manifold}%
   {connection, affine manifold}%
\Symb
   {$D$\Hyph affine connection coefficients on manifold}%
   {D affine connection coefficients, manifold}%
\Symb
   {holonomic coordinates of connection}%
   {holonomic coordinates of connection}%
\Symb
   {Cartan connection}%
   {overbrace Gamma i kl}%
\Symb
   {set of sections of bundle}%
   {set of sections of bundle}%

\SetIndexSpace
\Symb
   {inverse operator to operator $\psi_l$}%
   {inverse operator to operator psi l}%
\Symb
   {inverse operator to operator $\psi_r$}%
   {inverse operator to operator psi r}%

\SetIndexSpace
\Symb
   {anholonomity object}%
   {anholonomity object}%

\SetIndexSpace
\Symb
   {left basic operator of Lie group over algebra $A$}%
   {L basic operator of Lie group over algebra A}%
\Symb
   {left basic operator of group Lie}%
   {Lie Basic Operator L}%
\Symb
   {left basic operator of Lie 1-parameter group}%
   {Lie Basic Operator L, 1-Parameter Group}%
\Symb
   {left basic operator of Lie 1-parameter group over algebra $A$}%
   {Lie Basic Operator L, 1-Parameter Group, algebra}%
\Symb
   {right basic operator of group Lie}%
   {Lie Basic Operator R}%
\Symb
   {right basic operator of Lie 1-parameter group}%
   {Lie Basic Operator R, 1-Parameter Group}%
\Symb
   {right basic operator of Lie 1-parameter group over algebra $A$}%
   {Lie Basic Operator R, 1-Parameter Group, algebra}%
\Symb
   {right basic operator of Lie group over algebra $A$}%
   {R basic operator of Lie group over algebra A}%

\SetIndexSpace
\Symb
   {Lie group composition law}%
   {Lie group composition law}%

\SetIndexSpace
\Symb
   {Cartan derivative}%
   {overbrace nabla_l}%
\Symb
   {derivative}%
   {overline nabla_l, definition 1}%

\SetIndexSpace
\Symb
   {restriction of correspondence $\Phi$ to set $C$}%
   {restriction of correspondence}%

\SetIndexSpace
\Symb
   {Cartesian product of bundles}%
   {Cartesian product of bundles, definition 2}%
\Symb
   {Cartesian product of total spaces}%
   {Cartesian product of total spaces, definition 2}%
\Symb
   {direct product of division rings $D_i$, $i\in I$}%
   {direct product of division rings}%
\Symb
   {direct product of division rings $D_1$, ..., $D_n$}%
   {direct product of division rings, i 1 n}%
\Symb
   {direct product of $\RCstar D_i$\hyph vector spaces $\Vector V_i$, $i\in I$}%
   {direct product, rcd vector space}%
\Symb
   {direct product of $\RCstar D_i$\hyph vector spaces}%
   {direct product, rcd vector space, i 1 n}%
\Symb
   {product of groups $G_i$, $i\in I$}%
   {product of groups}%
\Symb
   {product of groups $G_1$, ..., $G_n$}%
   {product of groups, i 1 n}%
\Symb
   {product of objects $\{B_i,i\in I\}$ in category $\mathcal A$}%
   {product of objects in category}%
\Symb
   {product of objects $B_1$, ..., $B_n$ in category $\mathcal A$}%
   {product of objects in category, i 1 n}%
\Symb
   {reduced Cartesian product of bundles}%
   {reduced Cartesian product of bundles, definition 2}%
\Symb
   {reduced Cartesian product of total spaces}%
   {reduced Cartesian product of total spaces, definition 2}%

\SetIndexSpace
\Symb
   {fibered subset}%
   {fibered subset}%
\Symb
   {subbundle}%
   {subbundle}%

\CloseIndex